\documentclass[article]{article}

\usepackage{amssymb}
\usepackage{amsmath}
\usepackage{amsfonts}
\usepackage{amsthm}
\usepackage{stmaryrd}
\usepackage[all]{xy}
\usepackage{mathrsfs}
\usepackage{graphicx}
\usepackage{hyperref}
\usepackage{color}
\usepackage{multirow}

\usepackage{comment}
\usepackage[shortlabels]{enumitem}

\usepackage{scalerel}

\usepackage{etoolbox}
\patchcmd{\thebibliography}{\section*}{\section*}{}{}

\usepackage{stackengine}

\usepackage{bbm}
\DeclareMathAlphabet{\mathpzc}{OT1}{pzc}{m}{it}

\numberwithin{equation}{subsection}
\newtheorem{theorem}[subsubsection]{Theorem}

\newtheorem{conjecture}[subsubsection]{Conjecture}
\newtheorem{corollary}[subsubsection]{Corollary}
\newtheorem{lemma}[subsubsection]{Lemma}
\newtheorem{proposition}[subsubsection]{Proposition}
\newtheorem{definition}[subsubsection]{Definition}

\theoremstyle{definition}

\newtheorem{claim}[subsubsection]{Claim}

\newtheorem{example}[subsubsection]{Example}

\newtheorem{fact}[subsubsection]{Fact}

\newtheorem{remark}[subsubsection]{Remark}

\newcommand{\ra}{\rightarrow}

\def\BB{\mathbb{B}}
\def\CC{\mathbb{C}}

\def\FF{\mathbb{F}}

\def\LL{\mathbb{L}}

\def\NN{\mathbb{N}}

\def\QQ{\mathbb{Q}}

\def\TT{\mathbb{T}}

\def\ZZ{\mathbb{Z}}

\def\bfs{\mathbf{s}}

\def\bfS{\mathbf{S}}

\def\calC{\mathcal{C}}

\def\calF{\mathcal{F}}
\def\calG{\mathcal{G}}
\def\calH{\mathcal{H}}
\def\calI{\mathcal{I}}

\def\calL{\mathcal{L}}

\def\calO{\mathcal{O}}

\def\calS{\mathcal{S}}

\def\calX{\mathcal{X}}

\def\scrC{\mathscr{C}}

\def\rmH{\mathrm{H}}

\def\ul{\underline}
\def\wt{\widetilde}
\def\wh{\widehat}

\DeclareMathOperator{\Aut}{Aut}

\DeclareMathOperator{\Gal}{Gal}

\DeclareMathOperator{\Hom}{Hom}
\DeclareMathOperator{\Ext}{Ext}

\DeclareMathOperator{\rank}{rank}

\DeclareMathOperator{\Spec}{Spec}
\DeclareMathOperator{\Spf}{Spf}

\DeclareMathOperator{\Proj}{Proj}

\newcommand{\A}{\mathbb{A}}
\newcommand{\Q}{\mathbb{Q}}
\newcommand{\Z}{\mathbb{Z}}

\newcommand{\C}{\mathbb{C}}

\newcommand{\Sh}{\mathrm{Sh}}

\newcommand{\rra}{\longrightarrow}
\newcommand{\lmt}{\longmapsto}

\newcommand{\cH}{\mathrm{H}}

\newcommand{\frakm}{\mathfrak{m}}

\newcommand{\Set}{\textbf{Set}}

\newcommand{\Sing}{\textrm{Sing}}

\newcommand{\red}{{\mathrm{red}}}

\newcommand{\pO}{\wh{\calO}}
\newcommand{\et}{{\mathrm{\acute{e}t}}}
\newcommand{\pe}{{\mathrm{pro\acute{e}t}}}

\newcommand{\Spa}{\mathrm{Spa}}
\newcommand{\RK}{\mathrm{Rig}_K}
\newcommand{\Ret}{\mathrm{Rig}_{K,\acute{e}t}}
\newcommand{\Reh}{\mathrm{Rig}_{K,\acute{e}h}}
\newcommand{\Bl}{\mathrm{Bl}}

\newcommand{\Pd}{\mathrm{Perfd}}

\newcommand{\Pf}{\mathrm{Perf}}
\newcommand{\Spd}{\mathrm{Spd}}
\newcommand{\eh}{{\mathrm{\acute{e}h}}}
\newcommand{\ct}{{\mathrm{cont}}}
\newcommand{\Ainf}{{\mathrm{A_{inf}}}}

\newcommand{\Bdr}{{\mathrm{B_{dR}^+}}}
\newcommand{\ind}{{\mathrm{ind}}}
\newcommand{\an}{{\mathrm{an}}}
\newcommand{\colim}{{\mathrm{colim}}}

\newcommand{\cosk}{{\mathrm{cosk}}}
\newcommand{\Ab}{{\mathrm{Ab}}}

\newcommand{\PE}{{\mathrm{pro\acute{e}t}}}
\newcommand{\sgn}{{\mathrm{sgn}}}
\newcommand{\pfd}{{\mathrm{perfd}}}
\newcommand{\dR}{{\mathrm{dR}}}
\newcommand{\Fil}{{\mathrm{Fil}}}
\newcommand{\OBdr}{{\calO\BB_\dR}}
\newcommand{\OBdrp}{{\calO\BB^+_\dR}}

\newcommand{\iif}{{\mathrm{inf}}}

\newcommand{\gr}{{\mathrm{gr}}}

\begin{document}
\date{}
  \title{Hodge-Tate decomposition for non-smooth spaces}\maketitle

  \centerline{Haoyang Guo}

\begin{abstract}
	In this article, we generalize the Hodge-Tate decomposition of $p$-adic \'etale cohomology to non-smooth rigid spaces.
	Our strategy is to study pro-\'etale cohomology of rigid spaces introduced by Scholze, using the resolution of singularities and the simplicial method.

\end{abstract}
  
  \tableofcontents
  
\section{Introduction}

\subsection{Goals and main results}
Let $X$ be a compact complex manifold.
One of the most important invariants of $X$ is its singular cohomology group $\cH^n_\Sing(X,\C)$, defined in a transcendental way.
The de Rham Theorem allows us to compute those cohomology groups by differential forms, which states that there exists a natural isomorphism
\[
\cH^n_\Sing(X,\C)\cong\cH^n(X,\Omega_{X/\C}^\bullet),
\]
where $\Omega_{X/\C}^\bullet$ is the analytic de Rham complex of $X$ over $\C$.
Moreover, induced by the Hodge filtration of the de Rham complex, there exists an (Hodge-de Rham) spectral sequence
\[
E_1^{i,j}=\cH^j(X,\Omega^i_{X/\C})\Longrightarrow\cH^{i+j}(X,\Omega_{X/\C}^\bullet)
\]
computing the de Rham cohomology.
This spectral sequence often degenerates;
in fact we have a stronger result when $X$ is a compact K\"ahler manifold.

\begin{theorem}[Hodge decomposition]\label{Hodge dec}
	Let $X$ be a compact K\"ahler manifold. 
	Then its singular cohomology $\rmH^n_\Sing(X,\mathbb{C})$ admits a canonical decomposition into a direct sum of the vector spaces $\rmH^{i,j}(X)$ of harmonic $(i,j)$-forms, for $i+j=n$.
\end{theorem}
The vector space $\rmH^{i,j}(X)$ is of the same dimension as $\rmH^j(X,\Omega_{X/\C}^i)$.
As a consequence, the Hodge-de Rham spectral sequence for $X$
\[
E_1^{i,j}=\cH^j(X,\Omega^i_{X/\C})\Longrightarrow \cH^{i+j}(X,\Omega_{X/\C}^\bullet)
\]
degenerates at its first page.

The above assumption holds in particular when $X$ is the analytification of a smooth projective algebraic variety $Y$ over $\C$.
Moreover, in the aforementioned setting, by the GAGA principle, the above spectral sequence can be replaced by the algebraic Hodge-de Rham spectral sequence using the  algebraic differential forms of $Y$ over $\C$.
The singular cohomology can thus be computed using a purely algebraic method.

Now we turn to the non-archimedean geometry.
Let us fix a complete and algebraic closed $p$-adic field extension $K/\Q_p$.

The proper smooth rigid spaces, introduced by Tate in the 1960s, form a natural analogue of compact complex manifolds in this setting.
Here typical examples include analytifications of proper smooth algebraic varieties over $K$.
Since the non-archimedean field $K$ is totally disconnected, the singular cohomology is not a meaningful invariant of $X$.
Instead, the correct analogue should be the $p$-adic \'etale cohomology groups $\cH^i(X_\et,\Z_p)$, which are defined by associating an \'etale site $X_\et$ to the rigid space $X$, and then taking the inverse limit of the \'etale cohomology group $\cH^i(X_\et,\Z/p^n)$.
Moreover, analogous to the complex setting, we have the following result computing \'etale cohomology:

\begin{theorem}[Hodge-Tate filtration]\label{HT dec}
	Let $X$ be a smooth proper rigid space over $K$.
	Then we have a natural spectral sequence
	\[
	E_2^{i,j}=\cH^i(X,\Omega^j_{X/K})(-j)\Longrightarrow \cH^{i+j}(X_\et,\Z_p)\otimes_{\Z_p} K.
	\]
	The spectral sequence degenerates at its $E_2$-page.
\end{theorem}
Here pointed out by a referee, it worths to mention that different from Theorem \ref{Hodge dec}, we do not need any ``K\"ahler condition" in the $p$-adic setting.

When $X$ is an algebraic scheme, the result is proved in many cases by Faltings.
Later on Scholze constructs the spectral sequence above for rigid spaces, and proves the degeneracy for those $X$ defined over a discretely valued field. 
Here we note that in the latter setting when $X$ comes from a discretely valued subfield (with perfect residue field), the filtration above is Galois equivariant and admits a canonical splitting, which is called the \emph{Hodge--Tate decomposition}.
The case for general proper smooth rigid spaces that are not necessarily defined over a discrete subfield is handled in \cite{BMS} by using a spreading-out technique of Conrad-Gabber. 
Similar to the complex geometry, we note that one of the biggest advantages of this decomposition is that the left side is given by a coherent cohomology, which is more of an algebraic nature than the $p$-adic \'etale cohomology.

The goal of our article is to generalize the Hodge-Tate decomposition to the case when $X$ is not necessarily smooth. 
In scheme theory, one strategy to generalize from the smooth setting to the non-smooth setting is to use the $h$-topology,  introduced by Voevodesky \cite{Voe96}.
The $h$-topology is defined by refining the \'{e}tale topology (on the category of all finite type $K$-schemes over $X$) in a way that we allow not just \'{e}tale coverings, but also proper surjective maps.
In particular, thanks to the resolution of singularities in characteristic $0$, any h-covering can be refined by smooth schemes, which makes the $h$-topology ``locally smooth".
Later Geisser introduced the $\eh$-topology (cf. \cite{Gei06}), which is a variant of Veovodesky's theory and is defined by a smaller but more amenable collection of covering families.

Inspired by those ideas, we introduce the $\eh$-topology $X_\eh$ for a rigid space $X$, where coverings are generated by \'etale coverings, universal homeomorphisms, and coverings associated to blowups (see Section \ref{sec2}).
Similarly to the scheme theory, our $\eh$-topology is locally smooth, and there exists a natural morphism of sites $\pi_X:X_\eh\ra X_\et$.
In addition, by sheafifying the usual sheaf of continuous differential $\Omega^j_{/K}$ in the $\eh$-topology, we obtain the sheaf of $\eh$-differential $\Omega^j_{\eh,/K}$ on the $\eh$ site $X_\eh$ over $K$. (We will use $\Omega_{\eh}^j$ to abbreviate the notation, when the base is clear.)

Our main theorems are the following:
\begin{theorem}[Hodge-Tate decomposition]\label{HTss}
	Let $X$ be a proper rigid space over a complete algebraically closed non-archimedean field $K/\Q_p$ of characteristic $0$.
	There then exists a natural spectral sequence
	\[
	E_2^{i,j}=\cH^i(X_\eh,\Omega_{ \eh}^j)(-j)\Longrightarrow \cH^{i+j}(X_\et,\Q_p)\otimes_{\Q_p}K,
	\]
	where $\cH^i(X_\eh,\Omega_{ \eh}^j)(-j)$ is the $i$-th $\eh$ cohomology group of $\Omega_{ \eh}^j$, and is equipped with a Galois action by the Tate twist of weight $j$ when $X$ is defined over a discretely valued subfield.
	The spectral sequence satisfies the following:
	\begin{enumerate}[(i)]
		\item The cohomology group $\cH^i(X_\eh, \Omega_{ \eh}^j)(-j)$ is a finite dimensional $K$-vector space that vanishes unless $0\leq i,j\leq n$.
		\item The spectral sequence degenerates at the $E_2$-page.
		\item When $X$ is a smooth rigid space, $\cH^i(X_\eh, \Omega_{ \eh}^j)(-j)$ is isomorphic to $\cH^i(X, \Omega_{X/K}^j)(-j)$, and the spectral sequence is the same as the Hodge-Tate spectral sequence for smooth proper rigid spaces in Theorem \ref{HT dec}.
	\end{enumerate}
\end{theorem}

The $\eh$ cohomology group $\cH^i(X_\eh,\Omega_\eh^j)$ above is not an exotic construction; in fact this is an analogue of the classical Deligne--Du Bois cohomology in the rigid spaces setting, and can be computed by cohomologies of coherent sheaves over the rigid space.
To see this, we notice that the map of sites $\pi_X:X_\eh\ra X_\et$ provides us with a natural quasi-isomorphism
\[
R\Gamma(X_\et,R\pi_{X *}\Omega^j_\eh)\cong R\Gamma(X_\eh,\Omega_\eh^j).
\]
As proved in the Proposition \ref{coh}, each higher pushforward $R^s\pi_{X *} \Omega^j_\eh$ is a coherent sheaf over $X$.
So the associated Leray spectral sequence above is computed by rigid cohomologies of coherent sheaves.

Moreover, $\eh$ cohomology for a proper rigid space $X$ coincides with its algebraic version when $X$ comes from a proper algebraic variety.
Precisely, assume $X=Y^\an$ is the analytification of a proper variety $Y$ over $K$,  and the base field $K$ is isomorphic to $\CC$ abstractly.
We can then get a functorial isomorphism
\[
\cH^i(X_\eh,\Omega_\eh^j)\cong \cH^i(Y,\ul\Omega_Y^j),
\]
where $\ul\Omega^j_Y$ is the $j$-th graded piece of the Deligne-Du Bois complex $\ul\Omega^\bullet_Y$ for the complex algebraic variety $Y$ (\cite[Section 7.3]{PS08}).
Furthermore, the cohomology group $\cH^i(Y,\ul\Omega_Y^j)$ is isomorphic to the $j$-th graded factor for the Hodge filtration of the singular cohomology $\cH^i_\Sing(Y(\CC), \CC)$. 
In fact, there exists a filtered isomorphism between the $\eh$ cohomology of $X$ and the singular cohomology of $Y$, so that the above isomorphism is obtained by taking the $j$-th graded factor (see Subsection \ref{subsec ag} for details).

Here we also want to mention that the analogous h-cohomology for algebraic varieties has already appeared in the $p$-adic Hodge theory, back to Beilinson's work \cite{Bei12}.
In loc. cit., Beilinson gives a different proof of the Hodge-Tate decomposition for algebraic varieties using de Jong's alterations.
In the complex settings, the h-cohomology (and the $\eh$-cohomology) of differentials appears for example in \cite{Gei06}, \cite{Lee07},  and \cite{HJ14}, where the relation between h-cohomology and different types of singularities in complex algebraic geometry is studied.

At last, we comment that in Theorem \ref{HTss}, the rigid space $X$ may not be defined over a discretely valued subfield (in which case the Tate twist appearing there will not enter the picture).
However when $X=Y_K$ is defined over a discretely valued subfield, the decomposition theorem above can be obtained from a comparison between the \'etale cohomology and an $\eh$ version de Rham cohomology, generalizing the \'etale-de Rham comparison for the smooth case in \cite{Sch13}.
Precisely, let $\mathrm{B_{dR}}$ be Fontaine's de Rham period ring, which is equipped with a natural filtration (see the discussion in Paragraph \ref{para Bdr}).
We then have the following result:
\begin{theorem}[\'etale-$\eh$ de Rham comparison]
	Let $Y$ be a proper rigid space over a discretely valued subfield $K_0$ of $K/\Q_p$ that has a perfect residue field.
	Then there exists a $\Gal(K/K_0)$-equivariant filtered isomorphism 
	\[
	\cH^n(Y_{K\,\et},\Q_p)\otimes_{\Q_p} \mathrm{B_{dR}}=\cH^n(Y_\eh,\Omega_{ \eh, /K_0}^\bullet)\otimes_{K_0}\mathrm{B_{dR}},
	\]
	whose $0$-th graded piece is:
	\[
	\cH^n(Y_{K\,\et},\Q_p)\otimes_{\Q_p} K=\bigoplus_{i+j=n}\cH^i(Y_\eh,\Omega_{ \eh,/K_0}^j)\otimes_{K_0} K(-j).
	\]
	The isomorphism is functorial with respect to $Y$ over $K_0$.
	In particular, the $p$-adic Galois representation $\cH^n(Y_{K\,\et},\Q_p)$ is de Rham.
\end{theorem}
Here recall that the filtration on the left side above is the tensor product filtration for the trivial filtration of the \'etale cohomology and the natural filtration on $\mathrm{B_{dR}}$; the filtration on the right side is the tensor product filtration for the Hodge filtration of the $\eh$ de Rham complex and the filtration of $\mathrm{B_\dR}$.

\subsection{Ideas of proof}

In this subsection, we sketch the ideas of proof for the Theorem \ref{HTss}.

Let $K/\Q_p$ be a fixed complete algebraic closed non-archimedean field extension, and let $X$ be a rigid space over $K$.
One way to compute the $p$-adic \'etale cohomology $\cH^*(X_\et, \Q_p)$ is to use the pro-\'etale topology.
Scholze introduces the pro-\'etale site $X_\pe$ using the perfectoid geometry, together with a natural morphism of Grothendieck topologies
\[
\nu:X_\pe\rra X_\et.
\]
Assuming $X$ is proper and $K$ is complete and algebraically closed, Scholze (\cite{Sch12B}, \cite{Sch13}) shows that there exists an isomorphism
\[
\cH^*(X_\et,\Z_p)\otimes_{\Z_p}K\cong\cH^*(X_\pe,\pO_X),
\]
where $\pO_X$ is the completed pro-\'etale structure sheaf on $X_\pe$.
So by the Leray spectral sequence, the study of the $p$-adic \'{e}tale cohomology can be broken into two parts: the study of $R^i\nu_*\wh\calO_X$, and the understanding of its cohomology.

When $X$ is a smooth rigid space, 
we have $R^i\nu_*\pO_X=\Omega_{X/K}^i(-i)$ (see \cite[3.23]{Sch12B}).
Here the twist $``(-i)"$ means when $X=X_0\times_{K_0} K$ comes from a smooth rigid space $X_0$ over a finite extension $K_0/\Q_p$, the sheaf $R^i\nu_*\pO_X$ is isomorphic to the $(-i)$-th Tate twist of the $i$-th continuous differential sheaf over $K$, as a sheaf of module equipped with a Galois action.

In general when $X$ is not necessarily smooth over $K$, we use the aforementioned $\eh$-topology to extend the \'etale site to a locally smooth site $X_\eh$,
and denote the natural map of sites as $\pi_X:X_\eh \ra X_\et$.
We sheafify the presheaf of $K$-linear continuous differential forms to get the $\eh$-sheaf of differentials $\Omega_{ \eh}^i$.
Then we have the following result, connecting the higher direct image $R^i\nu_*\pO_X$ of the completed pro-\'etale structure sheaf with the $\eh$-sheaf of differentials:
\begin{theorem}[$\eh$-pro\'et spectral sequence]\label{main}
	Let $X$ be a rigid space over $K$.
	Then there exists an $E_2$ (Leray) spectral sequence of $\calO_X$-modules
	\[
	E_2^{i,j}=R^i\pi_{X *}\Omega^j_\eh(-j)\Rightarrow R^{i+j}\nu_*\wh{\calO}_X.
	\]
	When $X$ is smooth over $K$, the higher direct image $R^i\pi_{X *}\Omega^j_\eh$ vanishes for $j\in \NN$ and $i>0$, with 
	\[
	R^0\pi_{X *} \Omega_{ \eh}^j(-j)=\Omega_{X/K}^j(-j)
	\]
	being the continuous differential of $X$ over $K$, together with the $(-j)$-th Tate twist when $X$ is defined over a discrete valued subfield.
\end{theorem}
Here we note that in the scheme case, the analogous result holds for the Deligne-Du Bois complex: for a smooth algebraic variety $Y$ over $\CC$, the $j$-th graded factor $\ul\Omega_Y^j$ of the filtered Deligne-Du Bois complex is quasi-isomorphic to the $j$-th K\"ahler differential sheaf $\Omega_{Y/\CC}^j[-j]$ (cf. \cite{DB81}, \cite{Lee07}).

For the proof of the theorem, we first introduce the $v$-topology for the given rigid space $X$ (defined in \cite{Sch17}) in Section \ref{sec3}.
Here the $v$-topology  serves as a common extension of the pro-\'etale topology and the $\eh$-topology.
With the help of the descent result between the pro-\'etale site and the $v$-site, we reduce the problem to the study of the $\eh$-differentials for rigid spaces in the Theorem \ref{descent}. 
The rest of the proof will then be devoted to the study of the $\eh$-differential for smooth rigid spaces, in Section \ref{sec5}.
Here we follow the idea by Geisser (cf. \cite{Gei06}), showing the vanishing of a cone for the natural map
\[
\Omega_{X/K}^i\rra R\pi_{X *}\Omega_\eh^i
\]
by comparing $\eh$ and \'etale cohomologies, with the help of the covering structure in $X_\eh$ studied in Section \ref{sec2}.

\begin{remark}\label{proof of main, DVR case}
At such a moment, we mention that when the rigid space $X$ is defined over a discretely valued subfield of $K$,  the $\eh$-$\pe$ spectral sequence in the Theorem \ref{main} is Galois equivariant by the functoriality, and hence degenerates at its $E_2$-page by the result of Tate.
Together with the known finiteness for the $p$-adic \'etale cohomology for proper spaces, we in particular get the proof of Theorem \ref{HTss} in the case when $X$ comes from a discretely valued subfield.
For the reader's convenience, we also want to point out that this only needs the first half of the article, namely Section \ref{sec2} to Section \ref{sec4}.\footnote{Here to recover the Hodge--Tate decomposition for proper smooth rigid spaces as in Theorem \ref{main} (iii), it suffices to compare it with the known decomposition result for this special case as in \cite[Corollary 1.8]{Sch13}. This implies that under the assumption the cohomology of $\eh$-differentials coincides with that of continuous differential, and in particular does not require results in Section \ref{sec5}. We thank the referee for mentioning this to the author.  }
\end{remark}

Using Theorem \ref{main} we can show the coherence and cohomological boundedness of $R\nu_*\wh\calO_X$ in the non-smooth casem, in Section \ref{sec6}:
\begin{theorem}[Finiteness]\label{perfect}
	Let $X$ be a rigid space over $K$.
	Then the higher direct images $R^i\pi_{X *} \Omega^j_\eh$ are coherent and vanish unless $0\leq i,j\leq \dim(X)$.

\end{theorem}

After that, we study the general situation about the degeneracy of the $\eh$-$\pe$ spectral sequence.
As we saw in Remark \ref{proof of main, DVR case}, the Hodge--Tate decomposition follows from the Galois action when $X$ comes from a small subfield.
It is then natural to ask if the degeneracy holds for more general $X$, or even in the level of the derived category.
Our next main result confirms the splitting of the derived direct image $R\nu_*\wh\calO_X$ into its cohomology sheaves in the derived category, under an assumption of $X$ being strongly liftable (see Definition \ref{stronglift}).
The condition is satisfied for example when $X$ is defined over a discretely valued subfield of $K$ that has perfect residue field (Example \ref{discrete lift}), or when $X$ is proper over $K$ (Proposition \ref{proper lift}).

\begin{theorem}
	Assume $X$ is a quasi-compact and strongly liftable rigid space over $K$. 
	Then there exists a non-canonical quasi-isomorphism 
	\[
	R\nu_*\wh\calO_X\rra \bigoplus_j R\pi_{X *}(\Omega^j_\eh(-j)[-j]).
	\]
	In particular, the $\eh$-$\pe$ spectral sequence
	\[
	E_2^{i,j}=R^i\pi_{X *}\Omega^j_\eh(-j)\Rightarrow R^{i+j}\nu_*\wh{\calO}_X
	\]
	degenerates at the $E_2$-page.
\end{theorem}
In particular, when $X$ is smooth, the above decomposition degenerates into the following simpler form:
\begin{corollary}
	Assume $X$ is quasi-compact, strongly liftable and is smooth over $K$.
	Then there exists a non-canonical quasi-isomorphism
	\[
	R\nu_*\wh\calO_X\rra  \bigoplus_{i=0}^{\dim(X)} \Omega_{X/K}^i(-i)[-i].
	\]
	
\end{corollary}
In fact, the isomorphism is functorial among strong lifts of $X$; for the precise statement, we refer the reader to Theorem \ref{deg}.

When $X$ is smooth, we prove the theorem above by comparing the analytic cotangent complex of $X$ over $\Bdr$ with the derived direct image $R\nu_*\wh\calO_X$, which builds the bridge between the liftability and the degeneracy. 
For the non-smooth setting, we first generalize most of the previous results to smooth quasi-compact (truncated) simplicial rigid spaces, then we use the cohomological descent for the smooth $\eh$-hypercover.
The idea is to show that the above isomorphism for smooth rigid spaces is sufficiently functorial to enhance to simplicial settings, by constructing those enhanced maps by hand.
This is given in Section \ref{sec7}.

At the end of Section \ref{sec7}, we give an application of the degeneracy, on the higher direct image of the $\eh$-differential forms. 
By a recent work about the almost purity theorem in \cite{BS}, we can show that the derived direct image $R\nu_*\wh\calO_X$ lives in the cohomological degree $[0,\dim(X)]$, for any rigid space $X$ over $K$ (Proposition \ref{finite rev}).
Moreover, if we assume $X$ is proper over $K$, then the degeneracy of the $\eh$-pro\'et spectral sequence implies the vanishing of the $R^i\pi_{X *}\Omega_{ \eh}^j$ for $i+j>\dim(X)$.
Note that the vanishing of the cohomology sheaves is a local statement, so we can state the upshot as follows:
\begin{theorem}
	Let $X$ be a locally compactifiable rigid space over $K$.
	Then we have
	\[
	R^i\pi_{X *}\Omega_{ \eh}^j=0,~for~i+j>\dim(X).
	\]
\end{theorem}
We remark that the theorem gives an improvement of the cohomological boundedness in the Theorem \ref{perfect} for $X$ being locally compactifiable.

At the moment, it is natural to ask if the vanishing above is true for all rigid spaces.
Here we make a stronger conjecture as follows.
\begin{conjecture}
	Let $X$ be a rigid space over $K$.
	Then the $\eh$-pro\'et spectral sequence degenerates at the $E_2$-page.
	In particular, we have $R^i\pi_{X *}\Omega^j_\eh=0$ for $i+j> \dim(X)$.
\end{conjecture}
In the complex geometry, the analogous vanishing results hold for the graded pieces of the  Deligne--Du Bois complex, for a compact complex variety $X$ over $\CC$, as in \cite[Theorem 7.29]{PS08} (the statement is proved by Guillen-Navarro and Aznar-Puerta-Steenbrink separately, but it is more known as Kawamata-Viehweg type vanishing).
Moreover, with the help of the rigid GAGA theorem (\cite[Appendix A1]{Con06}), our result implies the vanishing of the Deligne-Du Bois complex for proper varieties over $\CC$.
\footnote{Our vanishing result on $\eh$ differentials directly implies the part (b) in \cite[Theorem 7.29]{PS08}. 
	To get the global vanishing in part (a), let $X$ be a projective variety of dimension $d$, $\calL$ the ample line bundle, $Y$ the affine cone of $X$, and $\Bl(Y)$ the blowup of $Y$ at the origin, where the latter is also the $\mathbb{A}^1$-bundle over $X$ defined by $\calL^{-1}$. We can then apply the distinguished triangle of the $\eh$-differential sheaves for the blowup square of $Y$. 
	Using  the projection map from $\Bl(Y)$ onto $X$ together with the vanishing of the Deligne--Du Bois complexes at $Y$ as in part (b), we can get part (a) via an induction argument starting from the top differential degree $d+1$.}
We want to mention that the proof here essentially makes use of the p-adic Hodge theory, while the proof in \cite{PS08} uses the mixed Hodge structure.
It is thus interesting to ask if we can produce more similar results from Hodge theory in the complex geometry, using tools from the $p$-adic Hodge theory instead.

In Section \ref{sec8}, we provide a comparison between the pro-\'etale cohomology and the $\eh$ de Rham cohomology for proper rigid spaces over a discretely valued subfield $K_0$, generalizing the smooth case developed in \cite{Sch13}.
The idea is to use the simplicial method and the cohomological descent developed in Section \ref{sec7}, together with the $\eh$-descent of differentials for smooth rigid spaces (Theorem \ref{descent}).
Our main result in this section is the following:
\begin{theorem}[(pro-\'etale)-$\eh$ de Rham comparison]
	Let $Y$ be a proper rigid space over $K_0$.
	Then there exists a $\Gal(K/K_0)$-equivariant filtered quasi-isomorphism
	\[
	R\Gamma(Y_\eh,\Omega_{\eh, /K_0}^\bullet)\otimes_{K_0} \mathrm{B_{dR}}\rra R\Gamma(Y_{K,\pe}, \BB_\dR),
	\]
	which generalizes the smooth case in \cite{Sch13}.
\end{theorem}
Here the left side above has the tensor product filtration, which is given by the Hodge filtration of the $\eh$ de Rham complex and the canonical filtration on $\mathrm{B_{dR}}$.
In the right side above, the cofficient $\BB_\dR$ is the de Rham sheaf for $X$ over the pro-\'etale site (see \cite[Definition 6.1]{Sch13}, and the Definition \ref{Bdr pe}) with its natural filtration.
As an upshot, we obtain the degeneracy for the $\eh$ version of the Hodge-de Rham spectral sequence at its $E_1$-page, assuming the condition of $Y$ above (Proposition \ref{HT-dR deg}):
namely, the following natural spectral sequence degenerates at the first page
\[
E_1^{p,q}=\cH^q(Y_\eh,\Omega_{ \eh, /K_0}^p)\Longrightarrow\cH^{p+q}(Y_\eh,\Omega_{ \eh, /K_0}^\bullet).
\]

Finally, we deduce the generalized Hodge-Tate spectral sequence and the decomposition in Section \ref{sec9}, using all of our results about the derived direct image $R\nu_*\wh\calO_X$ developed in this article.\\

In summary, the article is structured as follows.
In Section \ref{sec2}, we introduce the $\eh$-topology for rigid spaces.
Here we prove the local smoothness of this topology, and discuss in detail about the topological structure of $\eh$-coverings.
Then in Section \ref{sec3}, we introduce the pro-\'etale topology and $v$-topology, and discuss the necessary comparison theorems.
After these two sections, we connect those topologies together in Section \ref{sec4}.
Here we reduce the theorem of the $\eh$-pro\'et spectral sequence to the $\eh$-descent of differentials.
Section \ref{sec5} is devoted to the proof of the $\eh$-descent, with the help of which we get the comparison between $\eh$ cohomology and the singular cohomology, when the rigid space is coming from an algebraic variety.
In Section \ref{sec6}, we obtain the coherence and the cohomological boundedness using the $\eh$ hypercover.
In Section \ref{sec7}, we study the degeneracy phenomenon of the derived direct image $R\nu_*\wh\calO_X$, using the cotangent complex and the simplicial method.
As an application, we improve the cohomological boundedness of $R\pi_*\Omega_\eh^j$ for locally compactifiable spaces.
In Section \ref{sec8}, we give a comparison between the $\eh$ de Rham cohomology and the pro-\'etale cohomology, for proper rigid spaces over a discretely valued field.
Finally in Section \ref{sec9}, we explain the proof for the Hodge-Tate decomposition. 

As in the theory of perfectoid spaces, we use the language of adic spaces throughout the article.
We refer the reader to Huber's book \cite{Hu96} for basic results about adic spaces.

\section{\'Eh-topology}\label{sec2}
In this section, we introduce the $\eh$-topology and study its local structure.

\subsection{Rigid spaces}\label{subsec rig}
We first give a quick review about rigid spaces, following \cite{Hu96}.

Let $K$ be a complete non-archimedean extension of $\Q_p$.
Denote by $\RK$ the category of \emph{rigid spaces over $\Spa(K)$}; namely its objects consist of adic spaces that are locally of finite type over $\Spa(K,\calO_K)$.
Then for any $X\in \RK$, it can be covered by affinoid open subspaces, where each of them is of the form $\Spa(A,A^+)$ with $A$ being a quotient of the convergent power series ring $K\langle T_1,\ldots , T_n\rangle$ for some $n$. 
Here $A^+$ is an integrally closed open subring of $A$ that is of topologically finite type over $\calO_K$, and $A$ is complete with respect to the $p$-adic topology on $K$.
By the finite type condition it can be showed that any such $A^+$ is equal to $A^\circ$, the subring consisting of all power-bounded elements in $A$ (\cite[4.4]{Hu94}). 
So to simplify the notation we abbreviate $\Spa(A,A^\circ)$ as $\Spa(A)$ in this setting.
Unless otherwise mentioned, in the following discussion we always assume $X$ to be a rigid space over $\Spa(K)$.

For each adic space $X$, we can define two presheaves: $\calO_X$ and $\calO_X^+$, such that
when the affinoid space $U=\Spa(B,B^+)\subset X$ is open and $B$ is complete, we have
\[
\calO_X(U)=B,~\calO_X^+(U)=B^+.
\]
It is known that for any $X\in \RK$, both $\calO_X$ and $\calO_X^+$ are sheaves.
We could also define the coherent sheaves over rigid spaces, in a way that locally the category $\mathrm{Coh}(\Spa(B))$ of coherent sheaves over $\Spa(B)$ is equivalent to the category $\mathrm{Mod}_{fp}(B)$ of finitely presented $B$-modules (\cite[Theorem 2.3.3]{KL16}).
An important example of coherent sheaves are \emph{continuous differentials} $\Omega_{X/Y}^i$ for a map of rigid spaces $X\ra Y$, which is a coherent sheaf of $\calO_X$-module over $X$ (\cite[Section 1.6]{Hu96}).
\footnote{To simplify notations, we always use $\Omega^i$ to denote the continuous differential sheaves in our article, instead of algebraic ones. We will explicitly mention it when the latter comes up.}
Locally for a map of affinoid algebras $A\ra B$, it could be defined by taking the $p$-adic completion and inverting $p$ at the algebraic differential module of $B_0$ over $A_0$, where $A_0\ra B_0$ is a map of topologically of finite type rings of definition over $\calO_K$.

Recall that a $coherent~ideal$ is defined as a subsheaf $\calI$ of ideals in $\calO_X$ that is locally of finite presentation over $\calO_X$.
It is known when $X=\Spa(A)\in \RK$, there is a bijection between coherent ideals $\calI$ of $X$ and ideals of $A$, given by
\begin{align*}
\calI \lmt& \calI(X);\\
\wt I \reflectbox{$\lmt$}& I.
\end{align*}
Here $\wt I$ is the sheaf of $\calO_X$ module associated to
\[U\mapsto I\otimes_A \calO_X(U).
\]
For a coherent ideal $\calI$, we can define an $analytic~closed~subset$ of $X$,\footnote{It is also called the \emph{Zariski closed subset} in many literatures.
	Here we follow the one from \cite{BGR}.} by taking 
\[
Z:=\{x\in X~|~\calO_{X,x}\neq \calI_{X,x}\}=\{x\in X~|~|f(x)|=0,\forall f\in \calI\}.
\]
The subset $Z$ has a canonical adic space structure such that when $X=\Spa(A)$ and $\calI=\wt{I}$, we have $Z=\Spa(A/I, (A/I)^\circ)=:V(I)$.

\subsection{Blowups}\label{blowup}
Before we introduce the $\eh$-topology on $\RK$, we first recall the construction of blowup in rigid spaces, following \cite[4.1]{Con06}.

Let $X\in\RK$ be a rigid space, $\calI\subseteq \calO_X$ be a coherent sheaf, and $Z=V(\calI)$ be the analytic closed subset defined by $\calI$ as in Subsection \ref{subsec rig}. 
Following Conrad \cite[2.3, 4.1]{Con06}, we define the blowup of $X$ along $Z$ as follows:
\begin{definition}
	The $\mathrm{blowup}$ $\Bl_Z(X)$ of $X$ along $Z$ is the $X$-rigid space
	\[
	\Proj^\an(\bigoplus_{n\in \NN} \calI^n),
	\]
	which is the relatively analytified $\Proj$ of the graded algebra $\bigoplus_{n\in \NN} \calI^n$ over the rigid space $X$ (see \cite[2.3]{Con06}).
	
	It is called a $\mathrm{smooth~blowup}$ if the blowup center $Z$ is a smooth rigid space over $K$.
\end{definition}
\begin{remark}
	As a warning, our definition for the smooth blowup is different from some existing contexts, where both $X$ and $Z$ are required to be smooth.
\end{remark}

When $X=\Spa(A)$ is affinoid, the blowup of rigid space is in fact the ``pullback" of the schematic blowup $\Bl_I(\Spec(A))$ of $\Spec(A)$ at the ideal $I$ along the map $\Spa(A)\ra \Spec(A)$ of locally ringed spaces.
Precisely, consider the following natural diagram of locally ringed spaces
\[
\xymatrix{ \Bl_{V(I)}(\Spa(A)) \ar[r] \ar[d]& \Bl_{I}(\Spec(A)) \ar[d]\\
	\Spa(A) \ar[r]& \Spec(A).}
\]
Then following from the universal property of the relative analytification functor as in \cite[Theorem 2.2.5, Lemma 2.2.3]{Con06}, we have: for a given rigid space $Y$, there exists a functorial bijection between the collection of morphisms $h:Y\ra \Bl_{V(I)}(\Spa(A))$ of rigid spaces over $\Spa(A)$, and the collection of the following commutative diagram 
\[
\xymatrix{Y\ar[r]^{f~~~~~~~~~} \ar[d]_{g} & \Bl_{I}(\Spec(A))\ar[d] \\
	\Spa(A)\ar[r] & \Spec(A),}
\]
where $f$ is a map of locally ringed spaces and $g$ is a morphism of rigid spaces.

As what happens in the scheme theory, $\Bl_Z(X)$ satisfies the universal property (see \cite{Con06} after Definition 4.1.1): for any $f:Y\ra X$ in $\RK$ such that the pullback $f^*\calI$ is invertible, it factors uniquely through $\Bl_Z(X)\ra X$.
This leads to the isomorphism of the blowup map when it is restricted to the open complement $X\backslash Z$. 
Besides, it can be showed by universal property that rigid blowup is compatible with flat base change and analytification of schematic blowup (see \cite[2.3.8]{Con06}).
Precisely, for a flat map of rigid spaces $g:Y\ra X$ (i.e. $\calO_{Y,y}$ is flat over $\calO_{X,x}$ for any $y\in Y$ over $x\in X$), we have
\[
\Bl_{g^*\calI}(Y)=\Bl_Z(X)\times_X Y.
\]
When $X=X_0^\an$ is an analytification of a scheme $X_0$ of finite type over $K$, with $Z$ being defined by an ideal sheaf $\calI_0$ of $\calO_{X_0}$, we have
\[
\Bl_Z(X)=\Bl_{\calI_0}(X_0)^\an.
\]

We also note that the blowup map $\Bl_Z(X)\ra X$ is proper.
This is because by the coherence of $\calI$, locally $\calI$ can be written as a quotient of a finite free module, which (locally) produces a  closed immersion of $\Bl_Z(X)$ into a projective space over $X$, thus is proper over $X$.
Moreover, if both the center $Z$ and the ambient space $X$ are smooth over $K$, then the blowup itself $\Bl_Z(X)$ is also smooth.

\subsection{Universal homeomorphisms}
Another type of morphisms that will be used later is the universal homeomorphism.
\begin{definition}
	Let $f:X'\ra X$ is a morphism of rigid spaces over $K$.
	It is called a $\mathrm{universal}$ $\mathrm{homeomorphism}$ if for any morphism of rigid spaces $g:Y\ra X$, the base change $X'\times_X Y\ra Y$ is a homeomorphism.
\end{definition}
The following proposition gives a criterion  of universal homeomorphisms of rigid spaces:
\begin{proposition}\label{univ}
	Let $f:X\ra Y$ be a morphism of rigid spaces over $\Spa(K)$.
	Then it is a universal homeomorphism if and only if the following two conditions hold 
	\begin{enumerate}[(i)]
		\item $f$ is a finite morphism of rigid spaces.
		\item For any pair of affinoid open subsets $V=\Spa(A)\subset Y$ and $U=f^{-1}(V)=\Spa(B)$, the corresponding map of schemes
		\[
		\wt{f}:\Spec(B)\rra \Spec(A)
		\]
		is a universal homeomorphism of schemes.
	\end{enumerate}
\end{proposition}
\begin{proof}
	
	Assume $f$ is a universal homeomorphism.
	Let $x\in X$ be a rigid point.
	Since the map $f$ is quasi-finite, by \cite[1.5.4]{Hu96}  there exists open neighborhoods $U\subset X$ of $x$ and $V$ of $f(x)$ such that $f(U)\subset V$ and $f:U\ra V$ is finite.
	We may assume both $U$ and $V$ are connected.
	On the one hand, the finiteness of $f:U\ra V$ implies the image of $U$ is closed.
	On the other hand, as $f$ is a homeomorphism, $f(U)$ is open in $Y$, and thus open in $V$.
	Combine both of those, we see $V$ is exactly equal to $f(U)$ with $U=f^{-1}(V)$.
	So by the density of the rigid points $x\in X$, there exists an open covering $V_i$ of $Y$ such that $f^{-1}(V_i)$ is finite over $V_i$.
	Hence we get the finiteness of $f$.
	
	To check the universal homeomorphism for corresponding map of affine schemes, we recall from the Stack Project \cite[Tag 04DC]{Sta}, that $\wt{f}:\Spec(B)\ra \Spec(A)$ is a universal homeomorphism of schemes if and only if it is integral, universally injective and surjective.
	Since both $A$ and $B$ are $K$-algebras, where $K$ is an extension over $\Q_p$, it suffices to show the following claim.
	\begin{lemma}
		Let $f:\Spa(B) \ra \Spa(A)$ be a universal homeomorphism of affinoid rigid spaces.
		Then the induced map of affine schemes $\wt f:\Spec(B)\ra \Spec(A)$ is integral, bijective, and induces isomorphisms on their residues fields.
	\end{lemma}
\begin{proof}[Proof of Lemma]
	As we just showed above, the map of affinoid algebras $A\ra B$ is finite, thus $\wt{f}$ is a finite (hence integral) map of schemes.
	
	For the rest of the claim, we first consider its restrction on closed points.
	Let  $\wt y$ be a closed point of the scheme $\Spec(A)$, whose residue field $\kappa(\wt y)$ is a finite extension of the $p$-adic field $K$.
	The defining ideal of $\wt y$ in the scheme $\Spec(A)$ induces a unique rigid point $y$ of the rigid space $Y=\Spa(A)$, whose residue field is equal to $\kappa(\wt y)$.
	By assumption, the base change of the universal homeomorphism $f$ along the closed immersion $\{ y\} \ra Y$ induces a universal homeomorphism $X_y :=\Spa(\kappa(\wt y))\times_{\Spa(A)}\Spa(B) \ra \Spa(\kappa(\wt y))$, whose natural map to $X=\Spa(B)$ is a closed immersion.
	This implies that the reduced subspace of $X_y$ is a rigid point in $\Spa(B)$, and the corresponding closed subscheme inside of $\Spec(B)$ is supported at a unique closed point.
	Here we also notice that the residue field of $X_y$ is a finite separatble extension of $\kappa(\wt y)$.
	Moreover, applying the universal homeomorphism at the base change $X_y \times_{\Spa(\kappa(\wt y))} X_y \ra X_y$, we see the residue field of $X_y$ is isomorphic to $k(\wt y)$.
	As a consequence, the map $\wt f$ induces a bijection and isomorphisms of residues fields when restricted to their closed points.
	
	To finish the proof, it suffices to extend the claim for non-closed points.
	The bijection of $\wt f:\Spec(B)\ra \Spec(A)$ follows from the density of closed points.
	To see this, we may assume $\Spec(A)$ is irreducible.
	Then as $\wt f$ is a finite morphism whose image contains all closed points, we get the surjectivity of $\wt f$.
	For the injectivity, by the homeomorphism between $\Spa(B)$ and $\Spa(A)$, the scheme $\Spec(B)$ admits a unique irreducible component (hence a unique generic point), and has the same dimension as $\Spec(A)$.
	At last, as the induced map of $\wt f$ on the generic fields is finite and separable, its isomorphism follows from the bijection of points.
	So we are done.
	
\end{proof}

	Conversely, assume $f$ satisfies the two conditions as in the statement. 
	We first notice that both items in the statement are invariant under any base change of rigid soaces.
	We let $V=\Spa(A)$ and $f^{-1}(V)=\Spa(B)$ be two open affinoid open subsets of $Y$ and $X$ separately.
	Note that since $f$ is finite, for a morphism of affinoid rigid spaces $\Spa(C)\ra \Spa(A)$, the base change $\Spa(C)\times_{\Spa(A)}\Spa(B)$ is exactly $\Spa(B\otimes_A C)$ (\cite[1.4.2]{Hu96}).
	In particular, we see $\Spa(C)\times_{\Spa(A)}\Spa(B)\ra \Spa(C)$ is a finite morphism of rigid spaces, with the underlying map of schemes being a universal homeomorphism.
	As a consequence, both of the two conditions in the statement above are base change invariant, and to show $f:X\ra Y$ is a universal homeomorphism of rigid spaces, it suffices to show that $f$ itself is a homeomorphism.
	Moreover, by the finiteness, as the map $f$ is both closed and continuous, we are only left to show the bijectivity of $f$, as a map of rigid spaces.
	
	Now we pick any point $y\in Y$, and consider the completed residue field  with its valuation ring $(k(y),k(y)^+)$ of $y$.
	We take an open affinoid neighborhood $V=\Spa(A)$ of $y$ with $f^{-1}(V)=\Spa(B)$.
	Then the base change of the map $\Spec(B)\ra \Spec(A)$ of schemes gives
	\[
	\Spec(B\otimes_A k(y))\rra \Spec(k(y)),
	\]
	which is a universal homeomorphism by assumption.
	Here the target has exactly one point, and by finiteness we have $B\otimes_A k(y)=B\wh{\otimes}_A k(y)$.
	So by assumption the reduced subscheme $\Spec(B\wh\otimes k(y))_{\red}$ is equal to $k(y)$ (since they are of characteristic $0$).
	We then note that the adic spectrum 
	\[
	\Spa(B\wh{\otimes}_A k(y), B^\circ\otimes_{A^\circ} k(y)^+)
	\]
	 is exactly the preimage of $y$ in the rigid space $X$ along the morphism $f$.
	Notice that  the integral closure of $k(y)^+$ in $k(y)$ is contained in the quotient ring of  $B^\circ\otimes_{A^\circ} k(y)^+$ by its nilpotent elements, which has to be $k(y)^+$ itself (as the integral closure is contained in the field $k(y)$ and is finite over $k(y)^+$).
	In this way, the preimage $f^{-1}(y)$ has exactly one point $x$ whose residue field with valuation is equal to $(k(y),k(y)^+)$.
	Hence $f$ is bijective, and thus a homeomorphism.
\end{proof}
At last, when the target is assumed to be a smooth rigid space, there is no nontrivial universal homeomorphisms:
\begin{proposition}\label{uni homeo smooth}
	Let $X$ be a smooth rigid space, and $X'$ be a reduced rigid space.
	Then any universal homeomorphism $f:X'\ra X$ is an isomorphism.
\end{proposition}
\begin{proof}
	By the Proposition \ref{univ}, every universal homeomorphism $f:X'\ra X$ can be covered by morphisms of affinoid spaces $\Spa(B)\ra \Spa(A)$, where the underlying morphism of schemes $\Spec(B)\ra \Spec(A)$ is a universal homeomorphism. 
	So it suffices to show that when $X=\Spa(A)$ is a smooth affinoid rigid space over $\Spa(K)$, $A$ is a seminormal ring (so any universal homeomorphism $\Spec(B)\ra \Spec(A)$ from a reduced scheme is an isomorphism).
	But note that by the smoothness of $X$, $A$ is a regular ring (by \cite[1.6.10]{Hu96}, locally $X$ is \'etale over the adic spectrum of Tate algebras $K\langle T_i\rangle$, which is regular).
	So $A$ is normal, and thus seminormal.
\end{proof}
\subsection{$\mathrm{\acute{E}h}$-topology and its structure}
Now we can introduce the $\eh$-topology on $\RK$.
\begin{definition}\label{eh-top}
	The $\eh$-topology on the category $\RK$ is the Grothendieck topology such that the covering families are generated by the following types of morphisms:
	\begin{itemize}
		\item \'etale coverings;\\
		\item universal homeomorphisms;\\
		\item coverings associated to blowups: $\Bl_Z(X)\sqcup Z\ra X$, where $Z$ is a closed analytic subset of $X$.
	\end{itemize}
In the sense of Grothendieck pretopology in \cite[Expos\'e II.1]{SGA4}, this means that a family of maps $\{X_\alpha\ra X\}$ is in the set $Cov(X)$ if $\{X_\alpha\ra X\}$ can be refined by a finitely many compositions of the three classes of maps above.

	We denote by $\Reh$ the $\mathrm{big~\eh~ site}$ on $\RK$ given by the $\eh$-topology.
	For a given rigid space $X$ over $K$, we define $X_\eh$ as the localization  of $\Reh$ on $X$ (in the sense of \cite[Tag 00XZ]{Sta}, i.e. it is defined on the category of $K$-rigid spaces over $X$ with the $\eh$-topology.
\end{definition}
\begin{remark}\label{can-ref}
	\begin{enumerate}
		\item We notice that a covering associated to a blowup $\Bl_Z(X)\sqcup Z\ra X$ is always surjective: by the discussion in the Subsection \ref{blowup}, $\Bl_Z(X)\ra X$ is an isomorphism when restricted to $X\backslash Z$.
		\item Among the three classes of maps above, a covering associated to a blowup is not base change invariant in general.
		But note that for any morphism $Y\ra X$, the pullback of the blowup $X'=\Bl_Z(X)\ra X$ along $Y\ra X$ can be refined by the blowup
		\[
		\xymatrix{
			\Bl_{Y\times_X Z}(Y)\coprod Y\times_X Z \ar[rd] \ar[rdd] &&\\
			& Y\times_X \Bl_Z(X)\coprod Y\times_X Z \ar[r]\ar[d] & \Bl_Z(X)\coprod Z \ar[d]\\
			&Y\ar[r] &X.}
		\]
		We call $\Bl_{Y\times_X Z}(Y)\coprod Y\times_X Z$ the $canonical~refinement$ for~the~base~change~of~the~blow~up.
		
		\item Though denoted as $X_\eh$, this site is still a big site.
		As an extreme case, when $X=\Spa(K)$, the site $X_\eh$ is identical to $\Reh$.
	\end{enumerate}
\end{remark}
\begin{remark}
	Here we note that our definition of $\eh$-topology is different from $h$-topology.
	One of the main differences is that the $\eh$-topology excludes the ramified covering.
	
	For example, consider the $n$-folded cover map of the unit disc to itself $f:\BB^{1}\ra \BB^{1}$, which sends the coordinate $T$ to $T^n$.
	Then $f$ is a finite surjective map that is relatively smooth at all the other rigid points except at $T=0$, where it is ramified.
	If $f$ is an $\eh$-covering, by the Theorem \ref{ref} which we will prove later, $f$ can be refined by finite many compositions of coverings associated to smooth blowups and \'etale coverings.
	Notice that \'etale coverings are unramified maps that preserve the smoothness and dimensions.
	Moreover, smooth blowups of a one dimensional smooth rigid space are isomorphic to itself.
	In this way, such a finite composition will not produce a covering that is ramified at any rigid points, and we get a contradiction.
\end{remark}

\begin{example}\label{eg}
	Let $X$ be a rigid space.
		We take $X'=X_{\red}$ to be the reduced subspace of $X$. 
		Then $X'\ra X$ is a universal homeomorphism, which is then an $\eh$-covering.
		So in the $\eh$-topology, every space locally is reduced.

\end{example}
\begin{proposition}\label{irr}
	Let $X$ be a quasi-compact quasi-separated rigid space over $K$.
	Assume $X_i$ for $i=1,\ldots ,n$ are irreducible components of $X$ (see \cite{Con99}).
	Then the map
	\[
	\coprod_{i=1}^n X_i\ra  X
	\]
	is an $\eh$-covering.
\end{proposition}
\begin{proof}
	We first claim that the canonical map $\pi:\Bl_{X_1}(X)\ra X$ factors through $\bigcup_{i>1} X_i\ra X$; in other words, the image of $\pi$ is disjoint with $X_1\backslash (\bigcup_{i>1} X_i)$.
	
	Let $x\in X_1\backslash (\bigcup_{i>1} X_i)$ be any point.
	Take any open neighborhood $U\subset X_1\backslash (\bigcup_{i>1} X_i)$ that contains $x$.
	Then the base change of $\pi$ along the open immersion $U\ra X$ becomes
	\[
	\Bl_{U\bigcap X_1}(U)\ra U,
	\]
	by the flatness of $U\ra X$ and the discussion in the Subsection \ref{blowup}.
	But by our choice of $U$, the intersection $U\bigcap X_1$ is exactly the whole space $U$, which by definition leads to the emptiness of $\Bl_{U\bigcap X_1}(U)$.
	Thus the intersection of $\Bl_{X_1}(X)$ with $p^{-1}(U)$ is empty, and the point $x$ is not contained in the image of $\pi$.
	
	At last, note that the claim leads to the following commutative diagram
	\[
	\xymatrix{\Bl_{X_1}(X)\coprod X_1\ar[rr] \ar[rd]&& X\\
		& (\bigcup_{i>1} X_i)\coprod X_1 \ar[ru] &},
	\]
	which shows that the map $(\bigcup_{i>1} X_i)\coprod X_1\ra X$ is also an $\eh$-covering.
	Thus by induction on the number of components $n$, we get the result.
\end{proof}
Here we define a specific types of $\eh$-covering.
\begin{definition}
	For an $\eh$-covering $f:X'\ra X$ of rigid spaces, we say it is a \emph{proper $\eh$-covering} if $f$ is proper, and there exists a nowhere dense analytic closed subset $Z_{\red}\subset X_{\red}$ such that
	\[
	f|_{f^{-1}(X\backslash Z)_{\red}}:f^{-1}(X\backslash Z)_{\red}\rra X_{\red}\backslash Z_{\red}
	\]
	is an isomorphism.
\end{definition}
As an example, a covering associated to a blowup for the center being nowhere dense is a proper $\eh$-covering.

The idea of allowing blowups in the definition of the $\eh$ site is to make all rigid spaces $\eh$-locally smooth.
To make this explicit, we recall the Temkin's non-embedded disingularization: 
\begin{theorem}[\cite{Tem12}, 1.2.1, 5.2.2]\label{ROS}
	Let $X$ be a generically reduced, quasi-compact rigid space over $\Spa(K)$.
	Then there exists a composition of finitely many smooth blowups $X_n\ra X_{n-1}\ra\cdots \ra X_0=X$, such that $X_n$ is smooth.
\end{theorem}
\begin{corollary}[Local smoothness]\label{localsm}
	For any quasi-compact rigid space $X$, there exists a proper $\eh$-covering $f:X'\ra X_{\red}$ such that $X'$ is a smooth rigid space over $\Spa(K)$. 
	Moreover, $f$ is a composition of finitely many coverings associated to smooth blowups.
\end{corollary}
\begin{proof}
	By the Temkin's result, we may let $X_n\ra \cdots \ra X_0=X_{\red}$ be the blowup in that Theorem, such that the center of each $p_i:X_{i}\ra X_{i-1}$ is a smooth analytic subset $Z_{i-1}$ of $X_{i-1}$.
	Then by taking the composition of the covering associated to the blowup associated to each $p_i$, the map
	\[
	X':=X_n\coprod (\sqcup_{i=0}^{n-1} Z_i) \ra X_{\red}
	\]
	is a proper $\eh$-covering, such that $X'$ is smooth.
	So we get the result.
\end{proof}

At last, we give a useful result about the structure of the $\eh$-covering.
In order to do this, we need the embedded strong desingularization by Temkin:
\begin{theorem}[Embedded desingularization. Temkin \cite{Tem18}, 1.1.9, 1.1.13]
Let $X$ be a quasi-compact smooth rigid space over $\Spa(K)$, and $Z\subset X$ be an analytic closed subspace.
Then there exists a finite sequence of smooth blowups $X'=X_n\ra\cdots X_0=X$, such that the strict transform of $Z$ is also smooth.
\end{theorem}
\begin{corollary}\label{ref of blowup}
Any blowup $f:Y\ra X$ over a smooth quasi-compact rigid space $X$ can be refined by a composition of finitely many smooth blowups.
\end{corollary}
\begin{proof}
Assume $Y$ is given by $\Bl_Z(X)$, where $Z\subset X$ is a closed analytic subspace.
Then by the Embedded desingularization, we could find $g:X'\ra X$ to be a composition of finitely many smooth blowups such that the strict transform $Z'$ of $Z$ is smooth over $K$.
Here the total transform of $Z$ is $g^{-1}(Z)=Z'\cup E_Z$, where $E_Z$ is a divisor.
Next we could blowup $Z'$ in $X'$ and get $h:X''\ra X'$.
Note that $h$ itself is a smooth blowup.
In this way, the composition $h\circ g$ is a composition of finitely many smooth blowups that factorizes through $f:Y\ra X$, by the universal property of $f$ and the observation that the preimage of $Z$ along $h\circ g$ is the divisor 
\[
h^{-1}(Z')\cup h^{-1}(E_Z).
\]

\end{proof}
\begin{theorem}\label{ref}
	Let $X\in\RK$ be a quasi-compact smooth rigid space and $f:X'\ra X$ be an $\eh$-covering. 
	Then $f$ can be refined by a composition of finitely many \'etale coverings and coverings associated to smooth blowups  over $X$.
\end{theorem}
\begin{proof}
	By the definition of the $\eh$-topology, 
	 a given $\eh$-covering $f$ could be refined by a finitely many compositions of \'etale coverings, universal homeomorphisms, and coverings associated to blowups. 
	So up to a refinement we may write $f$ as $f:X'=X_n\ra X_{n-1}\ra\cdots\ra X_0=X$, where each transition map $f_i:X_i\ra X_{i-1}$ is one of the above three types of morphisms.
	
	Now we produce a refinement we want, by doing the following operations on $f$ starting from $i=1$:
	\begin{itemize}
		\item If $X_1\ra X_{0}$ is an \'etale morphism, then we are done for this $i=1$.
		\item If $X_1\ra X_0$ is a  universal homeomorphism, then by the Proposition \ref{uni homeo smooth} we may take the reduced subspace of $X_1$, which is isomorphic to $X_0$ and thus is smooth. 
		\item If $X_1\ra X_{0}$ is a covering associated to a blowup, then by the Proposition \ref{ref of blowup}, the associated blowup can be refined by finitely many compositions of smooth blowups.
		We let $X_1'\ra X_1$ be the disjoint union of that refinement with all of the centers.
		Then we take the base change of $X_n\ra \cdots \ra X_1$ along $X_1'\ra X_1$ and get a new coverings $X_n\times_{X_1} X_1'\ra \cdots \ra X_1'\ra X_0=X$, i.e.
		\[
		\xymatrix{
			X_n\times_{X_1} X_1' \ar[r] &\cdots \ar[r] & X_1' \ar[rd] \ar[d] &\\
			X_n \ar[r] & \cdots\ar[r] &X_1 \ar[r] & X_0=X.}
		\]
		Furthermore, starting at $j=2$, we do the following operation and increase $j$ by 1 each time:
		If $X_j\ra X_{j-1}$ is a covering associated to a blowup, we refine the map $X_j\times_{X_1} X_1'\ra X_{j-1}\times_{X_1} X_1'$ by its canonical refinement $X_j'\ra X_{j-1}\times X_1'$ (see the Remark \ref{can-ref}), and take the base change of the chain $X_n\times_{X_1} X_1'\ra\cdots \ra X_j\times_{X_1} X_1'$ along $X_j'\ra X_j\times_{X_1} X_1'$.\footnote{The covering associated to a blowup is not preserved under the base change, thus we need to adjust all of the maps in $X_n\times_{X_1} X_1'\ra\cdots \ra X_1'$ so that they will then become exactly those three types of morphisms.}
	\end{itemize}
After the discussion of the above three possibilities, $X_n\ra \cdots \ra X_0$ is refined by finitely many compositions $X_n'\ra\cdots X_1'\ra  X_0$ such that
\begin{itemize}
	\item $X_1'\ra X_0$ is a composition of finitely many \'etale coverings and coverings associated to smooth blowups;
	\item $X_n'\ra X_1'$ is a composition of $n-1$ $\eh$-coverings by \'etale coverings, coverings associated to blowups, or universal homeomorphisms.
\end{itemize}
In this way, we could do the above operation for $X_i'\ra X_{i-1}'$ and $i=2,\ldots$, each time get a new chain of coverings $X_n''\ra \cdots X_0$ such that $X_i''\ra X_0$ is a finite compositions of smooth blowups and \'etale coverings, and $X_n''\ra X_i''$ is a composition of $n-i$  coverings of three generating classes.
Hence after finitely many operations, we are done.
\end{proof}

\begin{corollary}
	Any $\eh$-covering of a quasi-compact rigid space $X$ can be refined by a composition
	\[
	X_2\ra X_1\ra X_0=X,
	\]
	where $X_1=X_{\red}$, the map $X_2\ra X_1$ is equal to finitely many compositions of \'etale-coverings and coverings associated to smooth blowups, and $X_2$ is smooth over $K$.
\end{corollary}
\begin{proof}
	Let $X'\ra X$ be a given $\eh$-covering.
	By the Example \ref{eg}, $X_1:=X_{\red}\ra X_0$ is an $\eh$-covering.
	And by the local smoothness of $\eh$-topology (Corollary \ref{localsm}), there exists a composition of finitely many coverings associated to smooth blowups $Y_1\ra X_1$, such that $Y_1$ is smooth.
	So the base change of $X'\times_X X_1\ra X_1$ along $Y_1\ra X_1$ becomes an $\eh$-covering whose target is smooth and quasi-compact.
	Hence by the Theorem \ref{ref} above, we could refine $X'\times Y_1\ra Y_1$ by $X_2\ra Y_1$, where the latter is a finitely many composition of \'etale coverings and coverings associated to smooth blowups.
	At last, notice that an \'etale map or a smooth blowup will not change the smoothness.
	Hence the composition $X_2\ra X_1\ra X_0$ satisfies the requirement.
\end{proof}

\section{Pro-\'etale topology and $v$-topology}\label{sec3}
In this section, we recall the pro-\'etale topology and $v$-topology over a given rigid space, in order to build the bridge between the pro-\'etale topology and the $\eh$-topology.
We follow mostly Scholze's foundational work \cite{Sch13} and \cite{Sch17}, together with the Berkeley's lecture notes \cite{SW20} by Scholze and Weinstein.\\

\subsection{Small $v$-sheaves}
Let $\Pd$ be the category of perfectoid spaces. 
They are adic spaces that have an open affinoid covering $\{\Spa(A_i,A_i^+),i\}$ such that each $A_i$ is a perfectoid algebra.
Since many of our constructions are large, we need to avoid the set-theoretical issure.
Following Section 4 in \cite{Sch17}, we fix an uncountable cardinality $\kappa$ with some conditions, and only consider those perfectoid spaces, morphisms, and algebras that are ``$\kappa$-small". 
We refer to Scholze's paper for details, and will follow this convention throughout the section. 

We first recall the $v$-topology defined on the category $\Pd$.
\begin{definition}[\cite{Sch17}, 8.1]
	The \emph{big $v$-site} $\Pd_v$ is the Grothendieck topology on the category $\Pd$, for which a collection $\{f_i:X_i\ra X,~i\in I\}$ of morphisms is a covering family if  for each quasi-compact open subset $U\subset X$, there exists a finite subset $J\subset I$ and quasi-compact open $V_i\subset X_i$, such that $|U|=\cup_{i\in J} f(|V_i|).$
	
	Here the index category $I$ is assumed to be $\kappa$-small.
\end{definition}
It is known that the $v$-site $\Pd_v$ are subcanonical; namely the presheaf represented by any $X\in \Pd$ is an $v$-sheaf.
Moreover, both integral and rational completed structure sheaves $\wh\calO^+:X\mapsto \wh\calO^+_X(X)$ and $\wh\calO:X\mapsto \wh\calO_X(X)$ are $v$-sheaves on $\Pd$ (\cite[8.6, 8.7]{Sch17}).

We then introduce a special class of $v$-sheaves that admits a geometric structure, generalizing the perfectoid spaces.
Consider the subcategory $\Pf$ of the category $\Pd$ consisting of perfectoid spaces of characteristic $p$.
We can equip $\Pf$ with the pro-\'etale topology and the $v$-topology to get two sites $\Pf_\pe$ and $\Pf_v$ separately.
\begin{definition}[\cite{Sch17}, 12.1]\label{sm-v}
	A \emph{small $v$-sheaf} is a sheaf $Y$ on $\Pf_v$ such that there is a surjective map of $v$-sheaves $X\ra Y$, where $X$ is a representabile sheaf of some $\kappa$-small perfectoid space in characteristic $p$.
\end{definition}
By the definition and the subcanonicality of the $v$-topology, any perfectoid space $X$ in characteristic $p$ produces a small $v$-sheaf.

Here is a non-trivial example.
\begin{example}[\cite{SW20}, 9.4]\label{until}
	Let $K$ be a $p$-adic extension of $\Q_p$, namely $K$ is complete with respect to a non-archimdean valuation extending that of $\mathbb{Q}_p$.
	Then we can produce a presheaf $\Spd(K)$ on $\Pf$, such that for each $Y\in \Pf$, we take
	\[
	\Spd(K)(Y):=\{isomorphism~classes~of~pairs~(Y^\sharp, \iota:(Y^\sharp)^\flat\ra Y)\},
	\]
	where $Y^\sharp$ is a perfectoid space (of characteristic $0$) over $K$, and $\iota$ is an isomorphism of perfectoid spaces identifying $Y^\sharp$ as an untilt of $Y$.
	It can be showed that $\Spd(K)$ is in fact a small $v$-sheaf.
	
	By the tilting correspondence, it can be showed that there is an equivalence between the category $\Pd_K$ of perfectoid spaces over $K$, and the category of perfectoid spaces $Y$ in characteristic $p$ with a structure morphism $Y\ra \Spd(K)$ (See \cite[9.4.4]{SW20}).
\end{example}

One of the main reasons we introduce small $v$-sheaves is that it brings both perfectoid spaces and rigid spaces into a single framework.
More precisely, we have the following fact:
\begin{proposition}[\cite{Sch17}, 15.5; \cite{SW20}, 10.2.3]\label{ad-dia}
Let $K$ be a $p$-adic extension of $\Q_p$ as in Example \ref{until}.
		There is a functor 
		\begin{align*}
		\{analytic~adic~spaces~over~\Spa(K)\}&\rra \{small~v-sheaves~over~\Spd(K)\} ;\\
		X&\lmt X^\diamond,
		\end{align*}
		such that when $X$ is a perfectoid space over $\Spa(K)$, the small $v$-sheaf $X^\diamond$ coincides with the representable sheaf for the tilt $X^\flat$.

		Moreover, the restriction of this functor to the subcategory of seminormal rigid spaces gives a fully faithful embedding:
		\[
		\{seminormal~rigid~spaces~over~\Spa(K)\}\rra \{small~v-sheaves~over~\Spd(K)\}.
		\]

\end{proposition}
Here we remark that every perfectoid space is seminormal (\cite[Theorem 3.7.4]{KL16}).

We can also define the ``topological structure on $X^\diamond$": in \cite{Sch17}, 10.1, Scholze defines the concept of being open, \'etale and finite \'etale for a morphism of pro-\'etale sheaves over $\Pd$.
In particular, for each small $v$-sheaf $X^\diamond$ coming from an adic space, we can define its small \'etale site $X^\diamond_\et$.
Those morphisms between small $v$-sheaves are compatible with maps of adic spaces, and we have
\begin{proposition}[\cite{Sch17}, 15.6]
	For each $X\in \RK$, the functor $Y\mapsto Y^\diamond$ induces an equivalence of small \'etale sites:
	\[
	X_\et\cong X^\diamond_\et,
	\]
	where the site on the left is the small \'etale site of the rigid space $X$ defined in \cite{Hu96}.
\end{proposition}
This generalizes the tilting correspondence of perfectoid spaces between characteristic $0$ and characteristic $p$.

\subsection{Pro-\'etale and $v$-topoi over $X$}\label{sec3.2}
In this subsection, we recall the small pro-\'etale site and the $v$-site associated to a given rigid space $X\in \RK$, for $K$ being a $p$-adic field.
Our goal is to produce a topology over $X$ that is large enough to include both pro-\'etale topology and $\eh$-topology together, and study the relation between their cohomologies.

We start by recalling basic concepts around the topology of small $v$-sheaves.

First recall that for a perfectoid space $X$, it is called \emph{quasi-compact} if every open covering admits a finite refinement; and it is called \emph{quasi-separated} if for any pair of quasi-compact perfectoid spaces $Y,Z$ over $X$, the fiber product $Y\times_X Z$ is also quasi-compact.

The concept of quasi-compactness and quasi-separatedness can be generalized to the pro-\'etale sheaves and  small $v$-sheaves.
A small $v$-sheaf $\calF$ is called \emph{quasi-compact} if for any family of morphisms $f_i:X_i \ra \calF,~i\in I$ such that $\coprod_{i\in I} X_i\ra \calF$ is surjective and $I$ is $\kappa$-small, it admits a finite subcollection $J\subset I$ such that $\coprod_{j\in J}X_j\ra \calF$ is surjective. 
Here $X_i$ are (pro-\'etale sheaves that are representable by) affinoid perfectoid spaces, 
The \emph{quasi-separatedness} for small $v$-sheaves is defined similarly as perfectoid spaces.

Now we are able to define the two topoi over a given rigid space $X$.
\begin{definition}\label{def-site}
	Let $X\in \RK$ be a rigid space over the $p$-adic field $K$.
	\begin{enumerate}[(i)]
		\item The \emph{small pro-\'etale~site} over $X$, denoted by $X_\pe$, is the Grothendieck topology on the full subcategory of pro-objects in $X_\et$ that are pro-\'etale over $X$, in the sense of \cite[Section 3]{Sch13}.
		Its covering families are defined as those jointly surjective pro-\'etale morphisms $\{f_i:Y_i\ra Y, ~i\in I\}$ such that for any quasi-compact open immersion $U\ra  Y$, there exists a finite subset $J\subset I$ and quasi-compact open $V_j\subset Y_j$ for $j\in J$, satisfying $|U|=\cup_{j\in J} f_j(|V_j|)$. 
		
		We call its topos the \emph{pro-\'etale topos} over $X$, denoted by $\Sh(X_\pe)$.
		
		\item The \emph{$v$-site over $X$} is defined as the site $\Pf_v|_{X^\diamond}$ of perfectoid spaces in characteristic $p$ over $X^\diamond$, with the covering structure given by the $v$-topology.
		Namely it is defined on the category of the category of pairs $(Y,f:Y\ra X^\diamond)$, where $Y$ is (the representable sheaf of) a perfectoid space in characteristic $p$, and $f:Y\ra X^\diamond$ is a map of $v$-sheaves over $\Pf_v$.
		A collection of maps $\{(Y_i, Y_i\ra X^\diamond) \ra (Y,Y\ra X^\diamond)\}$ is a covering in this site if $\{Y_i \ra Y\}$ is a covering in the $v$-site.
		
		We call its topos $\Sh(\Pf_v|_{X^\diamond})$ the \emph{$v$-topos over $X$}. 
	\end{enumerate}
\end{definition}
\begin{remark}\label{loc v}
		 The $v$-site $\Pf_v|_{X^\diamond}$ above is constructed as the localization (restriction) of the $v$-site $\Pf_v$ at the sheaf $X^\diamond$.
		 The general discussion of the localization of a site at a sheaf, which generalizes the localization of a site at an object, can be found for example in \cite[Tag 04GY]{Sta}.
		 Here we note that by \cite[Tag 0791]{Sta}, the $v$-topos over $X$ is isomorphic to the localization topos $\Sh(\Pf_v)|_{X^\diamond}$ of the $v$-topos $\Sh(\Pf_v)$ at the small $v$-sheaf $X^\diamond$.

\end{remark}
\begin{remark}\label{v in char 0/p}
Given a rigid space $X$, we can also form the characteristic zero analogue of the $v$-site $\Pd_v|_X$, on the category of perfectoid spaces over $X$ (cf. Definition \ref{def-site} (ii)).
The tilting correspondence and the definition of $X^\diamond$ induces a natural equivalence between the $v$-sites $\Pf_v|_{X^\diamond}$ and $\Pd_v|_X$, sending an affinoid perfectoid space $Z\ra X^\diamond$ onto the associated tilt $Z^\sharp \ra X$.
\end{remark}

Let $X\in \RK$ be a rigid space.
Then there is a natural morphism of topoi $\lambda=(\lambda^{-1}, \lambda_*): \Sh(\Pf_v|_{X^\diamond})\ra \Sh(X_\PE)$.
The inverse functor $\lambda^{-1}$ is computed via the functor $(-)^\diamond$ as in Proposition \ref{ad-dia}.
Precisely, when $Y\in X_\pe$ is affinoid perfectoid whose associated complete adic space is $\hat{Y}$, the inverse $\lambda^{-1}(Y)$ is the small $v$-sheaf $\hat Y^\diamond$ over $X^\diamond$, representable by the tilt $\hat{Y}^\flat$.
As affinoid perfectoid objects form a basis in $X_\pe$ (\cite[Proposition 4.8]{Sch13}), this allows us to extend $\lambda^{-1}$ to the whole category $X_\pe$.
In particular, by using the Galois descent as in \cite[Proposition 15.4]{Sch17}, for a rigid space $X'$ that is finite \'etale over $X$, we have $\lambda^{-1}(X')=X'^\diamond$.
Here we remark that by the loc. cit. the functor $\lambda^{-1}$ realizes a pro-\'etale presentation into an actual limit of $v$-sheaves: when $Y$ is affinoid perfectoid with a pro-\'etale presentation $\{Y_i\}$ over $X$, we have $Y^\diamond\cong \varprojlim Y_i^\diamond$. 

When $\{Y_i\ra Y\}$ is  a pro-\'etale covering of affinoid perfectoid objects over $X$, the inverse image $\{\hat Y_i^\diamond\ra \hat Y^\diamond\}$ forms  a $v$-covering of representable $v$-sheaves over $X^\diamond$.
For a general pro-\'etale sheaf $\calF$ over $X_\PE$, the functor $\lambda^{-1}$ sends $\calF$ to the $v$-sheaf associated to the presheaf
\[
Z\lmt \varinjlim_{\substack{Z\ra \hat W^\diamond~in~\Sh(\Pf_v|_{X^\diamond}),\\ \mathrm{affinoid~perfectoid}~W\in X_\PE}}
\calF(W).
\]
Here we note that when $Z$ is equal to the small $v$-sheaf $\hat Y^\diamond$ for $\hat Y$ a perfectoid space underlying a pro-\'etale object $Y$ over $X$, the above direct limit is $\calF(Y)$.	
On the other hand, the functor $\lambda_*$ is the direct image functor, given by
\[
\lambda_*\calG(Y)=\calG(\hat Y^\flat),~\mathrm{affinoid~perfectoid}~Y\in X_\PE.
\]

We define the \emph{untilted completed structure sheaves} $\wh\calO_v$ and $\wh\calO_v^+$ on $\Pf_v|_{X^\diamond}$, by sending $Z\ra X^\diamond$ to the following
\begin{align*}
\wh\calO_v(Z) &:= \wh\calO(Z^\sharp),\\
\wh\calO_v^+(Z) &:= \wh\calO^+(Z^\sharp),
\end{align*}
where $Z^\sharp$ is the untilt of $Z$ given by the map $Z\ra X^\diamond \ra \Spd(K)$, as in Proposition \ref{ad-dia}.
By \cite[Theorem 8.7]{Sch17}, both of them are sheaves on $\Pf_v|_{X^\diamond}$.
Here we notice that under the (tilting) equivalence in Remark \ref{v in char 0/p}, the sheaves $\wh\calO_v$ and $\wh\calO_v^+$ are sent to the completed structure sheaves $\wh\calO$ and $\wh\calO^+$ over $\Pd_v|_X$ in characteristic zero.

Furthermore, we have the following comparison result on completed pro-\'etale structure sheaves.
\begin{proposition}[(pro-\'etale)-v comparison]\label{pe-v}
	The direct image map induces the following canonical isomorphism of sheaves on $X_\PE$:
	\[
	\lambda_*\wh\calO_v^+ \rra  \wh\calO_X^+.
	\]
	Moreover, for $i>0$ the sheaf $R^i\lambda_*\wh\calO_v^+$ is almost zero.
	
	By inverting $p$, the similar results hold for $\lambda_*\wh\calO_v$ and $R^i\lambda_*\wh\calO_v$.
	In particular, the pro-\'etale cohomology of $\wh\calO_X$ satisfies the $v$-hyperdescent.
\end{proposition}
Here we follow the convention of the almost mathematics as in \cite[Section 3]{Sch17}.
\begin{proof}
	We first recall that for any quasi-compact analytic adic space $Y$ over $K$, there exists a pro-\'etale covering of $Y$ by perfectoid spaces (\cite[Lemma 15.3]{Sch17}).
	In particular, the pro-\'etale site $X_\PE$ admits a basis given by affinoid perfectoid spaces that are pro-\'etale over $X$.
	So it suffices to check the above isomorphism and vanishing condition for $Y\in X_\PE$ that are affinoid perfectoid.

The direct image of the untilted integral complete structure sheaf is the pro-\'etale sheaf associated to
\[
Y\lmt \Gamma(Y,~\lambda_*\wh\calO_v^+)=\Gamma(\hat Y^\diamond,~\wh\calO_v^+),
\]
where $Y\in X_\PE$ is affinoid perfectoid. 
But note that since $\hat Y^\diamond\cong \hat  Y^\flat$ is the representable sheaf of an affinoid perfectoid space over $X^\diamond$, by construction of $\wh\calO_v^+$ we have
\[
\Gamma(\hat Y^\diamond,\wh\calO_v^+)=\Gamma((\hat Y^\flat)^\sharp,\wh\calO^+).
\]
Here $\hat Y$ is the perfectoid space associated to the object $Y \in X_\PE$, and $\hat Y^\flat$ is the tilt of $\hat Y$.
So by the isomorphism of perfectoid spaces $(\hat Y^\flat)^\sharp\cong \hat Y$, we see the $\lambda_*\wh\calO_v^+$ is the pro-\'etale sheaf associated to the presheaf
\[
Y\lmt \Gamma(Y,\wh\calO^+),
\]
which is exactly the completed pro-\'etale structure sheaf over $X_\pe$.
Thus we get the equality.

For the higher direct image, we first note that $	R^i\lambda_*\wh\calO_v^+$ is the pro-\'etale sheaf on $X_\PE$ associated to the presheaf
\[
Y\lmt \cH^i_v(\hat Y^\flat,\wh\calO_v^+)
\]
for $Y$ being affinoid perfectoid in $X_\PE$.
By the construction of $\wh\calO_v^+$, the tilting correspondence $\Pf_v|_{X^\diamond} \cong \Pd_v|_X$ in Remark \ref{v in char 0/p} identifies the sheaf $\wh\calO_v^+$ over $\Pf_v|_{X^\diamond}$ with $\wh\calO^+$ over $\Pd_v|_X$.
In particular, we have the natural isomorphism of cohomology 
\[
\rmH^i_v(\hat Y^\flat, \wh\calO_v^+) \cong \rmH^i_v((\hat Y^\flat)^\sharp, \wh\calO^+) \cong \rmH^i_v(\hat Y, \wh\calO^+),
\]
which is almost zero by \cite[Propositino 8.8]{Sch17} and the assumption on $Y$.
So we are done.

\end{proof}

\section{\'Eh-pro\'et spectral sequence}\label{sec4}
In this section, we connect all of the topologies we defined together, and consider the $\eh$-pro\'et spectral sequence.

Let $X$ be a rigid space over $K$, for $K$ a complete and algebraically closed $p$-adic field.
We denote by $X_\eh$ the localization of the big $\eh$ site $\Reh$ at $X$.
Then the functor $Y\mapsto Y^\diamond$ induces a morphism of topoi
\[
\alpha:\Sh(\Pf_v|_{X^\diamond})\rra \Sh(X_\eh),
\]
where $\alpha^{-1}Y=Y^\diamond$ for $Y\ra X$ being a representable sheaf of an adic space.
We let $X_\et$ be the small \'etale site of $X$, consisting of rigid spaces that are \'etale over $X$.

Consider the following commutative diagram of topoi over $X$
\[
\xymatrix{
	\Sh(X_\PE) \ar[r]^\nu & \Sh(X_\et) \\
	\Sh(\Pf_v|_{X^\diamond}) \ar[u]^\lambda \ar[r]_\alpha& \Sh(X_\eh) \ar[u]_{\pi_X}.}
\]
Here we note that the diagram is functorial with respect to $X$.
In particular when $X=X_0\times_{K_0} K$ is a pullback of $X_0$ along a non-archimedean field extension $K/K_0$, the diagram is then equipped with a continuous action of $\Aut(K/K_0)$.

Now by the pro\'et-v comparison (Proposition \ref{pe-v}), we have
\begin{align*}
R\nu_*\wh{\calO}_X&=R\nu_*R\lambda_*\wh{\calO}_v\\
&=R\pi_*R\alpha_*\wh{\calO}_v.
\end{align*}
This induces a Leray spectral sequence
\[
E_2^{i,j}=R^i\pi_*R^j\alpha_*\wh{\calO}_v\Rightarrow R^{i+j}\nu_*\wh{\calO}_X.
\]

We then notice that by the above comparison again, the sheaf $R^j\alpha_*\wh\calO_v$ is the $\eh$-sheafification of the presheaf
\begin{align*}
X_\eh\ni Y \lmt &\cH^j(Y^\diamond_v,\wh{\calO}_v)\\
=&\cH^j(Y_\pe,\wh{\calO}_Y).
\end{align*}
When $Y$ is smooth, it is in fact the i-th continuous differential:
\begin{fact}[\cite{Sch12B}, 3.23]
	Let $Y$ be a smooth affinoid rigid space over $K$.
	Then we have a canonical isomorphism:
	\[
	\cH^j(Y_\pe,\wh{\calO}_Y)=\Omega^j_{Y/K}(Y)(-j),
	\]
	where the $\Omega^j_{Y/K}$ is the sheaf of the j-th continuous differential forms.
	Here the $``(-j)"$ means that the cohomology is equipped with an action of the Galois group $\Gal(K/K_0)$ by the Tate twist of weight $j$, when $Y=Y_0\times_{K_0} K$ is a base change of a smooth rigid space $Y_0$ over a complete discretely valued field $K_0$ whose residue field is perfect, and $K$ is complete and algebraically closed.
\end{fact}
In this way, by the local smoothness of the $\eh$-topology (Proposition \ref{localsm}) and the functoriality, the sheaf $R^j\alpha_*\wh\calO_v$ on $X_\eh$ is the $\eh$-sheaf associated to 
\[
smooth~Y\lmt \Omega^j_{Y/K}(Y)(-j).
\]
We call the $\eh$-sheaf associated to $Y\mapsto \Omega^j_{Y/K}(Y)$ the \emph{j-th $\eh$-differential over} $K$, denoted by $\Omega_{\acute{e}h,/K}^j$.
When the base field $K$ is clear, we use $\Omega_\eh^j$ to simplify the notation.
So substitute this into the spectral sequence above, we get the \emph{$\eh$-pro\'et spectral sequence}
\[
R^i\pi_{X *}\Omega_{eh}^j(-j)\Rightarrow R^{i+j}\nu_*\wh{\calO}_X.
\]

Our first main result is the following descent result of the $\eh$-differential, which we will prove in the next section.
\begin{theorem}[$\eh$-descent]\label{descent}
	Assume $X\in \RK$ is a smooth rigid space over $\Spa(K)$.
	Then for each $j\in \NN$, we have
	\[
	R\pi_{X *}\Omega_{eh}^j=R^0\pi_{X *}\Omega_{eh}^j[0]=\Omega_{X/K}^j.
	\]
\end{theorem}

\begin{remark}
	When $i=j=0$, the section $\calO_\eh(X)$ of the $\eh$-structure sheaf on any rigid space $X$ is $\calO(X^{sn})$, where $X^{sn}$ is the semi-normalization of $X_{\red}$.
	In other words, $\calO_\eh=\calO^{sn}$.
	This follows from \cite[10.2.3]{SW20}.
\end{remark}

\section{\'Eh-descent for the differentials}\label{sec5}
In this section, we prove the descent for the $\eh$-differential of a smooth rigid space $X\in \RK$, where $K$ is any $p$-adic field (not necessarily algebraically closed).
At the end of the section, we apply the $\eh$-descent to the case when $X$ is coming from an algebraic variety, to relate the $\eh$ cohomology to the Deligne-Du Bois complex (cf. \cite{DB81}, \cite[Section 7]{PS08}).

\subsection{$\mathrm{\acute{E}h}$-descent}
We will follow the idea in \cite{Gei06}, showing the vanishing of the cone $C$ for $\Omega_{X/K}^j\ra R\pi_*\Omega^j_{\eh}$ by comparing the \'etale cohomology and $\eh$ cohomology.

We first show the long exact sequence of continuous differentials for coverings associated to blowups. 
\begin{proposition}
	Let $f:X'\ra X$ be a blowup of a smooth rigid space $X$ along a smooth and nowhere dense closed analytic subset $i:Y\subset X$, with the pullback $g:Y'=X'\times_X Y\ra Y$.
	Then the functoriality of K\"ahler differentials induces the following distinguished triangle in the derived category of $X$:
	\[
	\xymatrix{
		\Omega_{X/K}^j\ar[r] &Rf_*\Omega_{X'/K}^j\bigoplus i_*\Omega_{Y/K}^j\ar[r] &i_*Rg_* \Omega_{Y'/K}^j.} \tag{$\ast$}
	\]
\end{proposition}
\begin{proof}
	We first note that since the argument is local on $X$, it suffices to show for any given rigid point $x\in X$, there exists a small open neighborhood of $x$ such that the result is true over that.
	So we may assume $X=\Spa(A)$ is affinoid, admitting an \'etale morphism to $\BB^n_K=\Spa(K\langle x_1,\ldots,x_n\rangle)$ by \cite[1.6.10]{Hu96}, and $Y$ is of dimension $r$, given by the $\Spa(A/I)$ for an ideal $I$ of $A$.
	Moreover, by refining $X$ to a smaller open neighborhood of $x$ if necessary, we could choose a collection of local parameters $f_1,\ldots,f_r$ and $g_1,\ldots, g_{n-r}$ at $x$, such that $\{g_l\}$ locally generates the ideal defining $Y$ in $X$. 
	In this way, by the differential criterion for \'etaleness (see \cite[1.6.9]{Hu96}), we could assume $Y\ra X$ is an \'etale base change of the closed immersion
	\[
	\BB^r_K\rra \BB^n_K.
	\]
	In particular, the blowup diagram for $\Bl_Y(X)\ra X$ locally is the \'etale base change of $\Bl_{\BB^r_K}(\BB^n_K)$ along $X\ra \BB^n_K$.
	
	Then we notice that the blowup of $\BB^n$ along $\BB^r$ is equivalent to the generic fiber of the $p$-adic (formal) completion of the blowup
	\[
	\A^r_{\calO_K}\rra \A^n_{\calO_K}.
	\]
	Furthermore,  as proved in  \cite[IV. Theorem 1.2.1]{Gro85}, there exists a natural distinguished triangle as follows
		\[
	\xymatrix{
		\Omega_{\A^n/\calO_K}^j\ar[r] &Rf_*\Omega_{\Bl_{\A^r}(\A^n)/\calO_K}^j\bigoplus i_*\Omega_{\A^r/\calO_K}^j\ar[r]  &i_*Rg_{*} \Omega_{\Bl_{\A^r}(\A^n)\times_{\A^n} \A^r/\calO_K}^j.}\tag{$\ast\ast$}
	\]
	Now we make the following claim:
	\begin{claim}
		The sequence $(\ast)$ for $(X,Y)=(\BB^n_K,\BB^r_K)$ can be given by the generic base change of the derived $p$-adic completion of the distinguished triangle $(\ast\ast)$.
	\end{claim}
Granting the Claim, since both derived completion and the generic base change are exact functors, we are done.
\begin{proof}[Proof of the Claim]
	We first notice since $\A^n_{\mathcal{O}_K}=\Spec(\calO_K[T_1,\ldots,T_n])$ is $p$-torsion free (thus flat over $\calO_K$), by \cite[Tag 0923]{Sta}, for a complex $C\in D(\A^n_{\mathcal{O}_K})$ its $p$-adic derived completion is given by 
	\[
	R\varprojlim(C\otimes^L_{\calO_K}\calO_K/p^m\calO_K).
	\]
	Moreover, note that differentials of $\A^n_{\mathcal{O}_K},\A^r_{\mathcal{O}_K},\Bl_{\A^r_{\mathcal{O}_K}}(\A^n_{\mathcal{O}_K})$ and $\Bl_{\A^r_{\mathcal{O}_K}}(\A^n_{\mathcal{O}_K})\times_{\A^n_{\mathcal{O}_K}} \A^r_{\mathcal{O}_K}$ over $\calO_K$ are all flat over $\mathcal{O}_K$.
	We use the notations $\A^n_m$ and $\A^r_m$ to abbreviate the schemes $\A^n_{\mathcal{O}_K/p^m}$ and $\A^r_{\mathcal{O}_K/p^m}$ separately.
	Then the derived base change of $(\ast\ast)$ along $\mathcal{O}_K\ra \calO_K/p^m$ can be written as the following 
	\[
	\xymatrix{
		\Omega_{\A^n_m/(\calO_K/p^m)}^j\ar[r] &Rf_*\Omega_{\Bl_{\A^r_m}(\A^n_m)/(\calO_K/p^m)}^j\bigoplus i_*\Omega_{\A^r_m/(\calO_K/p^m)}^j\ar[r]  &i_*Rg_{*} \Omega_{\Bl_{\A^r_m}(\A^n_m)\times_{\A^n_m} \A^r_m/(\calO_K/p^m)}^j.}\tag{$\ast\ast\ast$}
	\]
	Here we use the formula $\Omega^j_{Y/\calO_K}\otimes^L_{\mathcal{O}_K} \calO_K/p^m=\Omega^j_{Y_m/(\mathcal{O}_K/p^m)}$ for a smooth $\mathcal{O}_K$-scheme $Y$, 
	together with the derived base change formula for a proper morphism (\cite[Tag 07VK]{Sta}).
	Hence the derived $p$-adic completion of $(\ast\ast)$ is then computed by the derived limit of $(\ast\ast\ast)$ for $m\in \NN$.

	At last we discuss those derived limits term by term.
	For $C_m=\Omega_{\A^n_m/(\calO_K/p^m)}^j$ or $C_m=i_*\Omega_{\A^r_m/(\calO_K/p^m)}^j$, since their transition maps are surjective, the derived limit has no higher cohomology and we have
	\[
	R\varprojlim \Omega_{\A^n_m/(\calO_K/p^m)}^j=\Omega_{\wh\A^n/\calO_K}^j,~R\varprojlim i_*\Omega_{\A^r_m/(\calO_K/p^m)}^j=\Omega_{\wh\A^r/\calO_K}^j.
	\]
	
	For $C_m=Rf_*\Omega_{\Bl_{\A^r_m}(\A^n_m)/(\calO_K/p^m)}^j$ or $i_*Rg_{*} \Omega_{\Bl_{\A^r_m}(\A^n_m)\times_{\A^n_m} \A^r_m/(\calO_K/p^m)}^j$, 
	recall we have the formula of the derived functors 
	\[
	Rf_*R\varprojlim_m=R\varprojlim_mRf_*.
	\]
	Apply the formula, we get
	\begin{align*}
	R\varprojlim_m Rf_*\Omega_{\Bl_{\A^r_m}(\A^n_m)/(\calO_K/p^m)}^j&\cong Rf_*R\varprojlim_m \Omega_{\Bl_{\A^r_m}(\A^n_m)/(\calO_K/p^m)}^j\\
	&=Rf_*\Omega_{\Bl_{\wh\A^r}(\wh\A^n)/\calO_K}^j.
	\end{align*}
	The analogous holds for $C_m=i_*Rg_{*} \Omega_{\Bl_{\A^r_m}(\A^n_m)\times_{\A^n_m} \A^r_m/(\calO_K/p^m)}^j$.
	
		In this way, the derived limit of $(\ast\ast\ast)$ is isomorphic to
	\[
	\xymatrix{
		\Omega_{\wh\A^n/\calO_K}^j\ar[r] &Rf_*\Omega_{\Bl_{\wh\A^r}(\wh\A^n)/\calO_K}^j\bigoplus i_*\Omega_{\wh\A^r/\calO_K}^j\ar[r]  &i_*Rg_{*} \Omega_{\Bl_{\wh\A^r}(\wh\A^n)\times_{\wh\A^n} \wh\A^r/\calO_K}^j.} 
	\]
	
	At last, we take the base change of this distinguished  triangle along $\Z_p\ra \Q_p$, then we get $(\ast)$ for the pair of discs.

\end{proof}

\end{proof}
\begin{corollary}\label{leset}
	Under the above notation for the smooth blowup, we get a natural long exact sequence of \'etale cohomology of continuous differentials:
	\[
	\xymatrix{
		\cdots\ar[r]& \cH^j(X_\et,\Omega_{X/K}^i)\ar[r]& \cH^j(X'_\et,\Omega_{X'/K}^i)\bigoplus \cH^j(Y_\et,\Omega_{Y/K}^i)\ar[r]& \cH^j(Y'_\et,\Omega_{Y'/K}^i)\ar[r]& \cdots.}
	\]
\end{corollary}

Similarly there exists a long exact sequence of the covering associated to a blowup for the $\eh$-cohomologies:
\begin{proposition}\label{leseh}
	Let $f:X'\ra X$ be a morphism of rigid spaces over $\Spa(K)$, $Y\subset X$ be a nowhere dense analytic closed subspace, and $Y'=Y\times_X X'$ be the pullback.
	Let $X$ be separated.
	Assume they satisfy one of the following two conditions:
	\begin{enumerate}[(i)]
		\item $X'\ra X$ is a blowup along $Y$.
		\item $X$ is quasi-compact, $Y$ is an irreducible component of $X$, and $X'$ is the union of all the other irreducible components of $X$
	\end{enumerate}
	Then the functoriality of differentials induces a natural long exact sequence of cohomologies:
	\[
	\xymatrix{
		\cdots\ar[r]& \cH^j(X_\eh,\Omega_\eh^i)\ar[r]& \cH^j(X'_\eh,\Omega_\eh^i)\bigoplus \cH^j(Y_\eh,\Omega_\eh^i)\ar[r]& \cH^j(Y'_\eh,\Omega_\eh^i)\ar[r]& \cdots,}
	\]
	where $\Omega_\eh^i$ is the $\eh$-sheafification of the $i$-th continuous differential forms.
\end{proposition}
\begin{proof}
	For the rigid space $Z\in \RK$, we denote by $h_Z$ the abelianization of the $\eh$-sheafification of the representable presheaf 
	\[
	W\mapsto \Hom_{\Spa(K)}(W,Z).
	\]
	Then for an $\eh$ sheaf of abelian group $\calF$, we have
	\[
	\cH^j(Z_\eh,\calF)=\Ext^j_\eh(h_Z,\calF),
	\]
	since $h_Z$ is the final object in the category of sheaves of abelian groups over $Z_\eh$.
	So back to the question, it suffices to prove the exact sequence of $\eh$-sheaves
	\[
	0\rra h_{Y'}\rra h_{X'}\bigoplus h_Y \rra h_X \rra 0
	\]
	for above two conditions.
	
	Assume $\alpha:Z\ra X$ is a $K$-morphism.
	Then since $X'\coprod Y\ra X$ is an $\eh$-covering (see \ref{eh-top}, \ref{irr}), the element $\alpha\in h_X(Z)$ is locally given by a map $Z\times_X (X'\coprod Y)\ra (X'\coprod Y)$, which is an element in $h_{X'}(Z\times_X X')\bigoplus h_Y(Z\times_X Y) $, so we get the surjectivity.
	
	Now assume $(\sum n_r\beta_r,\sum m_s\gamma_s)$ is an element in $h_{X'}(Z)\bigoplus h_Y(Z)$ whose image in $0$ is $h_X(Z)$. 
	After refining $Z$ by an admissible covering of quasi-compact affinoid open subsets if necessary, we may assume $Z$ is quasi-compact affinoid.
	By taking a further $\eh$-covering of $Z$, we may also assume $Z$ is smooth and connected, given by $Z=\Spa(A)$ for $A$ integral.
	
	Then we look at the composition of those maps with $(f,i):X'\coprod Y\ra X$.
	\begin{itemize}
		\item Assume $f\circ\beta_1=f\circ \beta_2$ for some elements $\beta_i$.
		
		In the first setting of the Proposition, since $X'\ra X$ is a blowup along a nowhere dense (Zariski) closed subset, the restriction of $\beta_1$ and $\beta_2$ on the open subset $Z\backslash f^{-1}(Y)$ coincides.
		So by the assumption that $Z$ is integral (thus equal-dimensional), we see either the closed analytic subset $f^{-1}(Y)$ is the whole $Z$ and both $\beta_1$ and $\beta_2$ comes from $Z\ra Y\times_X X'=Y'$, or $f^{-1}(Y)$ is nowhere dense analytic in $Z$.
		If $f^{-1}(Y)$ is nowhere dense in $Z$, $\beta_1$ and $\beta_2$ agrees on an Zariski-open dense subset of $Z$.
		So by looking at open affinoid subsets of $X'$, the separatedness assumption implies that $\beta_1=\beta_2$ (see \cite[Chap. II Exercise 4.2]{Har77}).
		
		In the second setting, note that $X'\ra X$ is a closed immersion.
		So $f\circ \beta_1=f\circ \beta_2$ implies $\beta_1=\beta_2$.
		\item Assume $i\circ \gamma_1=i\circ \gamma_2$ for some elements $\gamma_i$.
		Then we get the identity of $\gamma_1$ and $\gamma_2$ again by the injectivity of the closed immersion $i:Y\ra X$.
		\item Assume there exists a equality $f\circ \beta_i=i\circ \gamma_j$.
		Since the composition $f\circ \beta_i$ is mapped inside of the analytic subset $Y\subset X$, the map $\beta_i:Z\ra X'$ factors through $Z\ra X'\times_X Y=Y'$.
		So $\beta_i$ comes from $h_{Y'}(Z)$, and by the injectivity of $Y\ra X$ again $\gamma_j$ comes from $h_{Y'}(Z)$.
	\end{itemize}
	In this way, by combining all of those identical $\beta_i$ and $\gamma_j$ and canceling the coefficients, the rest of $(\sum n_r\beta_r,\sum m_s\gamma_s)$ are all coming from $h_{Y'}(Z)$, thus the middle of the short sequence is exact.

	At last, injectivity of $h_{Y'}\ra h_{X'}\bigoplus h_Y$ follows from the closed immersion of $Y'\ra X'$.
	So we are done.
\end{proof}

\begin{remark}
	The part (ii) of the Proposition \ref{leseh} can be regarded as an $\eh$-version Mayer-Vietories sequence.
\end{remark}

\begin{proof}[Proof of the Theorem \ref{descent}]
	Now we prove the descent for the $\eh$-differential.
	
	Let $\Ret$ be the big \'etale site of rigid spaces over $K$.
	Namely it consists of rigid spaces over $K$, and its topology is defined by \'etale coverings.
	Then there exists a natural map of big sites $\pi:\Reh\ra \Ret$ that fits into the diagram
	\[
	\xymatrix{
		\Ret\ar[r]^{\iota_X} & X_\et\\
		\Reh \ar[u]^\pi  \ar[r]&X_\eh \ar[u]_{\pi_X}.}
	\]
	The sheaf $\Omega_\eh^i$ on $\Reh$ is defined as the $\eh$-sheafification of the continuous differential, which leads to the equality 
	\[
	\Omega_\eh^i=\pi^{-1} \Omega_{/K}^i,
	\]
	where $\Omega_{/K}^i$ is the $i$-th continuous differential on $\Ret$. 
	Besides, for any $Y\in \RK$, the direct image along $\Ret\ra Y_\et$ and $\Reh \ra Y_\eh$ are exact. 
	So it is safe to use $\calF|_{Y_\et}$ (resp. $\calF|_{Y_\eh}$) to denote the direct image of a sheaf on $\Ret$ (resp. $\Reh$) along those restriction maps, either in the derived or non-derived cases.
	
	Let $C$ be a cone of the adjunction map $\Omega_{/K}^i\ra R\pi_*\pi^{-1}\Omega_{/K}^i=R\pi_*\Omega_\eh^i$.
	It suffices to show the vanishing of $C$ when restricted to a smooth $X$; in other words, for each $X$ smooth over $K$, we want 
	\[
	\calH^j(C)|_{X_\et}=0, \forall j.
	\]
	We also note that as both $\Omega_{/K}^i$ and $R\pi_*\pi^*\Omega_{/K}^i$ has trivial cohomology of negative degrees, we have $\calH^j( C)|_{X_\et}=0$ for $j\leq -2$.
	In particular, $C$ is left bounded.
	
	Now we prove the above statement by contradiction.
	Assume $C$ is not always acyclic when restricted to the small site $X_\et$ for some smooth rigid space $X$ over $K$.
	By the left boundedness of $C$, we let $j$ be the smallest degree such that $\calH^j(C)|_{X_\et}\neq 0$ for some smooth $X$.
	Then $\calH^{j-l}(C)|_{Y_\et}=0$ for any $l>0$ and any smooth $Y$ over $K$.
	As this is a local statement, we fix an $X$ to be a smooth, connected, quasi-compact quasi-separated rigid space of the smallest possible dimension such that $\calH^j(C)|_{X_\et}\neq 0$.
	So by our assumption, there exists a nonzero element $e$ in the cohomology group
	\[
	\cH^0(X_\et,\calH^j(C))=\cH^j(X_\et, C).
	\]
	Here the equality of those two cohomologies follows from the vanishing assumption for $\calH^{j-l}(C)|_{X_\et}$ for $l>0$.
	
	We apply the preimage functor $\pi^{-1}$ to the triangle
	\[
	\Omega_{/K}^i\rra R\pi_*\pi^{-1}\Omega_{/K}^i\rra C,
	\]
	and get a distinguished triangle on $D(\Reh)$
	\[
	\pi^{-1}\Omega_{/K}^i\rra \pi^{-1}R\pi_*\pi^{-1}\Omega_{/K}^i\rra \pi^{-1}C.
	\]
	Note that since $\pi^{-1}$ is exact and the adjoint map $\pi^{-1}\ra \pi^{-1}\circ\pi_*\circ\pi^{-1}$ is an isomorphism, by taking the associated derived functors we get a canonical isomorphism
	\[
	\pi^{-1}\Omega_{/K}^i\cong \pi^{-1}R\pi_*\pi^{-1}\Omega_{/K}^i.
	\]
	So $\pi^{-1}C$ is quasi-isomorphic to $0$, and there exists an $\eh$-covering $X'\ra X$ such that $e$ will vanish when pullback to $X'$.

	Next we use the covering structure of the $\eh$-topology (Proposition \ref{ref}).
	By taking a refinement of $X'\ra X$ if necessary, we assume $X'\ra X$ is the composition 
	\[
	X'=X_m\ra X_{m-1}\ra \cdots \ra X_0=X,
	\]
	where $X_l\ra X_{l-1}$ is either a covering associated to a smooth blowup or an \'etale covering.
	
	Now we discuss the vanishing of the nonzero element $e$ along those pullbacks $X'=X_m\ra \cdots \ra X$.
	Assume $e|_{X_{l-1,\et}}$ is not equal to $0$ (which is ture when $l=1$).
	If $X_l\ra X_{l-1}$ is an \'etale covering, then since $e|_{X_{l-1,\et}}\in \cH^0(X_{l-1,\et},\calH^j(C))$ is a global section of a nonzero \'etale sheaf $\calH^j(C)$ on $X_{l-1,\et}$, the restriction of $e$ onto this \'etale covering will not be zero by the sheaf axioms. 
	If $X_l\ra X_{l-1}$ is a covering associated to a smooth blowup, we then make the following claim:
	\begin{claim}
		Under the above assumption, 
		the restriction $e|_{X_{l,\et}}$ in $\cH^0(X_{l,\et},\calH^j(C))=\cH^j(X_{l,\et},C)$ is not equal to $0$.
	\end{claim}
    Granting the Claim, since $X'\ra X$ is a finite composition of those two types of coverings, the pullback of $e$ to the cohomology group $\cH^0(X'_\et,\calH^j(C))$ cannot be $0$, and we get a contradiction.
    Hence $C|_{X_\et}$ must vanish in the derived category $D(X_\et)$ for smooth quasi-compact rigid space $X$, and we get the natural isomorphism
    \[
    \Omega_{X/K}^i\ra R\pi_{X *}\Omega_\eh^i,~\forall i.
    \]
    
    \begin{proof}[Proof of the Claim]
    	By assumption, since $e|_{X_{l-1,\et}}$ is nonzero, it suffices to show that the map of cohomology groups 
    	\[
    	\cH^j(X_{l-1,\et},C)\rra \cH^j(X_{l,\et},C)
    	\]
    	is injective.
    	
    	To simplify the notation, we let $X=X_{l-1}$, and $X'=X_l$ be the covering  ${\Bl_Y(X)}\sqcup Y\ra X$ associated to the blowup at the smooth center $Y\subset X$.
    	We let $Y'$ be the pullback of $Y$ along ${\Bl_Y(X)}\ra X$.
    	Since $X'\ra X$ is a covering associated to a blowup along a smooth subspace $Y$ of smaller dimension, by the two long exact sequence of cohomologies for differentials (Proposition \ref{leset}, \ref{leseh}), we get
    	\[
    	\xymatrix{
    		& \ar[d] & \ar[d]& \ar[d] &\\
    				\cdots\ar[r]& \cH^j(X_\et,\Omega_{X/K}^i)\ar[r] \ar[d]^{\tau_x}& \cH^j({\Bl_Y(X)}_\et,\Omega_{{\Bl_Y(X)}/K}^i)\bigoplus \cH^j(Y_\et,\Omega_{Y/K}^i)\ar[r]\ar[d]^{\tau_Y}_{\tau_{{\Bl_Y(X)}}}& \cH^j(Y'_\et,\Omega_{Y'/K}^i)\ar[r]\ar[d]^{\tau_{Y'}}& \cdots \\
    				\cdots\ar[r]& \cH^j(X_\eh,\Omega_\eh^i)\ar[r]\ar[d]& \cH^j({\Bl_Y(X)}_\eh,\Omega_\eh^i)\bigoplus \cH^j(Y_\eh,\Omega_\eh^i)\ar[r]\ar[d]& \cH^j(Y'_\eh,\Omega_\eh^i)\ar[r]\ar[d]& \cdots\\
    				\cdots\ar[r] & \cH^j(X_\et,C) \ar[r]\ar[d]& \cH^j({\Bl_Y(X)}_\et,C)\bigoplus\cH^j(Y_\et,C) \ar[r] \ar[d]& \cH^j(Y'_\et,C) \ar[r] \ar[d]&\cdots.\\
    				&  & & &}
    	\]
    	By the assumption of the $j$, since the cohomology sheaf $\calH^{j-l}(C)|_{Y',\et}=0$ for $l>0$, we have
    	\[
    	\cH^{j-1}(Y'_\et,C)=0.
    	\]
    	Besides, note that $\dim(X)$ is the smallest dimension such that $C|_{X_\et}$ is not quasi-isomorphic to $0$.
    	So both of $\cH^j(Y_\et,C)$ and $\cH^j(Y'_\et,C)$ are zero.
    	In this way, the third row above becomes an isomorphism between the following two cohomologies
    	\[
    	\cH^j(X_\et,C) \ra \cH^j({\Bl_Y(X)}_\et,C)\bigoplus0=\cH^j(X'_\et,C),
    	\]
    	and we get the injection.
    \end{proof}
\end{proof}

\begin{remark}
	In fact, the proof above works in a coarser topology, generated by rigid topology, universal homeomorphisms and coverings associated to blowups.
	This is because all we need are the local smoothness and the distinguished triangles for cohomology of differentials, which is a coherent cohomology theory.
	Moreover,  results here can be deduced from the pullback of this coarser topology to the $\eh$ topology. 
\end{remark}

\subsection{Application to algebraic varieties}\label{subsec ag}
Let $K$ be the field $\CC_p$ of $p$-adic complex numbers.
We fix an abstract isomorphism of fields between $\CC_p$ and $\CC$.
Our goal in this subsection is to relate the $\eh$ cohomology to the singular cohomology, when the rigid space  comes from an algebraic variety.

More precisely, we have:
\begin{theorem}\label{eh-sing comp}
	Let $Y$ be a proper algebraic variety over $K=\CC_p$, and let $X=Y^\an$ be its analytification, as a rigid space over $K$.
	Then there exists a functorial isomorphism
	\[
	\cH^i(X_\eh,\Omega_\eh^j)\cong \gr^j\cH^i_\Sing(Y(\CC),\CC),
	\]
	where $\cH^i_\Sing(Y(\CC))$ is the $i$-th singular cohomology of the complex manifold $Y(\CC)$ equipped with the Hodge filtration.
\end{theorem}
\begin{proof}
	
	Let $\rho:Y_\bullet\ra Y$ be a map from a simplicial smooth proper algebraic varieties over $K$ onto $Y$, such that each $Y_n\ra (\cosk_n Y_{\leq n})_{n+1}$ is a finite compositions of $\eh$ coverings associated to smooth blowups (but with algebraic varieties instead of rigid spaces in the Definition \ref{eh-top}).
	Then the analytification $\rho^\an:X_\bullet\ra X$ is an $\eh$ hypercovering of $X$ by smooth proper rigid spaces $X_n=Y^\an_n$.
	Moreover, the continuous differential sheaves $\Omega_{X_n/K}^j$ of $X_n$, which is a vector bundle over $X_n$, is canonically isomorphic to the sheafification of the differential sheaves $\Omega_{Y_n/K}^j$ of the algebraic variety $Y_n$ over $K$.
	
	Next we  apply the cohomological descent (Paragraph \ref{sec7.4 coh des}), and get the following natural quasi-isomorphism
	\[
	R\pi_{X *}\Omega_\eh^j\cong R\rho_*^\an R\pi_{X_\bullet *}\Omega_\eh^j.
	\]
	As each $X_n$ is smooth over $K$, by the Theorem \ref{descent} we have
	\[
	R\pi_{X_n *}\Omega_\eh^j\cong \Omega_{X_n/K}^j.
	\]
	In particular, the derived pushforward $R\pi_{X *}\Omega_\eh^j$ can be computed as
	\begin{align*}
	R\pi_{X *}\Omega_\eh^j&\cong R\rho_*^\an \Omega_{X_\bullet/K}^j\\
	&\cong R\rho_*^\an (\Omega_{Y_\bullet/K}^j)^\an.
	\end{align*}
	
	We then take the derived global section, to get
	\[
	R\Gamma(X_\eh, \Omega_\eh^j)\cong R\Gamma(Y^\an, R\rho_*^\an (\Omega_{Y_\bullet/K}^j)^\an).
	\]
	As all of the algebraic varieties $Y$ and $Y_n$ are proper over $K$ with each $\Omega_{Y_n/K}^j$ being coherent, by the rigid GAGA theorem (\cite[Appendix A1]{Con06}), we obtain a natural isomorphism
	\[
	R\Gamma(X_\eh, \Omega_\eh^j)\cong R\Gamma(Y,R\rho_* \Omega_{Y_\bullet/K}^j).
	\]
	
	Now by the construction, the map from the simplicial varieties $Y_\bullet\ra Y$ is a \emph{smooth h-hypercovering} in the sense of \cite{HJ14}.
	In particular, as proved in \cite[Theorem 7.12]{HJ14}, the complex $R\rho_*\Omega_{Y_\bullet/K}^j$ is naturally isomorphic to the $j$-th graded piece of the Hodge filtration of the Deligne-Du Bois complex $\ul\Omega_Y^\bullet$.
	So we may replace the derived pushforward, to get the isomorphism of cohomology groups as below
	\[
	\cH^i(X_\eh,\Omega_\eh^j)\cong \cH^i(Y,\ul\Omega_Y^j).
	\]
	In this way, as the right side is isomorphic to the $j$-th graded piece $\gr^j\cH^i_\Sing(Y(\CC),\CC)$  of the Hodge filtration of the singular cohomology (\cite[7.3.1]{PS08}), we get the isomorphism
	\[
	\cH^i(X_\eh,\Omega_\eh^j)\cong \gr^j\cH^i_\Sing(Y(\CC),\CC).
	\]
\end{proof}
We note that in the proof above, the comparison is compatible with the differential maps on both side.
So the above leads to a comparison between the $\eh$ de Rham cohomology and the singular cohomology, when $X$ is coming from an algebraic variety.
\begin{corollary}\label{eh-sing whole}
	Let $Y$ be a proper algebraic variety over $K=\CC_p$, and let $X=Y^\an$ be its analytification, as a proper rigid space over $K$.
	Then there exists a functorial filtered isomorphism
	\[
	\cH^i(X_\eh,\Omega_\eh^\bullet)\cong \cH^i_\Sing(Y(\CC),\CC),
	\]
	where $\cH^i_\Sing(Y(\CC),\CC)$ is the $i$-th singular cohomology of the complex manifold $Y(\CC)$, equipped with the Hodge filtration. 
\end{corollary}

\begin{remark}
	Let $X=Y^\an$ be the analytification of a proper algebraic variety $Y$ over $\CC_p$ as above.
	The proof of the Theorem \ref{eh-sing comp} in fact implies that the $\eh$ cohomology $\cH^i(X_\eh,\Omega_\eh^j)$ of $\Omega_\eh^j$ is isomorphic to the $h$ cohomology $\cH^i(Y_h,\Omega_h^j)$ (via \cite[Corollary 6.16]{HJ14}), for the $h$ cohomology of the scheme $Y$ introduced in \cite{HJ14}.
	So every computation  for proper algebraic variety $Y$ in \cite{HJ14} can be used to compute the $\eh$ cohomology of the rigid space $Y^\an$.
\end{remark}

\section{Finiteness}\label{sec6}
In this section, we prove a finiteness result about the $R\nu_*\wh{\calO}_X$ for $X$ being a rigid space, namely the coherence and the cohomological boundedness of $R\nu_*\wh\calO_X$, where $K$ is an arbitrary $p$-adic field.

Assume $X$ is a rigid space over $K$.
\begin{proposition}[Coherence]\label{coh}
	The sheaf of $\calO_X$-module $R^n\nu_*\wh{\calO}_X$ is coherent.
\end{proposition}
\begin{proof}
	By the $\eh$-pro\'et spectral sequence $R^i\pi_{X *}\Omega^j_\eh\Rightarrow R^{i+j}\nu_*\wh{\calO}_X$, it suffices to show the coherence for each $R^i\pi_{X *}\Omega^j_\eh(-j)$.
	We then assume $X$ is reduced, since the direct image along $X_{\red}\ra X$ preserves the coherence of modules.
	Then by the local smoothness (Proposition \ref{localsm}), there exists an $\eh$-hypercover $s:X_\bullet\ra X$, such that each $X_k$ is smooth and the map $s_k:X_k\ra X$ is proper.
	Here we notice that each $R\pi_{X_k *}\Omega^j_{\eh}=\Omega^j_{X_k/K}$ is coherent on $X_k$ by the assumption of $X_k$ and Theorem \ref{descent}.
	So the properness of $s_k:X_k\ra X$ implies that each $R^qs_{k *}\Omega^j_{X_k/K}$ is coherent over $\calO_X$.
	On the other hand, thanks to the cohomological descent (see the discussion later in Paragraph \ref{sec7.4 coh des})), the derived direct image $Rs_{*}R\pi_{X_\bullet *}\Omega_{\eh}^j$ along the $\eh$-hypercover $X_\bullet\ra X$ is quasi-isomorphic to the $R\pi_{X *}\Omega^j_\eh$.
	In this way, the spectral sequence associated to the simplicial object $s:X_\bullet\ra X$ provides 
	\begin{align*}
	E_1^{p,q}=R^qs_{p *}\Omega_{X_p/K}^j &\Rightarrow \calH^{p+q}(Rs_*R\pi_{X_\bullet *}\Omega_{\eh}^j)\\
	&=R^{p+q}\pi_{X *}\Omega^j_\eh,
	\end{align*}
	where each term on the left side is coherent over $X$.
	Hence the sheaf $R^{p+q}\pi_{X *} \Omega^j_\eh$ is coherent on $X$.
\end{proof}

Next we consider the cohomological boundedness of the derived direct image.

\begin{theorem}[Cohomological boundedness]\label{coh-bound0}
		For a quasi-compact rigid space $X$, the cohomology 
		\[
		\cH^i(X_\eh,\Omega_{\eh}^j)
		\]
		vanishes except $0\leq i,j\leq \dim(X)$.

\end{theorem}
\begin{remark}
	The analogous statement about the boundedness of the Hodge numbers for varieties over the complex number $\CC$ is proved by Deligne (\cite[Theorem 8.2.4]{Del74}).
\end{remark}
\begin{proof}
		We do this by induction on the dimension of the $X$.
	When $X$ is of dimension $0$, the reduced subspace $X_{\red}$ is a finite disjoint union of $\Spa(K')$ with $K'/K$ finite, which is smooth over $\Spa(K)$.
	So by the local reducedness of the $\eh$-topology and the vanishing of the higher direct image of  $\iota:X_{\red}\ra X$, the case of dimension $0$ is done by the Theorem \ref{descent}.
	
	We then assume the result is true for all quasi-compact rigid spaces of dimensions strictly smaller than $\dim(X)$.
	By the local smoothness (Proposition \ref{localsm}) and the vanishing of the higher direct image along $X_{\red}\ra X$ again, we may assume $X$ is reduced and there exists a composition of finitely many blowups at smooth centers
	\[
	X'=X_n\rra \cdots \rra X_1\rra X_0=X,
	\]
	such that $X'=X_n$ is smooth over $\Spa(K)$. 
	Observe that by the property of the $\eh$ differential for smooth spaces, the sheaf $R^i\pi_{X_n *}\Omega^j_\eh$ is zero except when $i=0$ and $0\leq j\leq \dim X_n=\dim X$.
	So the claim is true for $X_n$ as $\cH^i(X_{n,\eh},\Omega_\eh^j)=\cH^i(X_n,\Omega_{X_n/K}^j)$.
	Moreover, to prove the claim for $X$, it suffices to show that if the result is true for $X_{l+1}$, then it is true for $X_{l}$, where $X_{l+1}\ra X_l$ is the $l$-th blowup at a nowhere dense analytic subspace.
	
	To simplify the notation, we let $X'=X_{l+1}$, $X=X_l$, $f:X'\ra X$ be the blowup, $i:Y\ra X$ be the inclusion map of the blowup center, and $Y'$ be the preimage $X'\times_X Y$ with the map $g:Y'\ra X'$.
	By the assumption, $\cH^i(X'_\eh,\Omega_\eh^j)$  vanishes unless $i,j\leq \dim(X)=\dim(X')$.
	Furthermore, thanks to the induction hypothesis we have $\cH^i(Y'_\eh,\Omega_\eh^j)=0$  unless $i,j\leq \dim(Y')<\dim(X)$.
	Now we consider the distinguished triangle of the $\eh$-cohomology (Proposition \ref{leseh})
	\[
	\xymatrix{
		\cdots\ar[r]& \cH^i(X_\eh,\Omega_\eh^j)\ar[r]& \cH^i(X'_\eh,\Omega_\eh^j)\bigoplus \cH^i(Y_\eh,\Omega_\eh^j)\ar[r]& \cH^i(Y'_\eh,\Omega_\eh^j)\ar[r]& \cdots.}
	\]
	We discuss all possible cases:
	\begin{itemize}
		\item If $j>\dim(X)$, then since the blowup center $Y$ is nowhere dense in $X$, we have $j \dim(X)> \dim(Y')$.
		So by induction hypothesis on dimensions, both $\cH^i(Y_\eh,\Omega_\eh^j)$ and $\cH^{i-1}(Y'_\eh,\Omega_\eh^j)$ vanish for every $i$.
		Moreover by the assumption on $X'$, we know the vanishing of $\cH^i(X'_\eh,\Omega_\eh^j)$, $i\in \NN$.
		So the long exact sequence leads to the vanishing for $\cH^i(X_\eh,\Omega_\eh^j)$ if $j>\dim(X)$.
		\item If $i>\dim(X)$, then since $i-1> \dim(X)-1\geq \dim(Y')$, by induction hypothesis on dimensions again we have $\cH^{i-1}(Y'_\eh,\Omega_\eh^j)$ and $\cH^i(Y_\eh,\Omega_\eh^j)$ are zero.
		Similarly we have the vanishing of the $\cH^i(X'_\eh,\Omega_\eh^j)$ by the assumption on $X'$.
		In this way, the long exact sequence implies that the cohomology $\cH^i(X_\eh,\Omega_\eh^j)$ is zero for $i>\dim(X)$ and any $j\in \NN$.
		
	\end{itemize}
	
\end{proof}
\begin{corollary}\label{coh-bound}
	
	Let $X$ be a rigid space over $K$.
	Then unless $0\leq i,j\leq \dim(X)$, the higher direct image $R^i\pi_{X *}\Omega_{\eh}^j$ vanishes. 
\end{corollary}
\begin{proof}
	We only need to note that the sheaf $R^i\pi_{X *}\Omega_{\eh}^j$ is the coherent sheaf on $X$ associated to the presheaf
	\[
	U\lmt \cH^i(U_\eh,\Omega_{\eh}^j),
	\]
	for $U\subseteq X$ open and quasi-compact.

\end{proof}

\begin{corollary}\label{coh-bound2}
	Let $X$ be a rigid space over $K$.
	Then each cohomology sheaf $R^i\nu_*\wh{\calO}_X$ is coherent over $\mathcal{O}_X$, and vanishes unless $0\leq i\leq 2\dim(X)$.
\end{corollary}

We will improve the above two corollaries for locally compactifiable rigid spaces in the Proposition \ref{finite rev} and the Proposition \ref{finite impro}, using the degeneracy result developed in the next section and the almost purity theorem in \cite{BS}.

\section{Degeneracy Theorem}\label{sec7}
In this section, we show the degeneracy of the spectral sequence $R^i\pi_{X *}\Omega_\eh^j(-j)\Rightarrow R^{i+j}\nu_*\wh{\calO}_X$, under  the condition that $X$ is strongly liftable.
More precisely, we use the cotangent complex to show the existence of a quasi-isomorphism
\[
R\nu_*\wh\calO_X\cong \bigoplus_j R\pi_{X *}(\Omega_{\eh}^j(-j)[-j]),
\]
assuming $X$ is strongly liftable (see Definition \ref{stronglift}).
The condition is satisfied if it is  proper over $K$ (Proposition \ref{proper lift}), or defined over a discretely valued subfield that has a perfect residue field (Example \ref{discrete lift}).

\subsection{Cotangent complex for formal schemes and adic spaces}\label{cot-def}
In this subsection, we first recall basics about the cotangent complex for adic spaces. 
The detailed discussion about the analytic cotangent complexes of formal schemes and adic spaces can be found in \cite[Section 7.1-7.3]{GR03} (over $\mathcal{O}_K$ and $K$) and \cite[Section 5.1, 5.2]{Guo20} (over $\Ainf/\xi^e$ and $\Bdr/\xi^e$).

Let $R_0$ be a $p$-adically complete ring; namely there exists a continuous morphism of adic rings $\Z_p\ra R_0$ with $R_0$ being $p$-adically complete.
Recall for a map of complete $R_0$-algebras $A\ra B$ that are $p$-torsion free, we can define its complete cotangent complex $\wh\LL_{B/A}$ as the term-wise $p$-adic completion of the usual cotangent complex $\LL_{B/A}$.
Here $\LL_{B/A}$ is given by the corresponding complex of the simplicial $B$-module  $\Omega_{P_\bullet(B)/A}^1\otimes_{P_\bullet(B)} B$, where $P_\bullet(B)$ is the standard $A$-polynomial resolution of $B$.
The image of $\wh\LL_{B/A}$ in the derived category of $B$-modules is the $p$-adic derived completion of $\LL_{B/A}$, which lives in cohomological degrees $\leq 0$ such that
\[
\cH^0(\wh\LL_{B/A})=\wh\Omega_{B/A}^{1},
\]
where $\wh\Omega_{B/A}^{1}$ is the continuous differential of $B$ over $A$ and is defined as the $p$-adic completion of the algebraic K\"ahler differential $\Omega_{B/A}^1$.
We note that the construction of the complex $\wh\LL_{B/A}$ is functorial with respect to complete $R_0$-algebras $A\ra B$.
So when $\calX\ra \Spf(R_0)$ is an $R_0$-formal scheme that is $p$-torsion free, we can construct a complex of presheaves, which assigns the complex $\wh\LL_{B/R_0}$ to an affinoid open subset $\Spf(B)$ in $\calX$.
The \emph{complete cotangent complex} $\wh\LL_{\calX/R_0}$ for a $p$-torsion free $R_0$-formal scheme $\calX$ is the actual complex of sheaves defined by sheafifying the above complex of presheaves term-wisely.

Now following the construction as in \cite[Section 7.2]{GR03},
for a map of $p$-adic affinoid Huber pairs $(A,A^+)\ra (B,B^+)$, 
we define its \emph{analytic cotangent complex $\LL^\an_{(B,B^+)/(A,A^+)}$} as the colimit
\[
\underset{\substack{A_0\ra B_0\\ A_0,B_0~open~bounded}}{\colim} \wh\LL_{B_0/A_0}[\frac{1}{p}],
\]
where the colimit is indexed over the set of all maps of rings of definition $A_0\ra B_0$ in $A^+\ra B^+$,
and $\wh\LL_{B_0/A_0}$ is the complete cotangent complex for a map of $p$-complete rings as above.
We often use the notation $\wh\LL_{B/A}^\an$ instead of $\LL^\an_{(B,B^+)/(A,A^+)}$ to simplify the notation, when the choice of the rings $A^+$ and $B^+$ is clear from the context.
The construction is functorial with respect to the pair $(A,A^+)\ra (B,B^+)$, and we can sheafify it to define the analytic cotangent complex $\LL^\an_{X/Y}$ for a map of  adic spaces $X\ra Y$.
Here the complex $\LL^\an_{X/Y}$ is a complex of sheaves of $\mathcal{O}_X$-modules that lives in non-positive cohomological degrees, 
such that 
\[
\cH^0(\LL_{X/Y}^\an)=\Omega_{X/Y}^1,
\]
with the latter $\Omega_{X/Y}^1$ is the continuous differential for the map of rigid spaces $X\ra Y$.

\begin{remark}\label{simplified cot}
	In many cases where the base ring is fixed, the colimit in the construction above can be simplified.
	
	For example, let $(R,R^+)$ be either a reduced topologically finite type algebra over a $p$-adic field, or $(\Ainf[\frac{1}{p}],\Ainf)$ (the definition of $\Ainf$ will be recalled in the next paragraph \ref{para Bdr}),
	and let $R_0$ be the fixed ring of definition $R^+$ (this is guaranteed by the reducedness of $A$, and the boundedness of $R^\circ$ by for example \cite[6.2.4/1]{BGR}) or $\Ainf$ separately.
	Then for an affinoid $R$-algebra $(B,B^+)$, we have the following natural quasi-isomorphism
	\[
	 \underset{\substack{R_0\ra  B_0\\ B_0~open~bounded}}{\colim} \wh\LL_{B_0/R_0}[\frac{1}{p}] \rra \wh\LL_{B/R}^\an,
	\]
	where the colimit ranges only among rings of definition of $(B,B^+)$.
	This is because in both cases the ring $R_0$ is the largest ring of definition, so the index systems of colimits are cofinal to the one in the original definition.
	
	Moreover, if in addition the integral subring $B^+$ of the Huber pair $(B,B^+)$ is bounded, then the above colimit can be further simplified into one single complex by the following quasi-isomorphism
	\[
	\wh\LL_{B^+/R_0}[\frac{1}{p}] \rra \wh\LL_{B/R}^\an,
	\]
	which follows from the same reason about the index system.
\end{remark}

\begin{remark}
	The construction of the analytic cotangent complexes here are slightly different from the one used in \cite{GR03} and \cite{Guo20}: the colimit in the definition of $\LL^\an_{(B,B^+)/(A,A^+)}$ above is over the set of \emph{all} rings of definitions, while the ones in loc. cit. are over the set of \emph{topologically of finite type} rings of definitions.
	The reason we include all rings of definition here is to extend the construction to perfectoid algebras, which are almost never topologically of finite type.
	
	To see those two constructions of analytic cotangent complexes for topologically finite type algebras $A$ over $\mathrm{B_{dR}^+}/\xi^e$ coincide, it suffices to notice that any ring of definitions $A_0$ of $A$ is contained in a ring of definition $A_1$ that is topologically finite type over $\Ainf/\xi^e$.
	When $A$ is reduced (hence is topologically of finite type over $K$), the subring of power-bounded elements $A^\circ$ is the largest ring of definition which is topologically of finite type over $\calO_K$ (apply \cite[6.4.1/5]{BGR} at a surjection $K\langle T_i\rangle\ra A$).
	
	For the general case when $A$ is not necessarily reduced, this can be seen as follows:
	Let $A_0$ be the given ring of definition, $I_0$ be the nilpotent radical of $A_0$, and $A_1$ be a ring of definition that is topologically of finite type over $\Ainf/\xi^e$ whose quotient by its nilpotent radical $I_1$ is $(A_\red)^\circ$.
	Here we note that by the $p$-torsion-freeness of $A_1/I_1$ and \cite[Lemma 13.4, (iii, b)]{BMS}, the ideal $I_1$ is finitely generated,  and we can assume the ideal $I_1$ is generated by a finite set of elements $g_j$, $1\leq j\leq m$.
	Moreover, the subring $A_0$ of $A$ is contained in the union of open subrings $\bigcup_{n\in \NN} A_1[\frac{1}{p^n} I_1]$, as the latter is the preimage of $(A_\red)^\circ$ in $A$ along the surjection $A \ra A_\red$, and $(A_0)_\red\subseteq (A_\red)^\circ$ by the last paragraph for reduced rings.
	By assumption the subring $A_0$ is bounded, so we could choose an integer $n$ large enough such that $A_0\subset A_1[\frac{1}{p^n}I_1]$.
	Therefore the claim follows as the ring of definition $A_1[\frac{1}{p^n}I_1]$ admits a surjection from $\Ainf/\xi^e\langle T_i, S_j\rangle$, where the map extends a surjection $\Ainf/\xi^e\langle T_1,\ldots, T_l\rangle \ra A_1$ and sends $S_j$ onto $\frac{1}{p^n} g_j$.
	Here we remind the reader that the construction makes sense as each $\frac{1}{p^n}g_j$ is nilpotent and in particular is topologically nilpotent.
\end{remark}

\paragraph{Lifting obstruction}\label{para Bdr}
One of the most important applications of the cotangent complex is the deformation problem.

Let $(R,R^+)$ be a $p$-adically complete Huber pair over $\mathbb{Q}_p$. 
Assume $I$ is a closed ideal in $R^+$.
We define $S$ as the adic space $\Spa(R/I,\wt{R^+/I})$, and $S'$ as the adic space $\Spa(R/I^2,\wt{R^+/I^2})$, where  $\wt{R^+/I}$ and $\wt{R/{I^2}}$ are integral closures.
Let $X$ be a flat $S$-adic space.
Then a \emph{deformation} of $X$ along $S\ra S'$ is defined as a closed immersion $i:X\ra X'$ of $S'$-adic spaces with $X'$ being flat over $S$, such that the defining ideal is $i^*I$.
Namely we have the following cartesian diagram
\[
\xymatrix{
	X\ar[r]^i \ar[d] & X'\ar[d]\\
	S \ar[r] & S'.}
\]

We now focus on the case where the coefficient $(R,R^+)$ is specified as below.
Assume $K$ is a complete and algebraically closed $p$-adic field, and let $X$ be a quasi-compact rigid space over $\Spa(K)$. 
Recall that the ring $\Ainf$ is defined as the ring of the Witt vectors $W(\varprojlim_{x\mapsto x^p} \calO_K)$.
There is a canonical surjective continuous map $\theta:\Ainf\ra \calO_K$, with kernel being a principal ideal $\ker(\theta)=(\xi)$ for some fixed element $\xi\in \Ainf$.
We then recall that the de Rham period ring $\Bdr$ is defined to be the $\xi$-adic completion of $\Ainf[\frac{1}{p}]$.
Here we abuse the notation and denote by $\theta:\Bdr\ra K$ the canonical surjection, which is induced from $\theta:\Ainf\ra \calO_K$ as above.
Note that for $n\geq 1$, we have $\Bdr/(\xi)^n=\Ainf[\frac{1}{p}]/(\xi)^n$, which is a complete Tate ring over $\Q_p$ (\cite[Section 1.1]{Hu96}).
In particular the deformation of any rigid space $X/K$ along $(\Bdr/\xi^n,\Ainf/\xi^n)\ra (K,\calO_K)$ is the same as the deformation along $(\Ainf[\frac{1}{p}]/\xi^n,\Ainf/\xi^n)\ra (K,\calO_K)$.

We then note that the deformation theory along $(\Bdr,\Ainf)\ra (K,\calO_K)$ only depends on the $p$-adic topology.
Precisely, we have the following observation:
\begin{lemma}
	Let $X$ be a topologically of finite type, $p$-torsion free formal scheme over $\Ainf/\xi^N$ for some $N\in\NN$.
	Let $X_n$ be the base change of $X$ along $\Ainf\ra \Ainf/\xi^{n+1}$.
	Then we have the following quasi-isomorphism
	\[
	 \wh\LL_{X/\Ainf}\rra R\varprojlim_n \wh\LL_{X_n/(\Ainf/\xi^{n+1})},
	\]
	where $\wh\LL$ is denoted by the $p$-adic complete cotangent complex given at the beginning of this section.
\end{lemma}
\begin{proof}
	We may assume $X$ is affinoid.
	Let $T_n$ be the $p$-adic formal scheme $\Spf(\Ainf/\xi^{n+1})$, and let $T$ be the $p$-adic formal scheme $\Spf(\Ainf)$.
	Consider the following sequence of $p$-adic formal schemes
	\[
	X_n\rra T_n\rra T.
	\]
	Then by taking the distinguished triangle of transitivity for usual cotangent complexes, we get
	\[
	\LL_{T_n/T}\otimes^L_{\calO_{T_n}} \calO_{X_n} \rra \LL_{X_n/ T}\rra \LL_{X_n/ T_n}. \tag{$\ast_n$}
	\]
	Here we note that the triangle remains distinguished in $D(\calO_{X_n})$ after the derived $p$-adic completion. 
	
	We then take the derived inverse limit (with respect to $n$) of the $p$-adic derived completion of $(\ast_n)$.
	When $n\geq N$, we have $\wh\LL_{X_n/ T}=\wh\LL_{X/T}$.
	Besides, since $X_n$ is the base change of $X$ along $T_n\ra T$, as complexes we have the equality $\wh\LL_{T_n/T}\otimes^L_{\calO_{T_n}} \calO_{X_n}=\wh\LL_{T_n/T}\otimes^L_{\calO_T} \calO_X$.
	So by taking the inverse limit with respect to $n$, we get
	\[
	R\varprojlim_n (\wh\LL_{T_n/T}\otimes^L_{\calO_T} \calO_X)\rra \wh\LL_{X/T} \rra R\varprojlim_n \wh\LL_{X_n/ T_n}.
	\]
	
	But note that the inverse system $\{\wh\LL_{T_n/T}\otimes^L_{\calO_T} \calO_X\}_n$ is in fact acyclic.
	This is because the cotangent complex $\wh\LL_{T_n/T}$ is isomorphic to $(\xi)^n/(\xi^{2n})[1]$, while the composition of transition maps
	\[
	(\xi)^{2n}/(\xi^{4n})[1]\rra (\xi)^{2n-1}/(\xi^{4n-2})[1]\rra \cdots \rra (\xi)^n/(\xi^{2n})[1]
	\]
	is $0$.
	In this way, by the vanishing of its $R\varprojlim_n$, we get the quasi-isomorphism we need.
\end{proof}
This Lemma allows us to forget the complicated natural topology on $\Bdr$ when we look at the deformation along $\Bdr\ra K$.
So throughout the article, we will consider the adic space $\Spa(\Ainf[\frac{1}{p}],\Ainf)$ that is only equipped with the $p$-adic topology, and any cotangent complex that has $\Ainf$ or $\Ainf[\frac{1}{p}]$ as the base will be considered $p$-adically.

Let $S$ and $S'$ be the adic space $\Spa(\Ainf[\frac{1}{p}]/\xi)$ and $\Spa(\Ainf[\frac{1}{p}]/\xi^2)$ separately. 
Here we note that $S$ is also equals to $\Spa(K)$.
Denote by $i$ the map $X\ra \Spa(\Ainf[\frac{1}{p}],\Ainf)$.
We let $\calO_X(1)$ be the free $\calO_X$ module of rank one, defined by
\[
i^* (\xi\mathcal{O}_S)=\calO_X\otimes_{\Ainf[\frac{1}{p}]} \xi \Ainf[\frac{1}{p}] = \xi/\xi^2\calO_X.
\]
When $X$ is defined over a discretely valued subfield, it admits the Galois action of the Hodge--Tate weight $(-1)$.

Our first result is about the relation between the deformation of $X$ and the splitting of the cotangent complex.

\begin{proposition}\label{ob}
	Let $X$ be a rigid space over $S=\Spa(K)$.
	Then a flat lifting $X'$ of $X$ along $S\ra S'$ induces a section $s_X$ of $\LL_{S/S'}^\an\otimes_{S} \calO_X\ra \LL_{X/S'}^\an$ in the distinguished triangle for the transitivity 
	\[
	\LL_{S/S'}^\an\otimes_{S} \calO_X\rra \LL_{X/S'}^\an\rra \LL_{X/S}^\an.\]
	
	Moreover, assume $X'\ra Y'$ is a map of flat adic spaces over $S'$, which lifts the map $f:X\ra Y$ of rigid spaces over $K$.
	Then the induced sections above are functorial, in the sense that the following natural diagram of sections commute:
	\[
	\xymatrix{
		\LL_{Y/S'}^\an \ar[r]^{s_Y~~~~~~~} \ar[d] & \LL_{S/S'}^\an\otimes_{S'} \calO_Y \ar[d]\\
		Rf_*\LL_{X/S'}^\an \ar[r]^{Rf_*(s_X)~~~~~~~~} & Rf_*(\LL_{S/S'}^\an\otimes_{S} \calO_X).}
	\]
\end{proposition}
\begin{proof}
	The base change diagram 
	\[
	\xymatrix{
		X \ar[r] \ar[d] & X' \ar[d]\\
		S \ar[r] & S'}
	\]
	induces the following two sequences of maps
	\begin{align*}
	X\ra S\ra S',\\
	X\ra X'\ra S'.
	\end{align*}
	We take their distinguished triangles of transitivity (\cite[Corollary 5.2.18]{Guo20}), and get 
	\[
	\xymatrix{
		& \LL_{X/S}^\an &\\
		\LL_{X'/S'}^\an\otimes_{\calO_{X'}} \calO_X \ar[r] & \LL_{X/S'}^\an \ar[r]\ar[u] & \LL_{X/X'}^\an\\
		& \LL_{S/S'}^\an\otimes_{\calO_{S}} \calO_X \ar[u], &}
	\]
	where both the vertical and the horizontal are distinguished. 

	Following \cite[Tag 09D8]{Sta}, we could extend the above to a bigger diagram
	\[
	\xymatrix{
		\LL_{X'/S'}^\an\otimes_{\calO_{X'}} \calO_X \ar[r]	& \LL_{X/S}^\an \ar[r] & E\\
		\LL_{X'/S'}^\an\otimes_{\calO_{X'}} \calO_X \ar[r] \ar@{=}[u]& \LL_{X/S'}^\an \ar[r]^{\alpha_X}\ar[u] & \LL_{X/X'}^\an \ar[u] \\		
		& \LL_{S/S'}^\an\otimes_{\calO_{S}} \calO_X \ar[u] \ar@{=}[r] &\LL_{S/S'}^\an\otimes_{\calO_{S}} \calO_X \ar[u]_{\beta_X},} \tag{$\ast$}
	\]
	where $E$ is the cone of 
	\[
	(\LL_{X'/S'}^\an\otimes_{\calO_{X'}} \calO_X\oplus  \LL_{S/S'}^\an\otimes_{\calO_{S'}} \calO_X) \rra  \LL_{X/S'}^\an,
	\]
	which fits into the diagram such that all of the vertical and horizontal triangles are distinguished.
	
	We then make the following claim.
	\begin{claim}
		The cone $E$ is isomorphic to $0$ in the derived category.
	\end{claim}
\begin{proof}[Proof of the Claim]
	By construction, since the right vertical triangle above is distinguished, it suffices to show that
	\[
	\beta_X:\LL_{S/S'}^\an \otimes_{\calO_{S}} \calO_X\rra \LL_{X/X'}^\an
	\]
	is a quasi-isomorphism.
	
	We may assume $X=\Spa(B,B^+)$ and $X'=\Spa(A,A^+)$ is affinoid, such that $A\otimes_{\Ainf/\xi^2} \calO_{K}=A/\xi=B$.
	Then by construction, the above map can be rewrite as 
	\[
	\LL_{S/S'}^\an \otimes_{K} B\rra \underset{\substack{A_0\ra B_0\\ A_0,B_0~open~bounded}}{\colim} \wh\LL_{B_0/A_0}[\frac{1}{p}],
	\]
	for $A_0\ra B_0$ being all pairs of the rings of definition of $(A,A^+)$ and $(B,B^+)$ separately.
	
	We then note that for a single pair $A_0\ra B_0$ such that $B_0\cong A_0/\xi$, the map
	\[
	\rho: \LL_{S/S'}^\an \otimes_{K} B \rra  \wh\LL_{A_0/B_0}[\frac{1}{p}].
	\]
	is a quasi-isomorphism:
	by the surjectivity of $\Ainf/\xi^2\rra \Ainf/\xi=\calO_K$ and $B_0\ra A_0$, applying \cite[7.1.31]{GR03} \footnote{Though the statement there is for topologically finite type algebras over $\calO_K$, the proof works for topologically finite type, p-torsion free algebras over $\Ainf/\xi^n$ as well. 
	}, we have
	\begin{align*}
	\wh\LL_{\calO_K/(\Ainf/\xi^2)}&\cong \LL_{\calO_K/(\Ainf/\xi^2)};\\
	\wh\LL_{A_0/B_0}& \cong \LL_{A_0/B_0}.
	\end{align*}
	So under the choice of $A_0$ and $B_0$ the map $\rho$ can be rewritten as a map of usual cotangent complexes
	\[
	\rho: \LL_{\calO_K/(\Ainf/\xi^2)}\otimes_{\calO_K} B_0[\frac{1}{p}]\rra \LL_{A_0/B_0}[\frac{1}{p}].
	\]
	But by assumption, $A=A_0[\frac{1}{p}]$ is flat over $\Ainf/\xi^2[\frac{1}{p}]$, while $B=B_0[\frac{1}{p}]$ is given by $A\otimes_{\Ainf/\xi^2} \calO_K$.
	Hence by the flatness of inverting $p$ and the flat base change of the usual cotangent complexes, we see $\rho$ is a quasi-isomorphism.
	
	At last, we only need to note that the collection of rings of definition $\{A_0\ra B_0\}$ such that $B_0=A_0/\xi$ is cofinal with the collection of all of the $A_0\ra B_0$ (since any given $B_0$ is a subring of $B$ that is topologically finite type over $\calO_K$, we can pick the generators and lift them to $A$ along the surjection $A\ra B$).
	So we get
	\[
	\underset{\substack{A_0\ra B_0\\ A_0,B_0~open~bounded}}{\colim} \wh\LL_{B_0/A_0}[\frac{1}{p}]\cong \underset{\substack{ A_0~open~bounded\\}}{\colim} \wh\LL_{(A_0/\xi)/A_0}[\frac{1}{p}],
	\]
	and the latter is quasi-isomorphic to $\LL_{S/S'}^\an\otimes_K B$.
	So we are done.
\end{proof}
    In this way, since $E$ is constructed so that the top horizontal and the right vertical triangles in $(\ast)$ are distinguished, we see under the assumption, $E$ is quasi-isomorphic to $0$.
    This allows us to get the section 
    \[
    s_X:\LL_{X/S'}^\an \rra \LL_{S/S'}^\an\otimes_{S} \calO_X,
    \]
    defined as the composition of the $\alpha_X$ and the $\beta_X^{-1}$ in $(\ast)$.
    
    At last, we check the functoriality.
    Consider the map between two lifts
    \[
    \xymatrix{
    	X\ar[r]^f\ar[d] &Y\ar[d]\\
    	X'\ar[r]\ar[d]\ar[r]&Y'\ar[d]\\
    	S'\ar@{=}[r]& S'.}
    \]
    Then since each term in the big diagram $(\ast)$ is functorial with respect to $X\ra X'$, the map of lifts induces a commutative diagram from the $(\ast)$ for $Y$ to the derived direct image of $(\ast)$ for $X$ along $f:X\ra Y$.
    In particular, this implies the commutativity of the following:
    \[
		\xymatrix{
			\LL_{Y/S'}^\an \ar[r]^{\alpha_Y} \ar[d] & \LL_{Y/Y'}^\an  \ar[d]\ar[r]^{\beta_Y^{-1}}& \LL_{S/S'}^\an\otimes \calO_Y\ar[d]\\
			Rf_*\LL_{X/S'}^\an \ar[r]^{Rf_*(\alpha_X)} & Rf_*\LL_{X/X'}^\an\ar[r]^{Rf_*(\beta_X^{-1})~~~~~~~} & Rf_*(\LL_{S/S'}^\an\otimes \calO_X).}
	\]
	So by combining them, we get the map from $s_Y$ to $Rf_*(s_X)$.
	
\end{proof}

At last, we note the following relation between $\LL_{X/S'}^\an$ and $\LL_{X/\bfS}^\an$.
\begin{lemma}\label{trun}
	Let $X$ be a smooth rigid space over $\Spa(K)$, and $\bfS=\Spa(\Ainf[\frac{1}{p}],\Ainf)$ be the $p$-adic complete adic space.
	Then the sequence of maps $X\ra S' \ra \bfS$ induces the quasi-isomorphism
	\[
	\LL_{X/\Ainf[\frac{1}{p}]}^\an \rra \tau^{\geq -1} \LL_{X/(\Ainf[\frac{1}{p}]/\xi^2)}^\an=\tau^{\geq -1} \LL_{X/S'}^\an.
	\]
	This is functorial with respect to $X$.
\end{lemma}
\begin{proof}
	By taking the distinguished triangle for the transitivity, we get
	\[
	\LL_{S'/\bfS}^\an\otimes_{\calO_{S'}} \calO_X\rra \LL_{X/\bfS}^\an \rra \LL_{X/S'}^\an. \tag{$\ast$}
	\]
	Since $S'=\Spa(\Ainf[\frac{1}{p}]/\xi^2)$ is the closed subspace of $\bfS=\Spa(\Ainf)$ that is defined by the regular ideal $(\xi^2)$, we have
	\[
	\LL_{S'/\bfS}^\an\otimes_{\calO_{S'}} \calO_X=(\xi^2)/(\xi^4)\otimes^L_{\calO_{S'}}\calO_X[1].
	\]
	But note that by the distinguished triangle for $X\ra S \ra \bfS$, we have
	\[
	\LL_{S/\bfS}^\an\otimes_{\calO_{S}} \calO_X\rra \LL_{X/\bfS}^\an \rra \LL_{X/S}^\an,
	\]
	where $\LL_{S/\bfS}^\an\otimes_{\calO_{S}} \calO_X=(\xi)/(\xi^2)\otimes_{K}\calO_X[1]=\xi/\xi^2\calO_X[1]$, and $\LL_{X/S}^\an=\Omega_{X/K}^1[0]$ by the smoothness assumption (\cite[Theorem 7.2.42]{GR03}).
	In this way, since $\LL_{X/\bfS}^\an$ lives in cohomological degree $-1$ and $0$ and is killed by $\xi^2$, the image of $\LL_{S'/\bfS}^\an\otimes^L_{\calO_{S'}} \calO_X=(\xi^2)/(\xi^4)\otimes^L_{\calO_{S'}}\calO_X[1]$ in $\LL_{X/\bfS}^\an$ is $0$.
	Hence the sequence $(\ast)$ induces the quasi-isomorphism 
	\[
	\LL_{X/\Ainf[\frac{1}{p}]}^\an \rra \tau^{\geq -1} \LL_{X/(\Ainf[\frac{1}{p}]/\xi^2)}^\an,
	\]
	that lives in degree $-1$ and $0$.
	
	At last, note that since those two distinguished triangles are functorial with respect to $X$, so is the quasi-isomorphism.
\end{proof}

\subsection{Degeneracy in the smooth setting}
After the basics around the cotangent complex and the lifting criterion, we are going to show the degeneracy theorem for smooth rigid spaces, assuming the liftable condition to $\Bdr/\xi^2$.
We fix a complete and algebraically closed $p$-adic field $K$ as before.

We first prove a simple result about the cotangent complex over the $\Ainf$.
\begin{proposition}\label{Bdr-Zp0}
	Let $A$ be an $\Ainf$-algebra.
	Then  the following natural map of complete cotangent complexes is a quasi-isomorphism
	\[
	\wh\LL_{A/\Z_p}\rra \wh\LL_{A/\Ainf}.
	\]
\end{proposition}
\begin{proof}
	Consider the sequence of maps
	\[
	\Z_p\rra \Ainf\rra A.
	\]
	By basic properties of the usual cotangent complex of rings, we get a distinguished triangle in $D^-(A)$:
	\[
	\LL_{\Ainf/\Z_p}\otimes_{\Ainf} A\rra \LL_{A/\Ainf}\rra \LL_{A/\Z_p}.
	\]
	Apply the derived $p$-completion, we then get the following distinguished triangle
	\[
	(\LL_{\Ainf/\Z_p}\otimes_{\Ainf} A)^\wedge\rra \wh\LL_{A/\Ainf}\rra \wh\LL_{A/\Z_p}.
	\]
	By the derived Nakayama's lemma and the equality  $	(\LL_{\Ainf/\Z_p}\otimes_{\Ainf} A)^\wedge=(\wh\LL_{\Ainf/\Z_p}\otimes_{\Ainf} A)^\wedge$, 
	it suffices to prove the vanishing of the $\LL_{\Ainf/\Z_p}\otimes_{\Z_p} \Z_p/p$ (thus $\wh\LL_{\Ainf/\Z_p}$).
	Note that $\Ainf=W(\calO_K^\flat)$, where $\calO_K^\flat$ is a perfect ring in characteristic $p$. 
	As a consequence, since the cotangent complex of a perfect ring over $\mathbb{F}_p$ is quasi-isomorphic to zero (\cite{Bha17}, 3.1.6), we get
	\[
	\LL_{\Ainf/\Z_p}\otimes_{\Z_p} \Z_p/p=\LL_{\calO_K^\flat/\FF_p}\cong 0.
	\]
	Hence the $p$-adic completed cotangent complex $\wh\LL_{\Ainf/\Z_p}$ vanishes, and we thus obtain the quasi-isomorphism as in the statement.
	
\end{proof}
\begin{corollary}\label{Bdr-Zp}
	Let $X$ be an adic space over $\Spa(K,\calO_K)$, then the sequence of maps $\Q_p\ra\Ainf[\frac{1}{p}]\ra \calO_X$ induces a functorial quasi-isomorphism between analytic cotangent complexes
	\[
	\LL_{\calO_X/\Q_p}^\an\cong \LL_{\calO_X/\Ainf[\frac{1}{p}]}^\an.
	\]
\end{corollary}

\paragraph{Cotangent complex and derived direct image}
Now we are able to connect the cotangent complex with the $R\nu_*\wh\calO_X$.
Our first result is about the truncation of $R\nu_*\wh\calO_X$: 
\begin{theorem}\label{id}
	Let $X$ be a smooth rigid space over $\Spa(K)$.
	Then there exists a functorial quasi-isomorphism in the derived category of $\calO_X$-modules:
	\[
	\LL_{\calO_X/\Ainf[\frac{1}{p}]}^\an(-1)[-1]\cong \tau^{\leq 1}R\nu_*\wh\calO_X,
	\]
	
\end{theorem}
\begin{proof}
	In order to construct the isomorphism above, we will need the analytic cotangent complex for the complete pro-\'etale structure sheaf $\LL_{\wh\calO_X/R}^\an$, where $(R,R^+)$ are either $(\Ainf[\frac{1}{p}],\Ainf)$ or $(\Q_p,\Z_p)$.
	We will first work at the presheaf level and do the construction for affinoid perfectoid rings, and then show that the cotangent complex is in fact a twist of the complete structure sheaf.
	
\begin{itemize}
	\item[Step 1](Calculation at affinoid perfectoid)

	Denote by $X_\ind$ the indiscrete site on the category of affinoid perfectoid objects in $X_\pe$.
	Namely the category $X_{\ind}$ is the collection of affinoid perfectoid objects in $X_\pe$, and the topology is the one such that every presheaf on $X_{\ind}$ is a sheaf.
	Then there exists a canonical map of sites $\delta:X_{\pe}\ra X_\ind$. 
	We note that the inverse image functor $\delta^{-1}$ is an exact functor on abelian sheaves defined by the sheafification, and we have $L\delta^{-1}=\delta^{-1}$.
	
	Then we can define the completed structure sheaf $\wh\calO_\ind^+$, such that for $U\in X_\ind$ with its underlying perfectoid space $\Spa(A,A^+)$, we have
	\[
	\wh\calO_\ind^+(U)=A^+, \wh\calO_\ind(U)=A.
	\]
	Similarly we can define the cotangent complex as 
	\[
	\LL_{\ind,\wh\calO_X^+/R^+}^\an(U)=\varinjlim_{\substack{A_0\subset A^+\\A_0~open~bounded}}\wh\LL_{A_0/R^+}, ~~~\LL_{\ind,\wh\calO_X/R}^\an(U)=\LL_{\ind, \wh\calO_X^+/R^+}^\an[\frac{1}{p}].
	\]
	Here the cotangent complex for formal rings (adic rings) are the one introduced at the beginning of the subsection.
	We note that by the fact that a perfectoid algebra $(A,A^+)$ is uniform, we know $A^\circ$ is bounded in $A$  (\cite[1.6]{Sch12A}).
	In particular, the open subring $A^+$ of $A^\circ$ is also open bounded, and as complexes we have the equality
	\[
	\LL_{\ind, \wh\calO_X^+/R^+}^\an(U)=\wh\LL_{A^+/R^+},~~~\LL_{\ind, \wh\calO_X/R}^\an(U)=\LL_{A/R}^\an.
	\]
	
	We also note that by the Proposition \ref{Bdr-Zp0}, the sequence of sheaves $\Z_p\ra \Ainf\ra \wh\calO_\ind^+$ induces a quasi-isomorphism
	\[
	\LL_{\ind,\wh\calO_X/\Ainf[\frac{1}{p}]}^\an\cong \LL_{\ind, \wh\calO_X/\Q_p}^\an.
	\]
	Here to check this quasi-isomorphism it suffices to check sections at $U\in X_\ind$, since $X_\ind$ has only trivial coverings.

	Moreover, the map of rings $\Z_p\ra \calO_K\ra \wh\calO_\ind^+$ provides the natural distinguished triangle
	\[
	\wh\LL_{\calO_K/\Z_p}\wh\otimes_{\calO_K} \wh\calO_\ind^+ \ra \wh\LL_{\ind,\wh\calO_X^+/\Z_p}\ra \wh\LL_{\ind,\wh\calO_X^+/\calO_K}.
	\]
	Since the mod $p$ reduction of $\calO_K\ra \wh\calO_\ind^+$ is relatively perfect, and $\wh\LL_{\calO_K/\Z_p}$ is isomorphic to the Breuil-Kisin twist $\calO_K\{1\}[1]$ of weight $-1$, we have the quasi-isomorphism
	\[
	\wh\LL_{\ind,\wh\calO_X^+/\Z_p}\cong \wh\calO_\ind^+\{1\}[1].
	\]
	So by inverting $p$, we have
	\[
	\LL_{\ind, \wh\calO_X/\Q_p}^\an\cong\wh\calO_\ind(1)[1]. \tag{$\ast$}
	\]
	Here the same is true when we replace $\Q_p$ by $\Ainf[\frac{1}{p}]$.
	
\item[Step 2](pro-\'etale cotangent complex)

	Now we go back to the pro-\'etale topology.
	As above, let $(R,R^+)$ be either $(\Ainf[\frac{1}{p}],\Ainf)$ or $(\Q_p,\Z_p)$. 
	We first observe that the definition of (integral) analytic cotangent complex can be extended to the whole pro-\'etale site $X_\PE$, which is the complex of sheaves given by sheafifying the complex of presheaves that assigns each object $U$ with its underlying perfectoid space $\Spa(A,A^+)$ to
	\[
	\varinjlim_{\substack{A_0\subset A^+\\A_0~open~bounded}}\wh\LL_{A_0/R^+}, ~~\varinjlim_{\substack{A_0\subset A^+\\A_0~open~bounded}}\wh\LL_{A_0/R^+}[\frac{1}{p}],
	\]
	We denote those two as 
	\[
	\LL_{\wh\calO_X^+/R^+}^\an,~~\LL_{\wh\calO_X/R}^\an.
	\]
	Here we note that the definition is compatible with the one for rigid spaces (see the discussion at the beginning of the subsection \ref{cot-def}).
	In particular by the functoriality of the construction, the canonical map of ringed sites $\nu:X_\PE\ra X_\et$ induces a natural map
	\[
	\LL_{X/R}^\an \rra R\nu_*\LL_{\wh\calO_X/R}^\an.
	\]
	Moreover, as the collection of  affinoid perfectoid open subsets form a base for $X_\PE$, the pro-\'etale cotangent complex is equal to the inverse image of indiscrete cotangent complex along $\delta:X_\PE\ra X_\ind$, i.e.
	\[
	\LL_{\wh\calO_X/R}^\an=\delta^{-1}\LL_{\ind,\wh\calO_X/R}^\an.
	\]
	
	Now we take the (derived) inverse image $\delta^{-1}$ for the quasi-isomorphism $(\ast)$ to get the quasi-isomorphism
	\[
	\LL_{\wh\calO_X/R}^\an\cong \delta^{-1}\wh\calO_\ind (1)[1].
	\]
	On the other hand, for affinoid perfectoid $U$ (with its underlying perfectoid space $\Spa(A,A^+)$) in $X_\PE$, the section of complete structure sheaves at $U$ is known to be (\cite[4.10]{Sch13}) 
	\[
	\wh\calO_X^+(U)=A^+,~~\wh\calO_X(U)=A.
	\]
	In this way, the inverse image of indiscrete structure sheaves are identified with the complete structure sheaves over the pro-\'etale site
	\[
	\delta^{-1}\wh\calO_\ind^+=\wh\calO_X^+,~\delta^{-1}\wh\calO_\ind=\wh\calO_X.
	\]
	Thus  we get the natural quasi-isomorphism
	\[
	\LL_{\wh\calO_X/R}^\an\cong\wh\calO_X(1)[1]. \tag{$\ast\ast$}
	\]
	
\item[Step 3](Comparison)

	At last we consider  the statement in the Theorem.
	The map between ringed sites $\nu:(X_\PE,\wh\calO_X)\ra(X,\calO_X)$ induces a morphism of cotangent complexes
	\[
	\LL_{\calO_X/\Ainf}^\an\rra  R\nu_*\LL_{\wh\calO_X/\Ainf}^\an.
	\]
	By the Corollary \ref{Bdr-Zp} and the natural quasi-isomorphism $\LL_{\wh\calO_X/\Ainf}^\an=\wh\calO_X(1)[1]$ in the Step 2, the map above is isomorphic to the following
	\[
	\LL_{\calO_X/\Q_p}^\an\rra R\nu_*\wh\calO_X(1)[1].
	\]
	So in order to show the quasi-isomorphism in the Theorem \ref{id}, it suffices to show the quasi-isomorphism of the the map
	\[
	\LL_{\calO_X/\Q_p}^\an(-1)[-1]\rra \tau^{\leq 1}R\nu_*\wh\calO_X. \tag{$\ast\ast\ast$}
	\]

	Now, since the statement is local on $X$, we may assume $X$ is affinoid, admitting an \'etale morphism to $\TT^n$.
	Then we note that both sides of the above are invariant under the \'etale base change:
	the right side is a complex of \'etale coherent sheaf, while the base change of the left side is given by the vanishing of the relative cotangent complex for an \'etale map. 
	\footnote{
		This follows from the distinguished triangle $Lf^*\LL_{\TT^n/\mathbb{Q}_p}^\an \ra \LL_{X/\QQ_p}^\an \ra \LL_{X/\TT^n}^\an$ and the vanishing of $\LL_{X/\TT^n}^\an$ (\cite[Theorem 7.2.42]{GR03}), where $f:X\ra \TT^n$ is an \'etale morphism.
		Here we note that as neither $X$ or $\TT^n$ is topologically of finite type over $\mathbb{Q}_p$, we cannot apply $\cite{GR03}$ to get the triangle directly.
		To see the triangle, we first notice that as $X=\Spa(B,B^\circ)$ and $\TT=\Spa(A,A^\circ)$ are reduced and topologically of finite type over $K$, by Remark \ref{simplified cot} the analytic cotangent complex can be naturally computed as follows
		\begin{align*}
		\LL_{X/\mathbb{Q}_p}^\an=\wh\LL_{B^\circ/\mathbb{Z}_p}[\frac{1}{p}],~~
		\LL_{\TT^n/\mathbb{Q}_p}^\an=\wh\LL_{A^\circ/\mathbb{Z}_p}[\frac{1}{p}],~~
		\LL_{X/\TT^n}^\an =\wh \LL_{B^\circ/A^\circ}[\frac{1}{p}].
		\end{align*}
		Moreover, the pullback $Lf^*\LL_{\TT^n/\mathbb{Q}_p}^\an$, which is equal to $B\otimes_A^L (\wh\LL_{A^\circ/\mathbb{Z}_p} [\frac{1}{p}])$ is naturally isomorphic to $(B^\circ\otimes_{A^\circ}^L \LL_{A^\circ/\mathbb{Z}_p})^\wedge[\frac{1}{p}]$. 
		So the distinguished triangle we want can be given by taking the derived $p$-completion and then inverting $p$ at the distinguished triangle for the algebraic cotangent complex of $\mathbb{Z}_p\ra A^\circ \ra B^\circ$.
	}
	So it suffices to show the case when $X=\TT^n=\Spa(K\langle T_i^{\pm1}\rangle,\calO_K\langle T_i^{\pm1}\rangle)$.
	But notice that the map $(\ast\ast\ast)$ can be given by inverting $p$ at the sequence
	\[
	\wh\LL_{\calO_K\langle T_i^{\pm1}\rangle/\Z_p}\{-1\}[-1]\rra \tau^{\leq 1}R\nu_*\wh\calO_{\TT^n}^+,
	\]
	here $\{-1\}$ is the Breuil-Kisin twist of the weight $1$.
	In this way, by the local computation in \cite[8.15]{BMS}, the map above induces a quasi-isomorphism
	\[
	\wh\LL_{\calO_K\langle T_i^{\pm1}\rangle/\Z_p}\{-1\}[-1]\rra \tau^{\leq 1}L\eta_{\zeta_p-1}R\nu_*\wh\calO_{\TT^n}^+,
	\]
	which after inverting $p$ induces the quasi-isomorphism of analytic cotangent complexes
	\[
	\LL_{\TT^n/\Q_p}^\an(-1)[-1]\rra \tau^{\leq 1}R\nu_*\wh\calO_{\TT^n}.
	\]
	Hence we are done.
\end{itemize}
\end{proof}

\begin{corollary}\label{liftsm}
	Assume $X$ is a smooth rigid space over $K$ that admits a flat lift $X'$ along $\Bdr/\xi^2\ra K$.
	Then the lift $X'$ induces a splitting of $\tau^{\leq 1}R\nu_*\wh\calO_X$ into a direct sum of its cohomology sheaves in the derived category. 
	
	Moreover, the splitting is functorial with respect to the lift $X'$.
\end{corollary}
\begin{proof}
	By the above Theorem \ref{id}, we have the functorial quasi-isomorphism
	\[
	\tau^{\leq 1}R\nu_*\wh\calO_X=\LL_{X/\Ainf[\frac{1}{p}]}^\an(-1)[-1].
	\]
	Moreover, the Lemma \ref{trun} about the truncation provides us with a functorial quasi-isomorphism
	\begin{align*}
	\LL_{X/\Ainf[\frac{1}{p}]}^\an(-1)[-1]=(\tau^{\geq -1}(\LL_{X/(\Ainf[\frac{1}{p}]/\xi^2)}^\an))(-1)[-1].
	\end{align*}
	At last, by the Proposition \ref{ob}, the right side splits into the direct sum of its cohomology sheaves if $X$ can be lifted into a flat adic space over $S'=\Spa(\Ainf[\frac{1}{p}]/\xi^2,\Ainf/\xi^2)=\Spa(\Bdr/\xi^2,\Ainf/\xi^2)$,
	such that the splitting in Proposition \ref{ob} is functorial with respect to the lift.
	So we get the result.

\end{proof}
We then notice that the splitting of the derived direct image is in fact true without the truncation.
\begin{proposition}\label{splitsm}
	Assume $X$ is a smooth rigid space over $K$ that admits a flat lift $X'$ over $\Bdr/\xi^2$.
	Then the lift $X'$ induces the derived direct image $R\nu_*\wh\calO_X$ to split as $\bigoplus_{i\geq 0} \Omega_{X/K}^i(-i)[-i]$ in the derived category.
	
	Here the isomorphism is functorial with respect to lifts, in the sense that when $f':X'\ra Y'$ is an $\Bdr/\xi^2$ morphism between lifts of two smooth rigid spaces $f:X\ra Y$ over $K$, then the induced map $R\nu_*\wh\calO_Y\ra Rf_*R\nu_*\wh\calO_X$ is compatible with the map between the direct sum of differentials.
\end{proposition}
\begin{proof}
	By the Corollary \ref{liftsm} above, the given lift to $\Bdr/\xi^2$ induces an $\calO_X$-linear quasi-isomorphism
	\[
	\calO_X[0]\oplus \Omega_{X/K}^1(-1)[-1]\rra \tau^{\leq 1}R\nu_*\wh\calO_X.
	\]
	It is functorial in the sense that if $f':X'\ra Y'$ is $\Bdr/\xi^2$-morphism between lifts of a map of two smooth rigid spaces $f:X\ra Y$ over $K$, then the induced map
	\[
	\tau^{\leq 1}R\nu_*\wh\calO_Y\rra Rf_* \tau^{\leq 1}R\nu_*\wh\calO_X
	\]
	is compatible with the section map
	\[
	\xymatrix{\calO_Y[0] \ar[r]  &Rf_*\calO_X[0] \\
	 \tau^{\leq 1}R\nu_*\wh\calO_Y \ar[u]_{s_Y}\ar[r] & Rf^*\tau^{\leq 1}R\nu_*\wh\calO_X \ar[u]_{Rf_*(s_X)},}
	\]
	which are induced by the functoriality in the Proposition \ref{ob}, Lemma \ref{trun}, and the Theorem \ref{id}.
	
	We compose the decomposition with $\tau^{\leq 1}R\nu_*\wh\calO_X\ra R\nu_*\wh\calO_X$, and get
	\[
	\Omega_{X/K}^1(-1)[-1]\rra R\nu_*\wh\calO_X.
	\]
	Here $R\nu_*\wh\calO_X$ is a commutative algebra object in the derived category $D(\calO_X)$.
	Moreover, as in \cite{DI}, the above map can be lifted to a canonical map of commutative algebra objects in the derived category
	\[
	\bigoplus_{i\geq 0}\Omega^i_{X/K}(-i)[-i]\rra R\nu_*\wh\calO_X. \tag{$\ast$}
	\]
	This can be constructed as follows:
	For each $i\geq 1$, the quotient map $(\Omega_{X/K}^1)^{\otimes i}\ra \Omega_{X/K}^i$ admits a canonical $\calO_X$-linear section $s_i$, by
	\[
	\omega_1\wedge\cdots \wedge\omega_i\lmt \frac{1}{n!}\sum_{\sigma\in S_i} \sgn(\sigma)\omega_{\sigma(1)}\otimes\cdots\otimes\omega_{\sigma (i)}.
	\]
	This allows us to give a canonical $\calO_X$-linear map from $\Omega^i_{X/K}(-i)[-i]$ to $R\nu_*\wh\calO_X$, by the diagram
	\[
	\xymatrix{(\Omega_{X/K}^1(-1)[-1])^{\otimes_{\calO_X}^L i} \ar[r] & (R\nu_*\wh\calO_X)^{\otimes_{\calO_X}^L i} \ar[d]\\
		\Omega_{X/K}^i(-i)[-i] \ar[u]^{s_i} \ar@{.>}[r]& R\nu_*\wh\calO_X.}
	\]
	Here the right vertical map is the multiplication induced from that of $\wh\calO_X$. 
	We note that since $X$ is smooth over $K$, the derived tensor product of $\Omega_{X/K}^1$ over $\calO_X$ degenerates into the usual tensor product.
	Moreover, by construction the total map $\oplus_i s_i$ is multiplicative under the wedge products.
	
	Finally, it suffices to show that the isomorphism for the truncation $\tau^{\leq 1}$ can be extended to the map $(\ast)$ above.
	When $X$ is of dimension one, since $\Omega_{X/K}^i$ is zero for $i\geq 2$, we are done.
	In general, by taking open subsets if necessary we may assume $X$ is affinoid and admits an \'etale map from a closed unit disc $\mathbb{B}^n_K$ over $K$.
	Here $\mathbb{B}^n_K$ admits a natural lift $\mathbb{B}^n_{\Bdr/\xi^2}$, and by its smoothness over $\Bdr/\xi^2$ we can lift $\mathbb{B}^n_K \ra X$ to an \'etale morphism $\mathbb{B}^n_{\Bdr/\xi^2} \ra X'$.
	Then as the map $(\ast)$ commutes with the \'etale base change, it suffices to show this for the closed unit disc $\mathbb{B}^n_K$ with its canonical lift $\mathbb{B}^n_{\Bdr/\xi^2}$.
	Finally, notice that both sides of $(\ast)$ admit K\"unneth formula for products (the K\"unneth formula  for pro-\'etale cohomology can be found in \cite[Proposition 8.14]{BMS}), we thus deduce the isomorphism of $(\ast)$ from the curve case.
	
\end{proof}

\subsection{Simplicial generalizations}\label{sec7.3}
We now generalize results in the past two subsections to simplicial cases.

\paragraph{Simplicial sites and cohomology}
First we recall briefly the simplicial sites. 
The general discussion can be found in \cite[Chapter 09VI]{Sta}.\footnote{We want to mention that the discussion of the simplicial sites and cohomology in our article might become simpler if we use the language of the infinity categories, as the latter behaves better than the ordinary derived category when we consider a diagram of derived objects. }

Consider a non-augemented simplicial object of sites $\{Y_n\}.$
Namely for each nondecreasing map $\phi:[n]\ra [l]$ in $\Delta$, where $[n]$ (resp. $[l]$) is the totally ordered set of $n+1$ (resp. $l+1$) elements, there exists a morphism of sites $u_\phi:Y_l\ra Y_n$ satisfying the commutativity of diagrams induced from $\Delta$.
Then we can define its associated \emph{non-augmented simplicial site} $Y_\bullet$, following the definition of $\calC_{total}$ in \cite[Tag 09WC]{Sta}.
An object of $Y_\bullet$ is defined as an object $U_n \in Y_n$ for some $n\in \NN$, and a morphism $(\phi,f):U_l\ra V_n$ is given by a map $\phi:[n]\ra [l]$ together with a map of objects  $f:U_l\ra u_\phi^{-1}(V_n)$ in $Y_l$.
To give a covering of $U\in Y_n$, it means that to specify a collection of $V_i \in Y_n$, such that $\{V_i\ra U\}$ is a covering in the site $Y_n$.
It can be checked that the definition satisfies axioms of being a Grothedieck topology.
Moreover, by allowing $n$ to include the number $-1$, we can define the \emph{augmented simplicial site} $Y_\bullet$.
Here we remark that unless mentioned specifically, a simplicial site or a simplicial object in our article is always assumed to be non-augmented.
Similarly, by replacing $\Delta$ by the finite category $\Delta_{\leq m}$ and assume $n\leq m$, we get the definition of the \emph{ $m$-truncated simplicial site} $Y_\bullet$.

From the definition above, in order to give a (pre)sheaf on $Y_\bullet$, it is equivalent to give a collection of (pre)sheaves $\calF^n$ on each $Y_n$  together with the data that for any map $\phi:[n]\ra [l]$ in the index category, we specify a map of sheaves $\calF^n\ra u_{\phi*} \calF^l$ over $Y_n$ that is compatible with arrows in $\Delta$.
This allows us to define the derived category $D(Y_\bullet)$ of abelian sheaves on $Y_\bullet$.

We could also define the concept of simplicial ringed sites, which consists of pairs $(Y_\bullet,\calO_{Y_\bullet})$, for $Y_\bullet$ being a simplicial site and $\calO_{Y_\bullet}$ being a sheaf of rings on $Y_\bullet$, assuming $u_\phi:(Y_l,\calO_{Y_l})\ra (Y_n,\calO_{Y_n})$ being maps of ringed sites.

\begin{remark}
	In the level of the derived category, the category $D(Y_\bullet)$ is not equivalent to the category where objects are given by specifying one in each $D(Y_n)$ together with natural boundary maps, unless we replace derived categories by derived infinity categories and also consider the higher morphisms.
	This is the main reason why we need to reconstruct many objects in the simplicial level in this section, instead of using the known results for single site or space directly.
	The essential difference is that  an object in the simplicial level has much stronger functoriality than a collection of objects over each individual space.
\end{remark}

From the construction above, it is clear that there exists a map of sites $Y_\bullet \ra Y_n$.
The pushforward functor along this map is restriction functor, sending the collection of sheaves $(\calF_l)_l$ to its $n$-th component $\calF_n$ and is exact (\cite[Tag 09WG]{Sta}) for sheaves of abelian groups.
Here is a useful Lemma about the vanishing criterion of objects in the derived category of a simplicial site $D^+(Y_\bullet)$.
\begin{lemma}\label{sim-van}
	Let $K$ be an object in the derived category $D^+(Y_\bullet)$ of a $m$-truncated simplicial site $Y_\bullet$ for $m\in \NN\cup\{\infty\}$.
	Then $K$ is acyclic if and only if for each $n\leq m$, the restriction $K|_{Y_n}$ in $D^+(Y_n)$ is acyclic.
\end{lemma}
\begin{proof}
	If $K$ is acyclic, then since the restriction functor is an exact functor, we see $K|_{Y_n}$ is also acyclic.
	
	Conversely, assume $K|_{Y_n}$ is acyclic for each integer $n\leq m$.
	If $K$ is not acyclic, then by the assumption that $K$ lives in $D^+(Y_\bullet)$, we may assume $i$ is the least integer such that the $i$-th cohomology sheaf $\calH^i(K)\in \Ab(Y_\bullet)$ is nonvanishing.
	Then by definition, there exists some $n$ such that $\calH^i(K)|_{Y_n}$ is nonzero.
	But again by the exactness of the restriction, we have the equality
	\[
	\calH^i(K)|_{Y_n}=\calH^i(K|_{Y_n}),
	\]
	where the latter is zero by assumption.
	So we get a contradiction, and hence $K$ is acyclic.
\end{proof}
As a small upshot, we have
\begin{lemma}\label{sim-des}
	Let $\lambda_\bullet:X_\bullet\ra Y_\bullet$ be a morphism of two m-truncated simplicial sites for $m\in \NN\cup\{\infty\}$, such that for each integer $n\leq m$, the map $\lambda_n:X_n\ra Y_n$ is of cohomological descent.
	Namely the canonical map by the adjoint pair
	\[
	\calF\rra R\lambda_{n*}\lambda_n^{-1} \calF
	\]
	is a quasi-isomorphism for any $\calF\in \Ab(Y_n)$.
	Then $\lambda_\bullet$ is also of cohomological descent, namely for any abelian sheaf $\calF^\bullet$ on $Y_\bullet$, the counit map of this adjoint pair is a quasi-isomorphism
	\[
	\calF^\bullet\rra R\lambda_{\bullet*}\lambda_\bullet^{-1}\calF^\bullet.
	\]
\end{lemma}
\begin{proof}
Let $\calC$ be a cone of the map $\calF^\bullet\ra R\lambda_{\bullet*}\lambda_\bullet^{-1}\calF^\bullet$.
	It suffices to show the vanishing of the cone in the derived category $D(Y_\bullet)$.
	Then by the exactness of the restriction functor, for any integer $n\leq m$ the image $C|_{Y_n}$ in $D^+(Y_n)$ is also a cone of 
		\[
		\calF^n\rra (R\lambda_{\bullet*}\lambda_\bullet^{-1}\calF^\bullet)|_{Y_n}=R\lambda_{n*}\lambda_n^{-1}\calF^n,
		\]
		which vanishes by assumption.
	Since both $\calF^\bullet$ and $R\lambda_{\bullet*}\lambda_\bullet^{-1}\calF^\bullet$ are lower bounded, the cone $\calC$ is also in $D^+(Y_\bullet)$ and we can use the Lemma \ref{sim-van} above.
	So we get the result.
	
\end{proof}

\paragraph{Derived direct image for smooth simplicial spaces}
Next, we use simplicial tools above to generalize results of cotangent complexes and derived direct image to their simplicial versions.

Assume $f:X_\bullet\ra Y_\bullet$ is a morphism of (m-truncated) simplicial quasi-compact adic spaces over a $p$-adic Huber pair. 
Then we could define the $simplicial$ $analytic$ $cotangent$ $complex$ $\LL_{X_\bullet/Y_\bullet}^\an$ as an actual complex of sheaves on the simplicial site $X_\bullet$ such that the $n$-th term on the adic space $X_n$ is the analytic cotangent complex $\LL_{X_n/Y_n}^\an$, defined as in Subsection \ref{cot-def}.
In our applications, we will always assume $Y_\bullet$ to be a constant simplicial spaces associated to $T=\Spa(R,R^+)$, for some $p$-adic Huber pair $(R,R^+)$. 
We will use the notation $\LL_{X_\bullet/R}^\an$ or $\LL_{X_\bullet/T}^\an$ to indicate when the case is constant.

Here we emphasize that as in the definition of the analytic cotangent complex for $X_n/Y_n$ above, the complex $\LL_{X_\bullet/Y_\bullet}^\an$ is actual, namely it is defined in the category of complexes of abelian sheaves on $X_\bullet$, not just an object in the derived category.

Now let $X_\bullet$ be a ($m$-truncated) simplicial rigid space over $\Spa(K)$.
Then we can form the cotangent complex $\LL_{X_\bullet/\Ainf[\frac{1}{p}]}^\an$ over $\Ainf[\frac{1}{p}]$.
Moreover we can define the \emph{simplicial differential sheaf} $\Omega_{X_\bullet/K}^i$ on $X_\bullet$ in the way that on each $X_n$, the component of the sheaf is $\Omega_{X_n/K}^i$.

We first generalize the result about the obstruction of the lifting to the simplicial case:
\begin{proposition}\label{sim-lift}
	Let $X_\bullet$ be a ($m$-truncated) smooth quasi-compact simplicial rigid spaces over $\Spa(K)$.
	Then a flat lift $X_\bullet'$ of $X_\bullet$ along $\Bdr/\xi^2\ra K$ induces a splitting of $\LL_{X_\bullet/\Ainf[\frac{1}{p}]}^\an$ into the direct sum of its cohomological sheaves $\calO_{X_\bullet}(1)[1]\bigoplus \Omega_{X_\bullet/K}^1[0]$ in the derived category.
	
	The quasi-isomorphism is functorial with respect to $X_\bullet'$.
\end{proposition}
\begin{proof}
	This is the combination of the Proposition \ref{ob} and the Lemma \ref{trun}.
	We first notice that the sequence of maps
	\[
	X_\bullet\rra S'_\bullet=\Spa(\Ainf[\frac{1}{p}]/\xi^2)_\bullet\rra \bfS_\bullet=\Spa(\Ainf[\frac{1}{p}])_\bullet
	\]
	induces a map 
	\[
	\LL_{X_\bullet/\Ainf[\frac{1}{p}]}^\an \rra \LL_{X_\bullet/S'}^\an.
	\]
	But by the proof of the Lemma \ref{trun} and the vanishing Lemma \ref{sim-van}, the complex $\LL_{X_\bullet/\Ainf[\frac{1}{p}]}^\an$ lives only in degree $-1$ and $0$, which is isomorphic to the truncation of $\LL_{X_\bullet/S'}^\an$ at $\tau^{\geq -1}$.
	So we reduce to consider the splitting of $\LL_{X_\bullet/S'}^\an$.
	
	Now assume $X_\bullet$ admits a flat lift $X_\bullet'$ over $S'$
	The lift leads to the cartesian diagram
	\[
	\xymatrix{
		X_\bullet\ar[r] \ar[d] &X_\bullet '\ar[d]\\
		S_\bullet \ar[r] & S_\bullet ',}
	\]
	which induces the simplicial version of the diagram $(\ast)$ as in the proof of the Proposition \ref{ob}:
	\[
		\xymatrix{
		\LL_{X'_\bullet/S'}^\an\otimes_{\calO_{X'_\bullet}} \calO_{X_\bullet} \ar[r]	& \LL_{X_\bullet/S}^\an \ar[r] & E\\
		\LL_{X'_\bullet/S'}^\an\otimes_{\calO_{X'_\bullet}} \calO_{X_\bullet} \ar[r] \ar@{=}[u]& \LL_{X_\bullet/S'}^\an \ar[r]^{\alpha_\bullet}\ar[u] & L_{X_\bullet/X'_\bullet}^\an \ar[u] \\		
		& \LL_{S/S'}^\an\otimes_{\calO_{S}} \calO_{X_\bullet} \ar[u] \ar@{=}[r] &\LL_{S/S'}^\an\otimes_{\calO_{S}} \calO_{X_\bullet} \ar[u]_{\beta_\bullet},} 
	\]
	The vanishing of $E$ comes down to the vanishing of $E|_{X_n}$ by the Lemma \ref{sim-van}, which is true by assumption and the Proposition \ref{ob}.
	So we get a section map $\beta_\bullet^{-1}\circ \alpha_\bullet$, which splits $\LL_{X_\bullet/S'}^\an$ into the direct sum of $\LL_{X_\bullet/S}^\an$ and $\LL_{S/S'}^\an\otimes_{\calO_{S}} \calO_{X_\bullet}$ in the derived category.
	Note that since $X_\bullet$ is smooth, the cotangent complex $\LL_{X_\bullet/S}^\an$ is $\Omega_{X_\bullet/K}^1[0]$, while the truncation $\tau^{\geq -1}\LL_{S/S'}^\an\otimes_{\calO_{S}} \calO_{X_\bullet}$ is $\calO_{X_\bullet}(1)[1]$.
	Thus we get the result.
	
	At last, the quasi-isomorphism is functorial with respect to $X'_\bullet$, since the big diagram above is functorial with respect to lifts, as in the proof of the Proposition \ref{ob}.
\end{proof}

We then try to connect the simplicial version of cotangent complex with the derived direct image of the completed structure sheaves.

Let $X_\bullet$ be a ($m$-truncated) simplicial quasi-compact rigid spaces over $K$.
Then this induces the following commutative diagram of topoi of simplicial sites, as a simplicial version of the diagram in Section \ref{sec4}.
\[
\xymatrix{
	\Sh(X_{\bullet \PE}) \ar[r]^{\nu_\bullet} & \Sh(X_{\bullet\et}) \\
	\Sh(\Pf_v|_{X_\bullet^\diamond}) \ar[r]_{\alpha_\bullet} \ar[u]^{\lambda_\bullet} & \Sh(X_{\bullet\eh}) \ar[u]_{\pi_\bullet}.}
\]
We then define the complete pro-\'etale structure sheaf $\wh\calO_{X_\bullet}$ on the pro-\'etale simplicial site $X_{\bullet \PE}$, by assigning $\wh\calO_{X_n}$ on the pro-\'etale site $X_{n\, \pe}$.
Similarly we define the untilted complete structure sheaf $\wh\calO_v$ on the site $\Pf_v|_{X^\diamond_\bullet}$.
Here we notice the sheaf $\wh\calO_v$ satisfies the cohomological descent (Paragraph \ref{sec7.4 coh des}) along the canonical map $\lambda_\bullet:\Sh(\Pf_v|_{X^\diamond_\bullet})\ra \Sh(X_{\bullet \PE})$, by the Lemma \ref{sim-des} and comparison results (Proposition \ref{pe-v}).
This leads to the equality
\[
R\nu_{\bullet*}\wh\calO_{X_\bullet}=R\pi_{\bullet*}R\alpha_{\bullet*} \wh\calO_v.
\]
The restriction of this equality on each $X_n$ is the one in Section \ref{sec4}.

Define simplicial $\eh$-differential sheaves $\Omega_{\eh\bullet}^i$ on $X_{\bullet\eh}$ such that on each $X_{n\,\eh}$, the component of the sheaf is $\Omega_{\eh}^i$.
It is by the exactness of the restriction functor and the discussion in Section \ref{sec4} that
\[
R^j\alpha_{\bullet*}\wh\calO_v=\Omega_{\bullet\eh}^j(-j).
\]
When $X_n$ is smooth over $K$ for each $n$, we have
\[
R^j\nu_*\wh\calO_{X_\bullet}=\Omega_{X_\bullet/K}^j(-j),
\]
with
\[
R^i\pi_{\bullet*}R^j\alpha_{\bullet*}=\begin{cases}
0,~i>0;\\
\Omega_{X_\bullet/K}^j(-j),~i=0.
\end{cases}
\]
These are consequences of the $\eh$ differentials for smooth rigid spaces (Theorem \ref{descent}).

\begin{proposition}\label{sim-compare}
	Let $X_\bullet$ be a $m$-truncated simplicial smooth quasi-compact rigid spaces over $\Spa(K)$.
	Then there exists a canonical quasi-isomorphism 
	\[
	\tau^{\leq 1}R\nu_{\bullet*}\wh\calO_{X_\bullet}\cong \LL_{X_\bullet/\Ainf[\frac{1}{p}]}^\an(-1)[-1].
	\]
	The quasi-isomorphism is functorial with respect to $X_\bullet$.
\end{proposition}
\begin{proof}
	The Proposition is the simplicial version of the Theorem \ref{id}.
	
	We first notice that the map of simplicial sites $\nu_\bullet:X_\PE\ra X_{\bullet\et} $ induces a map of analytic cotangent complex
	\[
	\LL_{X_\bullet/\Ainf[\frac{1}{p}]}^\an \rra R\nu_{\bullet*} \LL_{\wh\calO_{X_\bullet}/\Ainf[\frac{1}{p}]}^\an.
	\]
	Meanwhile,  the triple $ \Ainf[\frac{1}{p}]\ra K\ra \wh\calO_{X_\bullet}$ provides us with a distinguished  transitivity triangles  
	\[
	\LL_{K/\Ainf[\frac{1}{p}]}^\an\otimes_K \wh\calO_{X_\bullet} \rra \LL_{\wh\calO_{X_\bullet}/\Ainf[\frac{1}{p}]}^\an\rra \LL_{\wh\calO_{X_\bullet}/K}^\an,
	\]
	where the vanishing of $\LL_{\wh\calO_{X_\bullet}/K}^\an$ follows from the Lemma \ref{sim-van} and the proof of the Step 3 in Theorem \ref{id}.
	So by combining the above two, we get the map
	\[
	\Pi:\LL_{X_\bullet/\Ainf[\frac{1}{p}]}^\an \rra R\nu_{\bullet*} \wh\calO_{X_\bullet}(1)[1].
	\]
	
	Then we consider the induced map of the $i$-th cohomology sheaves $\calH^i$.
	By the exactness of the restriction functor, the restricted map becomes
	\[
	\calH^i(\LL_{X_n/\Ainf[\frac{1}{p}]}^\an) \rra \calH^i(R\nu_{n*} \wh\calO_{X_n}(1)[1]),
	\]
	which is an isomorphism for $i=0,-1$ by the Theorem \ref{id}, and $\calH^i(\LL_{X_n/\Ainf[\frac{1}{p}]}^\an)$ is zero except $i=0,-1$.
	So by the vanishing of the cone, $\Pi$ induces a quasi-isomorphism 
	\[
	\LL_{X_\bullet/\Ainf[\frac{1}{p}]}^\an\rra \tau^{\geq -1}(R\nu_{\bullet*} \wh\calO_{X_\bullet}(1)[1])
	\]
	which leads to the result by a twist.

\end{proof}
Combining the Proposition \ref{sim-compare} and the Proposition \ref{sim-lift}, we get the simplicial version of the splitting for the truncated derived direct image:
\begin{corollary}
	Assume $X_\bullet$ is a ($m$-truncated) smooth quasi-compact  simplicial rigid space over $K$, which admits a flat simplicial lift $X'_\bullet$ to $\Bdr/\xi^2$.
	Then the lift $X'_\bullet$ induces $\tau^{\leq 1}R\nu_{\bullet*}\wh\calO_{X_\bullet}$ to split into the direct sum of its cohomology sheaves  $\calO_{X_\bullet}[0]\bigoplus \Omega_{X_\bullet/K}^1(-1)[-1]$ in $D(X_\bullet)$. 
	
	Here the splitting is functorial with respect to the lift $X_\bullet'$.
\end{corollary}
Moreover, similar to the Proposition \ref{splitsm}, the splitting can be extended without derived truncations.
\begin{corollary}\label{sim-split}
		Assume $X_\bullet$ is a ($m$-truncated) smooth quasi-compact simplicial rigid space over $K$ that admits a flat lift $X_\bullet'$ to $\Bdr/\xi^2$.
	Then the lift $X_\bullet'$ induces the derived direct image $R\nu_{\bullet*}\wh\calO_{X_\bullet}$ to split into 
	\[
	\bigoplus_{i\geq 0} \Omega_{X_\bullet/K}^i(-i)[-i]
	\] in the derived category, which is also isomorphic to 
	\[
	\bigoplus_{i\geq 0}R\pi_{X_\bullet *} ( \Omega_{\bullet\eh}^i(-i)[-i]).
	\]
	
\end{corollary}
\begin{proof}
	By the above Corollary, we have a quasi-isomorphism 
	\[
	\calO_{X_\bullet}[0]\oplus \Omega_{X_\bullet/K}^1(-1)[-1]\rra \tau^{\leq 1} R\nu_{\bullet*}\wh\calO_{X_\bullet}.
	\]
	Then by composing with $\tau^{\leq 1} R\nu_{\bullet*}\wh\calO_{X_\bullet}\ra  R\nu_{\bullet*}\wh\calO_{X_\bullet}$, similar to the proof of the Proposition \ref{splitsm} we may construct the  map below
	\[
	\bigoplus_{i\geq 0} \Omega_{X_\bullet/K}^i(-i)[-i]\rra R\nu_{\bullet*}\wh\calO_{X_\bullet},
	\]
	whose restriction on each $X_n$ is exactly the quasi-isomorphism in the Proposition \ref{splitsm}.
	Thus by the vanishing of the restriction of the cone, we see the above map is a quasi-isomorphism.
	
	At last, notice that by the smoothness of $X_\bullet$, we get the second direct sum expression.
\end{proof}

\subsection{Degeneracy in general}\label{sec7.4}
We then generalize the splitting of the derived direct image to the general case, without assuming the smoothness.
Our main tools are the cohomological descent and the simplicial generalizations in the last subsection. 

\paragraph{Strong liftability}
Before we prove the general degeneracy, we need to introduce a stronger version of the lifting condition, in order to make use of the cohomological descent. 

We first give the definition.
\begin{definition}\label{stronglift}
	Let $X$ be a quasi-compact rigid space over $K$. 
	We say $X$ is $\mathrm{strongly~liftable}$ if for each non-negative integer $n$, there exists an $n$-truncated augmented simplicial map of adic spaces $X'_{\leq n}\ra X'$ over $\Bdr/\xi^2$, where $X'$ and each $X'_i$ are flat and topologically of finite type over $\Bdr/\xi^2$, such that the pullback along $\Bdr/\xi^2\ra K$ induces an $n$-truncated smooth $\eh$-hypercovering of $X$ over $K$.

We call any such augmented $X'_{\leq n}\ra X'$ simplicial rigid space (or $X'_{\leq n}$ in short) a \emph{strong lift of length} $n$.
\end{definition}
\begin{example}\label{discrete lift}
	Let $k=\calO_K/\frakm_K$ be the residue field of $\mathcal{O}_K$, and we fix a section $i:k\ra \calO_K/p$ for the canonical surjection $\calO_K/p \ra k$ (whose existence is guaranteed by the formal smoothness of the perfect field $k$ over $\mathbb{F}_p$ (\cite[Tag 031Z]{Sta})).
	Note that this induces an injection of fields from $W(k)[\frac{1}{p}]$ to $K$ by the universal property of the Witt ring.
	Let $K_0$ be a subfield of $K$ that is finite over $W(k)[\frac{1}{p}]$, and let $X$ be a rigid space defined over $K_0$.
	We then claim that $X$ is strongly liftable.
	
	To see this, we first notice that as the resolution of singularities holds for rigid spaces over $K_0$, it suffices to show that any such field $K_0$ admits a map $K_0\ra \Bdr/\xi^2$ compatible with the inclusion $K_0\ra K$ above.
	Recall the ring $\Ainf$ is defined as $W(\mathcal{O}_{K^\flat})$, where $\mathcal{O}_K^\flat$ is the inverse limit $\underset{x\mapsto x^p}{\varprojlim} \mathcal{O}_K/p$.
	By the construction of $\mathcal{O}_{K^\flat}$ and the functoriality for the inverse limit and for Frobenius maps, the section $i:k \ra \calO_K/p$ induces a homomorphism $k\ra \mathcal{O}_{K^\flat}$, where the latter is a section to the canonical surjection $\calO_{K^\flat}\ra k$.
	In this way, thanks to the functoriality of the Witt vector functor, we can lift the section map to $W(k)\ra \Ainf=W(\mathcal{O}_{K^\flat})$.
	As an upshot, we get the following composition
	\[
	W(k)[\frac{1}{p}] \ra \Ainf[\frac{1}{p}]\ra \Ainf[\frac{1}{p}]/\xi^2=\Bdr/\xi^2,
	\]
	which lifts the map $W(k)[\frac{1}{p}]\ra K$ that we started with.
	At last, note that any finite field extension $K_0$ of $W(k)[\frac{1}{p}]$ is \'etale over $W(k)[\frac{1}{p}]$, while $\Bdr/\xi^2 \ra \mathcal{O}_K$ is an nilpotent extension of $W(k)[\frac{1}{p}]$ algebras.
	Hence $K_0$ admits a map to $\Bdr/\xi^2$ by the \'etaleness.

	Here we also note that this implies the strong liftability when $X$ is defined over a discretely valued subfield $L_0\subset K$ that is of perfect residue subfield $k'\subset k$, since any such $L_0$ is finite over $W(k')[\frac{1}{p}]$, while the latter is contained in $W(k)[\frac{1}{p}]$. 
	
\end{example}
\begin{example}
	Another example is the analytification of a finite type algebraic variety, by the spreading out technique.
	
	Let $Y$ be a finitely presented scheme over $K$.
	By \cite[8.9.1]{EGA-IV}, there exists a finitely generated $\Q$-subalgebra $A$ in $K$ together with a finitely presented $A$-scheme $Y_0$, such that $Y_0\times_{\Spec(A)} \Spec(K)=Y$.
	As the map $A\ra K$ factors through the fraction field of $A$, we may assume $A$ is a finitely generated field extension of $\mathbb{Q}$ and $Y_0$ is defined over $A$.
	Notice that the transcendental degree of $\mathbb{Q}_p$ over $\Q$ is infinite.
	So by embedding a transcendental basis of $A$ over $\mathbb{Q}$ into $\mathbb{Q}_p$, we may find a finite extension $K_0$ of $\mathbb{Q}_p$ such that $A$ can be embedded into $K_0$.
	In this way, we reduce the case to Example \ref{discrete lift}, as $Y^\an$ can be defined over a discrete valud subfield $K_0$ of $K$ that has a perfect residue field.

\end{example}
By the upcoming work of the spreading out of rigid spaces by Conrad-Gabber \cite{CG}, it turns out that $X$ is strongly liftable if it is a proper rigid space over $K$.

\begin{proposition}\label{proper lift}
	Let $X$ be a proper rigid space over $K$.
	Then it is strongly liftable.
\end{proposition}
\begin{proof}
	We follow the proof for the spreading out technique for rigid spaces by Bhatt-Morrow-Scholze in \cite{BMS} and study the structure of the deformation ring.
	However, instead of working on one rigid space, we need to work with a finite diagram of proper rigid spaces.
	Similar to Example \ref{discrete lift}, we fix a section $i:k=\calO_K/\frakm_K \ra \calO_K/p$ to the canonical surjection $\calO_K/p \ra k$, which induces an inclusion of $p$-adic fields $W(k)[\frac{1}{p}]\ra K$.
	
	Let $n$ be any non-negative integer.
	By the resolution of singularity (Theorem \ref{ROS}), we can always construct a (n-truncated) smooth $\eh$-hypercovering $X_{\leq n}\ra X$ over $K$, where each $X_i$ is proper over $K$ (\cite[Section 4]{Con03}).
	Then it suffices to show that there exists a proper $\eh$ hypercovering $X_{\leq n}\ra X$, together with a smooth rigid space $\calS$ over a subfield $K_0=W(k)[\frac{1}{p}]$ of $K$, such that the n-truncated simplicial diagram $X_{\leq n}\ra X$ can be lifted to a diagram of proper $K_0$-rigid spaces $\calX_{\leq n}\ra \calX$ over $\calS$.
	This is because the nilpotent extension $\Bdr/\xi^2\ra K$ is $K_0$ linear, so by the smoothness of $\calS$, the map $\Spa(K)\ra \calS$ can be lifted to a map $\Spa(\Bdr/\xi^2) \ra \calS$.
	Thus the base change of $\calX_{\leq n}$ along this lifting does the job.
	
	Now we prove the statement, imitating the proof of Proposition 13.15 and Corollary 13.16 in \cite{BMS}.
	We first deal with the formal lifting over the integral base.
	Let $W=W(k)$ be the ring of the Witt vector for the residue field of $K$, and $\scrC_W$ be the category of artinian, complete local $W$-rings with the same residue field $k$.
	We first make the following claim:
	\begin{claim}
		There exists an n-truncated smooth $\eh$-hypercovering $X_{\leq n}\ra X$ over $K$, such that it admits a lift to an n-truncated simplicial diagram of $p$-adically complete, topologically finite type $\calO_K$-formal schemes: 
		\[X^+_{\leq  n}\ra X^+.\]
	\end{claim}
	\begin{proof}
		Fix an $\calO_K$-integral model $X^+$ of $X$, whose existence is guaranteed by Raynaud's result on the relation between rigid spaces over $K$ and admissible formal schemes over $\calO_K$.
		We now construct inductively the required covering and the integral lift, following the idea of split hypercoverings (see \cite[Section 4]{Con03}, or \cite[Tag 094J]{Sta} for discussions).
		
		By the local smoothness of the $\eh$-topology (Corollary \ref{localsm}), pick $X_0\ra X$ to be a smooth $\eh$ covering that is proper over $X$.
		By Raynaud's result, there exists a morphism $X_0^+\ra X^+$ of $\calO_K$-formal schemes that lifts the $X_0\ra X$.
		This is the lift of the face map of the simplicial object at the degree $0$.
		
		Assume we already have an $n$-truncated smooth proper $\eh$-hypercovering $X_{\leq n}\ra X$ together with the integral lift $X^+_{\leq n}\ra X^+$ over $\calO_K$.
		Then recall from \cite[4.12, 4.14]{Con03}   that in order to extend $X_{\leq n}\ra X$ to a smooth proper $n+1$-truncated hypercovering whose $n$-truncation is the same as $X_{\leq n}$, it is equivalent to find a smooth proper $\eh$ covering of rigid spaces
		\[
		N\ra (\cosk_n X_{\leq n})_{n+1}.
		\]
		Under the construction, the degree $n+1$-term of the resulting $(n+1)$-truncated hypercovering will be 
		\[
		X_{n+1}:=N\coprod N',
		\]
		for $N'$ being some finite disjoint union of irreducible components of $X_i (0\leq i\leq n)$ (which is also smooth and proper over $K$).
		Such a smooth proper $\eh$-covering $N$ exists by the local smoothness of $X_\eh$.
		Furthermore, while we form this $n+1$-hypercovering of the rigid spaces, we also want to find the integral lift
		\[
		N^+\ra (\cosk_n X^+_{\leq n})_{n+1}
		\]
		of the morphism $N\ra (\cosk_n X_{\leq n})_{n+1}$.
		To do this, we use \cite[4.12]{Con03} and do the same formal construction for $N^+$ and $X^+_{\leq n}\ra X^+$ as above, and extends the latter to an $n+1$-truncated simplical formal schemes 
		\[
		X^+_{\leq n+1}\ra X^+,
		\]
		where $X^+_{n+1}=N^+\coprod N'^+$ is an $\calO_K$-model  of $X_{n+1}$.
		In this way, the generic fiber of $X^+_{\leq n+1}\ra X^+$ is  a $(n+1)$-simplicial object over $X$ whose $n$-truncation is $X_{\leq n}\ra X$, and whose $(n+1)$-th term is $X_{n+1}=N\coprod N'$, which is in fact a smooth proper $\eh$-covering of $(\cosk_n X_{\leq n})_{n+1}$.
		Hence 		by the induction hypothesis we are done.

	\end{proof}
	We then fix such an $\eh$-hypercovering $X_{\leq n}\ra X$ with its integral model $X^+_{\leq n}\ra X^+$ as in the claim.
	Define the functor of deformations of the special fiber $X^+_{\leq n,k}:=X^+_{\leq n}\times_W k$
	\[
	{Def}:\scrC_W\rra \Set,
	\]
	which assigns each $R\in \scrC_W$ to the isomorphism classes of lifts of the digrams $X_{\leq n,k}^+\ra X_k^+$ along $R\ra k$, such that each lifted rigid space is proper and flat over $R$.
	This functor is a deformation functor, and admits a versal deformation: to see this, we first note that as in \cite[Tag 0E3U]{Sta}, the functor $Def$ satisfies the Rim-Schlessinger condition (\cite[Tag 06J2]{Sta}).
	Then we made the following claim
	\begin{claim}
		The tangent space $TDef:=Def(k[\epsilon]/\epsilon^2)$ of the deformation functor  is of finite dimension.
	\end{claim}
	\begin{proof}[Proof of the Claim]
		Notice that there is a natural (forgetful) functor from $Def$ to the deformation functor of the morphism $Def_{Y_s\ra Y_t}$, where $Y_s$ is the disjoint union of  all the sources of arrows in the diagram $X^+_{\leq n,k}$, and $Y_t$ is the disjoint union of all of those targets.
		This induces a map between tangent spaces
		\[
		TDef\rra TDef_{Y_s\ra Y_t}.
		\]
		By the construction, both $Y_s$ and $Y_t$ are finite disjoint unions of proper rigid spaces, which are then proper over $k$.
		So from \cite[Tag 0E3W]{Sta}, we know the tangent space $TDef_{Y_s\ra Y_t}$ is finite dimensional.
		Furthermore, assume $D_1$ and $D_2$ are two lifted diagrams over $k[\epsilon]$.
		Then the difference of $D_1$ and $D_2$ is the collection of $k$-derivations $\calO_{X^+_{t(\alpha), k}}\ra \alpha_*\calO_{X^+_{s(\alpha), k}}$, satisfying certain $k$-linear relation so that those arrows in $D_1$ and $D_2$ commute.
		In particular, this consists of a subspace of $Der_k(\calO_{Y_t},u_*\calO_{Y_s})=\Hom_{Y_t}(\Omega_{Y_t/k}^1,u_*\calO_{Y_s})$, which by the properness again is finite dimensional.
		In this way, both the kernel and the target of the map $TDef\ra TDef_{Y_s\ra Y_t}$ are of finite dimensions, thus so is the $TDef$.
	\end{proof}
	By the above claim and \cite[Tag 06IW]{Sta}, the deformation functor $Def$ admits a versal object.
	In other words, there exists a complete artinian local $W$-algebra $R$ with the residue field $k$, and a diagram $\calX_{R,\leq n}\ra \calX_R$ of proper flat formal $R$ schemes deforming $X^+_{\leq n,k}\ra X^+_k$, such that the induced classifying map
	\[
	h_R:=\Hom_W(R,-)\rra Def
	\]
	is formally smooth.
	Moreover, by the proof of Proposition 13.15 in \cite{BMS}, we can take the ind-completion of $\scrC_W$ and extend $Def$ to a bigger category, which consists of local zero-dimensional $W$-algebras with residue field $k$ (not necessary to be noetherian).
	The category includes $\calO_K/p^m$, and since $X^+_{\leq n}\ra X^+$ is an $\calO_K$-lifting of $X^+_{\leq n,k}\ra X^+_k$, we see the diagram can be obtained by the base change of the universal family $\calX_{R,\leq n}\ra \calX_R$ along $R\ra \calO_K=\varprojlim_m\calO_K/p^m$.
	
	At last, we invert $p$ at the diagram 
	\[
	\xymatrix{
		X^+_{\leq n} \ar[r] \ar[d] & \Spf(\calO_K)\ar[d]\\
		\calX_{R,\leq n}\ar[r]& \Spf(R).}
	\]
	The diagram $X_{\leq n}\ra X$ then can be obtained  from a truncated simplicial diagram $\calX_{R,\leq n}[\frac{1}{p}]\ra \calX[\frac{1}{p}]$ of proper $K_0$-rigid spaces that are flat over $\calS=\Spa(R[\frac{1}{p}],R)$.
	By shrinking $S$ to a suitable locally closed subset if necessary, we may assume $S$ is smooth over $K_0$.
	So we are done.

\end{proof}

\paragraph{Cohomological desent}\label{sec7.4 coh des}
Another preparation we need is the hypercovering and the cohomological descent.

Assume we have a non-augmented simplicial site $Y_\bullet$ (truncated or not) and another site $S$.
Let $\{a_n:Y_n\ra S\}$ for $n\geq 0$ be a collection of morphisms to $S$ such that it is compatible with face maps and degeneracy maps in $Y_\bullet$.
Then we can define the $augmentation$ morphism $a:\Sh(Y_\bullet)\ra \Sh(S)$ between the topoi of $Y_\bullet$ and  $S$, such that for an abelian sheaf $\calF^\bullet$ on $Y_\bullet$, we have
\[
a_*\calF^\bullet=\ker(a_{0 *}\calF^0\rightrightarrows a_{1 *}\calF^1).
\]
It can be checked that the derived direct image $Ra_*$ can be written as the composition
\[
Ra_*=\bfs\circ Ra_{\bullet *},
\]
where $a_\bullet:Y_\bullet\ra S_\bullet$ is the morphism from $Y_\bullet$ to the constant simplicial site $S_\bullet$ associated to $S$, and $\bfs$ is the exact functor that takes a simplicial complex to its associated cochain complex of abelian groups.
Here we call the augmentation $a=\{Y_n\ra S\}$ is of $cohomological~ descent$ if the counit map induced by the adjoint pair $(a^{-1},a_*)$ is a quasi-isomorphism
\[
id\rra Ra_* a^{-1}.
\]

The augmentation allows us to compute the cohomology of sheaves on $S$ by the spectral sequence associated to the simplicial site.
\begin{lemma}[\cite{Sta}, 0D7A]\label{ss}
	Let $Y_\bullet$ be a simplicial site, or a $m$-truncated simplicial site for $m\geq 0$, and let $a=\{a_n:Y_n\ra S\}$ be an augmentation.
	Then for $K\in D^+(Y_\bullet)$, there exists a natural spectral sequence
	\[
	E_1^{p,q}=R^qa_{p *}(K|_{X_p})\Longrightarrow R^{p+q}a_*K,
	\]
	which is functorial with respect to $Y_\bullet\ra S$ and $K$.
	
	Moreover if we assume the $Y_\bullet$ is non-truncated, the augmentation $a$ is of cohomological descent, and $L\in D^+(S)$, then by applying the spectral sequence to $K=a^{-1}L$ we get a natural spectral sequence
	\[
	E_1^{p,q}=R^qa_{p *}a_p^{-1} L\Longrightarrow \calH^{p+q}(L).
	\]
\end{lemma}
We need another variant of this Lemma in order to use truncated hypercoverings to approximate the cohomology of $S$.
\begin{proposition}\label{sim-app}
	Let $\rho:Y_{\leq m}\ra S$ be a $m$-truncated simplicial hypercovering of sites for $m\in \NN$.
	Then for any $\calF\in\Ab(S)$, the cone for the natural adjunction map 
	\[
	\calF\ra R\rho_*\rho^{-1}\calF
	\]
	lives in the cohomological degree $\geq m-1$.

\end{proposition}
\begin{proof}
Let $\wt\rho:\cosk_m Y_{\leq m}\ra S$ be the $m$-th coskeleton of $\rho:Y_{\leq m}\ra S$.
We use the same symbols $\cosk_m Y_{\leq m}$ and $Y_{\leq m}$ to denote their associated simplicial  sites.
Then there exists a natural map of sites
\[
\iota:\cosk_m Y_{\leq m} \rra Y_{\leq m}.
\]
Those two augmentations induce maps of topoi
\[
\wt\rho:\Sh(\cosk_m Y_{\leq m})\ra \Sh(S),~\rho:\Sh(Y_{\leq m})\ra \Sh(S).
\]
By construction, as maps of topoi we have
\[
\wt\rho=\iota\circ\rho.
\]
So from this, for $\calF\in \Ab(S)$, we get the following commutative diagram
\[
\xymatrix{
	\calF\ar[rr] \ar[rd]& &R\wt\rho_*\wt\rho^{-1} \calF\\
	&R\rho_*\rho^{-1}\calF\ar[ru]&.}
\]

Now we let $\calC$ be a cone of $\calF\ra R\rho_*\rho^{-1} \calF$, and let $\wt\calC$ be a cone of $\calF\ra R\wt\rho_*\wt\rho^{-1} \calF$.
Then the above diagram induces the following commutative diagram of long exact sequence
\[
\xymatrix{
	\cdots \ar[r] & \calH^i(\calF) \ar[r] \ar[d] & R^i\wt\rho_*\wt\rho^{-1} \calF \ar[r] \ar[d]& \calH^i(\wt\calC)\ar[d] \ar[r]& \cdots\\
	\cdots \ar[r] & \calH^i(\calF) \ar[r]  & R^i\rho_*\rho^{-1} \calF \ar[r] & \calH^i(\calC) \ar[r]& \cdots.}
\]
By the Lemma \ref{ss} above and the commutative diagram, we have a map of the first pages of spectral sequences
\[
\xymatrix{R^q\wt\rho_{p*}\wt\rho_p^{-1} \calF \ar@{=>}[r]\ar[d]&  R^{p+q}\wt\rho_{*}\wt\rho^{-1}\calF \ar[d]\\
	R^q\rho_{p*}\rho_p^{-1} \calF \ar@{=>}[r] & R^{p+q}\rho_*\rho^{-1}\calF.}
	\]
But note that since $\wt\rho$ is the $m$-coskeleton of $\rho$, when $p+q\leq m$ the formation $R^q\wt\rho_{p*}\wt\rho_p^{-1}$ is the same as $R^q\rho_{p*}\rho_p^{-1}$.
So we get the isomorphism
\[
R^{p+q}\wt\rho_{*}\wt\rho^{-1}\calF \cong R^{p+q}\rho_*\rho^{-1}\calF,~p+q\leq m.
\]
Besides,  since $\rho$ is a $m$-truncated hypercovering, by Deligne the augmentation $\wt\rho$ of the coskeleton satisfies the cohomological descent.
So the map $\calF\ra R\wt\rho_*\wt\rho^{-1} \calF$ is a quasi-isomorphism.
In this way, the map $\calH^i(\calF)\ra R^i\rho_*\rho^{-1} \calF$ is an isomorphism for $i\leq m$, and hence $\calC$ lives in $D^{\geq m-1}(S)$.

\end{proof}

\paragraph{The degeneracy theorem}
Now we are able to state and prove our main theorem about the degeneracy.
\begin{theorem}\label{deg}
	Let $X$ be a quasi-compact, strongly liftable rigid space of dimension $n$ over $K$, and let the augmented truncated simplicial spaces $X'_{\leq m}$ be a strong lift of $X$ of length $m\geq 2n+2$.
	Then the strong lift $X'_{\leq m}$ induces a quasi-isomorphism
	\[
\Pi_{X'_{\leq m}}: R\nu_*\wh\calO_X\rra  \bigoplus_{i=0}^n R\pi_{X *}(\Omega_\eh^i(-i)[-i]).
\]

	The quasi-isomorphism $\Pi_{X'_{\leq m}}$ is functorial among strong lifts $X'_{\leq m}$ of rigid spaces of length $m\geq 2n+2$, in the sense that a map of $m$-truncated strong lifts $X'_{\leq m}\ra Y'_{\leq m}$ of $f:X\ra Y$ will induce the following commutative diagram in the derived category
	\[
	\xymatrix{Rf_*R\nu_*\wh\calO_X \ar[rr]^{Rf_*(\Pi_{X'_{\leq m}})~~~~~~~~~~~~~} && \bigoplus_{i=0}^{\dim(X)} Rf_*R\pi_{X *}(\Omega_\eh^i(-i)[-i]) \\
		R\nu_*\wh\calO_Y \ar[rr]_{\Pi_{Y'_{\leq m}}~~~~~~~~~~~~~~} \ar[u]& & \bigoplus_{i=0}^{\dim(Y)} R\pi_{Y*} (\Omega_\eh^i(-i)[-i])\ar[u],}
	\]
	where the right vertical map is induced by the functoriality of the K\:ahler differential.
		
\end{theorem}

\begin{proof}
	By assumption, we may assume $X_{\leq m}$ is a $m$-truncated smooth proper $\eh$-hypercovering of $X$ that admits a lift $X'_{\leq m}$ to a simplicial flat adic spaces over $\Bdr/\xi^2$.
	Denote by $\rho:X_\bullet\ra X$ the augmentation map.
	Then $X_{\leq m}$ is also an $m$-truncated $v$-hypercovering, and we have a natural map
	\[
	\wh\calO_v\rra R\rho_{v*}\rho^{-1}\wh\calO_v\cong R\rho_{v*}\wh\calO_{\bullet v},
	\]
	whose cone has cohomological degree $m-1\geq 2n+1$ by the Proposition \ref{sim-app}.
	
	We then apply derived direct image functors, and get a natural map
	\begin{align*}
	R\nu_*\wh\calO_X=R\pi_{X *}R\alpha_*\wh\calO_v&\ra R\pi_{X *}R\alpha_*R\rho_{v*}\wh\calO_{\bullet v}\\
	&=R\rho_*R\pi_{\bullet*}R\alpha_{\bullet*}\wh\calO_{\bullet v}\\
	&=R\rho_*R\nu_{\bullet*} \wh\calO_{X_\bullet}.
	\end{align*}
	Here the cone of the map lives in degree $\geq m-1\geq 2n+1$.

	Moreover, by the Corollary \ref{sim-split}, the strong lift $X'_{\leq m}$ induces a functorial (among strong lifts) quasi-isomorphism 
	\[
	R\nu_{\bullet*} \wh\calO_{X_\bullet}\rra  \bigoplus_{i\geq 0} R\pi_{\bullet*}(\Omega_{\bullet \eh}^i(-i)[-i]).
	\]
	So we get the following  distinguished triangle
	\[
	R\nu_*\wh\calO_X\rra R\rho_*R\pi_{\bullet*}(\Omega_{\bullet \eh}^i(-i)[-i])\rra \calC_1, \tag{1}
	\]
	where $\calC_1\in D^{\geq 2n+1} (X)$.
	
	Besides, by the Corollary \ref{sim-split} and the Proposition \ref{sim-app} again the truncated $\eh$-hypercovering $\rho$ induces a natural map
	\[
	\bigoplus_{i\geq 0} R\pi_*(\Omega_{ \eh}^i(-i)[-i])\rra\bigoplus_{i\geq 0} R\pi_*(R\rho_{\eh*}\rho_\eh^{-1}\Omega_{ \eh}^i(-i)[-i])= \bigoplus_{i\geq 0} R\rho_*R\pi_{\bullet*}(\Omega_{\bullet \eh}^i(-i)[-i]), \tag{2}
	\]
	whose cone $\calC_2$ lives in degree $\geq m-1\geq 2n+1$.
	
	At last, by combining $(1)$ and $(2)$, we get the following diagram that is functorial with respect to $X'_{\leq m}$, with both horizontal and vertical being distinguished:
	\[\xymatrix{
		& \calC_2&\\
		R\nu_*\wh\calO_X\ar[r]& R\rho_*R\pi_{\bullet*}(\Omega_{\bullet \eh}^i(-i)[-i])\ar[r]\ar[u]& \calC_1\\
		& \bigoplus_{i\geq 0} R\pi_*(\Omega_{ \eh}^i(-i)[-i]) \ar[u]&}
	\] 
	But note that since $\dim(X)=n$, by the cohomological boundedness (the Corollary \ref{coh-bound} and the Corollary \ref{coh-bound2}), both $R\nu_*\wh\calO_X$ and $\bigoplus_{i\geq 0} R\pi_*(\Omega_{ \eh}^i(-i)[-i])$ live in degree $\leq 2n$.
	Thus by taking the truncation $\tau^{\leq 2n}$, we get the quasi-isomorphism
	\[
	\xymatrix{ R\nu_*\wh\calO_X \ar[r]^{\sim~~~~~~~~~~~~~~} & \tau^{\leq 2n} (R\rho_*R\pi_{\bullet*}(\Omega_{\bullet \eh}^i(-i)[-i]))& \bigoplus_{i\geq 0}^n R\pi_*(\Omega_{ \eh}^i(-i)[-i]) \ar[l]_{~~~~~\sim}.} \tag{3}
	\]
	In this way, by taking $\Pi_{X'_{\leq m}}$ to be the quasi-isomorphism induced from $(3)$, we are done.
	
\end{proof}

\begin{corollary}\label{direct sum}
	Assume $X$ is a quasi-compact rigid space over $K$ that is either defined over a discretely valued subfield $K_0$ of perfect residue field, or proper over $K$.
	Then we have a non-canonical decomposition
	\[
	R\nu_*\wh\calO_X\cong \bigoplus_{i=0}^{\dim(X)} R\pi_{X *}(\Omega_\eh^i(-i)[-i]).
	\]
	In particular, the $\eh$-pro\'et spectral sequence (Theorem \ref{main}) degenerates at the $E_2$-page.
\end{corollary}

\subsection{Finiteness revisited}
In this subsection, we use the degeneracy of the derived direct image $R\nu_*\wh\calO_X$ to improve the cohomological boundedness results in Section \ref{sec6}. 

We first recall the recent work on the perfection and the almost purity in \cite{BS}.
\begin{theorem}[\cite{BS}, Proposition 8.5, Theorem 10.9]\label{apt}
	Let $A$ be a perfectoid ring, $B$ a finitely presented finite $A$-algebra, such that $\Spec(B)\ra \Spec(A)$ is finite \'etale over an open subset.
	Then there exists a perfectoid ring $B_{\pfd}$ together with a map of $A$-algebras $B\ra B_{\pfd}$, such that it is initial among all of the $A$-algebra maps $B\ra B'$ for $B'$ being perfectoid.
\end{theorem}

\begin{proposition}\label{finite rev}
	Let $X$ be a rigid space over $K$.
	Then $R^n\nu_*\wh\calO_X$ vanishes for $n> \dim(X)$.
\end{proposition}
Before we prove the statement, we want to mention that the proof of this Proposition will not need the $\eh$-pro\'et spectral sequence developed above.
\begin{proof}

	Since this is an \'etale local statement, and any \'etale covering of $X$ does not change the dimension, by passing $X$ to its open subsets if necessary, we may assume $X$ admits a finite surjective map onto a torus of the same dimension. \footnote{To see the existence of such surjections, we may argue as follows: as a unit disc is covered by finite many tori, it suffices to find a finite map from an affinoid rigid space $X=\Spa(A,A^+)$ onto a unit disc of the same dimension. Let $A_0$ be a ring of definition of $(A,A^+)$ that is topologically finite type over $\mathcal{O}_K$.
		Since $A_0/\frakm_K$ is a finite type algebra over the residue field $k=\mathcal{O}_K/\frakm_K$, by Noether's normalization lemma we could find a subalgebra $k[x_i]$ of $A_0/\frakm_K$ such that $A_0/\frakm_K$ is finite over $k[x_i]$. In this way, by lifting the map to a morphism $\mathcal{O}_K\langle x_i\rangle \ra A_0$, we get a finite surjective morphism from $X$ to a disc.}  
	
	We give them some notations.
	Denote by $X=\Spa(R,R^+)$ an affinoid rigid space over $\Spa(K)$.
	Assume there exists a finite surjective map $X\ra \TT_n=\Spa(K\langle T_i\rangle, \calO_K\langle T_i\rangle)$ onto the torus of dimension $n$.
	Let $\TT_n^\infty$ be the natural pro-\'etale cover of $\TT_n$ by extracting all $p^n$-th roots of $T_i$, and let $\hat{\TT}_n^\infty=\Spa(K\langle T_i^{\frac{1}{p\infty}} \rangle, \calO_K\langle T_i^{\frac{1}{p^\infty}}\rangle)$ be the underlying affinoid perfectoid space.
	Then the base change of $\TT_n^\infty$ along the map $X\ra \TT_n$ produces a pro-\'etale cover $X^\infty\ra X$ of $X$.
	Note that $\TT_n^\infty \ra \TT_n$ is an $\ZZ_p(1)^n$-torsor, so we have
	\[
	R\Gamma(X_\pe, \wh\calO_X)=R\Gamma_\ct(\ZZ_p(1)^n,R\Gamma_\pe(X^\infty, \wh\calO_X)).
	\]
	Thanks to the (pro-\'etale)-v comparison (Proposition \ref{pe-v}), the above is given by
	\[
	R\Gamma(X_\pe ,\wh\calO_X)=R\Gamma_\ct(\ZZ_p(1)^n, R\Gamma_v(\hat X^{\infty, \diamond},\wh\calO_v)),
	\]
	where $\hat X^{\infty,\diamond}$ is the small $v$-sheaf associated to the analytic adic space $\hat X^{\infty}$ as in the Proposition \ref{ad-dia}.
	Here we note that since $\hat \TT_n^\infty, \TT_n,$ and $X$ are all affinoid, we can write $\hat X^\infty$ as $\Spa(B[\frac{1}{p}],B)$ for some $p$-adic complete $\calO_K$-algebra $B$.
	
	We then recall that for a perfectoid space $Y$ of characteristic $p$ with a structure map to the $v$-sheaf $\Spd(K)$, and any $K$-analytic adic space $Z$, we have the following bijection (cf. \cite[10.2.4]{SW20}):
	\[
	\Hom_{\Spa(K)}(Y^\sharp,Z)=\Hom_{\Spd(K)}(Y,Z^\diamond),
	\]
	where $Y^\sharp$ is the unique untilt (as a perfectoid space over $\Spa(K)$) of $Y$ associated to the structure map $Y\ra \Spd(K)$ (Example \ref{until}).
	The bijection implies that as a $v$-sheaf over $\Pf_v$, the small $v$-sheaf $\hat X^{\infty, \diamond}$ associated to the adic space $\hat X^\infty$ is the pullback of the representable $v$-sheaf $\hat \TT_n^{\infty,\flat}$ along the map $X^\diamond\ra \TT_n^\diamond$. 
	On the other hand, given a perfectoid space $Y$ over $\Spd(K)$ together with the following commutative map
	\[
	\xymatrix{ Y^\sharp\ar[r] \ar[d]& \hat \TT^{\infty}_n\ar[d] \\
		X \ar[r]& \TT_n,}
	\]
	since $X\ra \TT_n$ is finite surjective of the same dimensions, the Theorem \ref{apt} implies that there exits a unique map of adic spaces $Y^\sharp\ra X^\infty_{\pfd}=\Spa(B_{\pfd}[\frac{1}{p}],B_{\pfd})$ that fits into the commutative diagram:
	\[
	\xymatrix{ Y^\sharp\ar@{.>}[rd] \ar[rrd]\ar[rdd] &&\\
		& X^\infty_{\pfd} \ar[r] \ar[d]& \hat \TT^{\infty}_n\ar[d] \\
		&X \ar[r]& \TT_n.}
	\]
	Compare the pullback $\hat X^{\infty,\diamond}$ with the universal affinoid perfectoid space $X_{\pfd}^\infty$, we see the $v$-sheaf $\hat X^{\infty,\diamond}$ is isomorphic to the representable $v$-sheaf $X_{\pfd}^{\infty, \flat}$, where the latter is given by the tilt of the perfectoid space $X_{\pfd}^\infty$.
	In particular, we get the equality 
	\[
	R\Gamma_v(\hat X^{\infty, \diamond},\wh\calO_v)=R\Gamma_v(X_{\pfd}^{\infty,\flat},\wh\calO_v).
	\]
	Since the higher $v$ (pro-\'etale) cohomology of the completed structure sheaf on affinoid perfectoid space vanishes, by combining equalities above we get
	\[
	R\Gamma(X_\pe ,\wh\calO)=R\Gamma_\ct(\ZZ_p(1)^n, B_{\pfd}[\frac{1}{p}]).
	\]
	At last we note that the above object lives in the cohomological degree $[0,n]$ in the derived category of abelian groups, for the continuous group cohomology of $\ZZ_p(1)^n$ can be computed by the Koszul complex of length $n$ (\cite[Section 7]{BMS}).
	Thus we are done.
\end{proof}

\begin{remark}
	Here we want to mention that the cohomological bound given here is stronger than the one in Corollary \ref{coh-bound2}.
\end{remark}

\begin{remark}
	In the proof above, the continuous group cohomology computing $R\Gamma(X_\pe, \wh\calO_X)$ can be defined concretely as $\left( R\varprojlim_m R\Gamma_{\mathrm{disc}}(\ZZ^n,R\Gamma(X_\pe^\infty, \wh\calO_X^+/p^m)) \right) [\frac{1}{p}]$, where $R\Gamma_{disc}(\ZZ^n, -)$ denotes the discrete group cohomology of $\ZZ^n$.
	This is because as $X$ is affinoid (thus quasi-compact and quasi-separated), we have $R\Gamma(X_\pe, \wh\calO_X)=\left(R\varprojlim_m R\Gamma(X_\pe, \wh\calO_X^+/p^m)\right)[\frac{1}{p}]$.
	Moreover, to compute each $R\Gamma(X_\pe, \wh\calO_X^+/p^m)$ we could use the \v{C}ech complex of $\wh\calO_X^+/p^m$ for the pro-\'etale covering  $X^\infty \ra X$.
	We at last note that as the covering is an $\ZZ_p(1)^n$-torsor,
 the \v{C}ech complex is equivalent to the  discrete group cohomology $R\Gamma_{\mathrm{disc}}(\ZZ^n,R\Gamma(X_\pe^\infty, \wh\calO_X^+/p^m))$, by the isomorphism in \cite[Lemma 7.3]{BMS} for $\Gamma=\mathbb{Z}_p(1)^n$.
\end{remark}

\begin{definition}\label{loc comp}
	Let $X$ be a rigid space over $K$.
	We say \emph{$X$ is locally compactifiable} if there exists an open covering $\{U_i\ra X\}_i$ of $X$, such that each $U_i$ admits an open immersion into a proper rigid space $Y_i$ over $K$.
\end{definition}
By definition, any proper rigid space over $K$ is locally compactifiable.
Moreover, by Nagata's compactification in algebraic geometry, any finite type scheme over $K$ admits an open immersion in a proper scheme over $K$.
So the analytification of any finite type scheme over $K$ is a locally compactifiable rigid space.

\begin{proposition}\label{finite impro}
	Let $X$ be a locally compactifiable rigid space over $K$.
	Then the higher direct image $R^i\pi_{X *}\Omega_\eh^j$ vanishes when $i+j> \dim(X)$.
\end{proposition}

\begin{proof}
	Since the vanishing of the higher direct image is a local statement, by taking an open covering, it suffices to assume $X$ admits an open immersion $f:X\ra X'$ for $X'$ being proper over $K$.
	Moreover, by dropping the irreducible components of $X'$ that have higher dimensions, we may assume $\dim(X')$ is the same as $\dim(X)$. 
	This is allowed as the dimension of an irreducible rigid space is not changed when we pass to its open subsets (see the discussion before 2.2.3 in\cite{Con03}).
	
	We then notice that the result is true for $X'$: by the Proposition \ref{finite rev}, we know $R^n\nu_{X' *}\wh\calO_{X'}$ vanishes for $n> \dim(X')$.
	On the other hand, by the degeneracy in the Corollary \ref{direct sum}, each $R^i\pi_{X' *}\Omega_{ \eh}^j(-j)$ is a direct summand of $R^{i+j}\nu_{X' *}\wh\calO_{X'}$.
	This implies that when $i+j>\dim(X')$, the cohomology sheaf $R^i\pi_{X' *}\Omega_{ \eh}^j$ vanishes.
	
	Finally, note that by the coherence proved in Section \ref{sec6}, since $R^i\pi_{X *} \Omega^j_\eh$ is the sheaf associated to the presheaf $U\mapsto \cH^i(U_\eh, \Omega^j_\eh)$ for open subsets $U$ inside of $X$, the preimage of $R^i\pi_{X' *} \Omega^j_\eh$ along $f$ is exactly $R^i\pi_{X *}\Omega^j_\eh$.
	In this way, by the equality of dimensions $\dim(X)=\dim(X')$, the vanishing of $R^i\pi_{X' *}\Omega_{ \eh}^j$ for $i+j>\dim(X')$ implies the vanishing of $R^i\pi_{X *}\Omega_{ \eh}^j$ for $i+j>\dim(X)$.
	So we get the result.
\end{proof}

\section{(Pro-\'etale)-\'eh de Rham comparison}\label{sec8}
In this section, we give a comparison theorem between the pro-\'etale cohomology and the $\eh$ de Rham cohomology, for proper rigid spaces that are defined over a discretely valued subfield.
The idea is to use $v$-hyperdescent for the de Rham period sheaf $\mathbb{B}_\dR^+$ and the simplicial resolution.

Throughout the section, we fix a complete algebraically closed $p$-adic field $K$ of rank one, and a discretely valued subfield $K_0$ that has a perfect residue field, whose associated ring of Witt vectors is denoted as $W$.

\paragraph{Period sheaves for simplicial spaces}
We first recall the period sheaves in the pro-\'etale topos and the $v$-topos, following mostly Section 6 in \cite{Sch13} and the \cite{Sch16}.
Here following the notations in Section \ref{sec7}, we extend construction to a simplicial rigid spaces directly.

Let $X_\bullet$ be a non-augmented simplicial rigid space over $K_0$, truncated or not.
Following the discussion in Subsection \ref{sec7.3}, we consider the commutative diagram of topoi of simplicial sites as below, as a simplicial version of the diagram in Section \ref{sec4}:
\[
\xymatrix{
	\Sh(X_{\bullet \PE}) \ar[r]^{\nu_\bullet} & \Sh(X_{\bullet\et}) \\
	\Sh(\Pf_v|_{X_\bullet^\diamond}) \ar[r]_{\alpha_\bullet} \ar[u]^{\lambda_\bullet} & \Sh(X_{\bullet\eh}) \ar[u]_{\pi_\bullet}.}
\]

We first extend the construction of various period sheaves in \cite{Sch13} and \cite{DLLZ} to the simplicial pro-\'etale site $X_{\bullet \PE}$ as follows.
\begin{definition}\label{Bdr pe}
	Let $X_\bullet$ be a non-augmented simplicial rigid space over $K_0$, truncated or not.
	Consider the following sheaves on $X_{\bullet \PE}$:
	\begin{enumerate}[(i)]
		\item The sheaf $\A_\iif:=W(\wh\calO_\pe^{\flat, +})$, together with a canonical specializing map $\theta:\A_\iif\ra \wh\calO_\pe^+$ extending that of $\Ainf\ra \calO_K$.
		\item The \emph{positive de Rham sheaf}
		\[
		\BB_\dR^+:=\varprojlim \A_\iif[\frac{1}{p}]/(\ker(\theta))^n.
		\]
		with a filtration defined by $\Fil^i\BB_\dR^+=(\ker(\theta))^i\BB_\dR^+$.
		\item The \emph{de Rham sheaf}
		\[
		\BB_\dR:=\BB_\dR^+[\frac{1}{t}],
		\]
		where $t$ is any generator of $\Fil^1\BB_\dR^+$. We equip the de Rham sheaf with the filtration
		\[
		\Fil^i\BB_\dR=\sum_{t\in \ZZ} t^{-j}\Fil^{i+j}\BB_\dR^+.
		\]
	\end{enumerate}
\end{definition}
Note that in the case when $X_\bullet$ is the truncated simplicial rigid space over $\Delta_{\leq 0}$, this recovers the non-simplicial version of the de Rham period sheaf for a rigid space $X=X_0$ as in \cite{Sch13}.
Moreover, the positive de Rham sheaf $\BB_\dR^+$ is filtered complete, with its $i$-th graded factor $\gr^i\BB_\dR^+$ equals to $t^i\wh\calO_{X_\PE}=\wh\calO_{X_\PE}(i)$ (the non-simplicial version is in \cite[Proposition 6.7]{Sch13}).

To associated a geometric structure on $\BB_\dR$ (namely the differential operator), we recall the de Rham period sheaf as follows.
\begin{definition}\label{OBdr pe}
	Let $X_\bullet$ be a simplicial smooth rigid space over $K_0$.
	Consider the following sheaves on $X_{\bullet \PE}$:
	\begin{enumerate}[(i)]
		\item The \emph{positive structure de Rham sheaf} $\OBdrp$ is defined as the sheaf associated to the presheaf sending $U\in X_{n, \pe}\subset X_{\bullet \pe}$ to the direct limit of the following rings
		\[
		\ker(\theta)-completion~of~\left(\left(\calO_\et^+(U_j)\wh\otimes_{W} \A_\iif (U)\right)[\frac{1}{p}]\right),
		\]
		where $\{U_j\}$ is a pro-\'etale presentation of $U$ as in \cite{Sch13}.
		It has a filtration given by $\Fil^i\OBdrp=(\ker(\theta))^i\OBdrp$.
		\item 
		Let $\calF$ be the sheaf $\OBdrp[\frac{1}{t}]$, with the filtration 
		\[
		\Fil^i\calF=\sum_{j\in \ZZ} t^{-j}\Fil^{i+j}\OBdrp,
		\]
		where $t$ is a generator of $\Fil^1\BB_\dR^+$.
		The \emph{structure de Rham sheaf} $\OBdr$ is then defined as the filtered completion of $\calF$, namely
		\[
		\calO\BB_\dR:=\varprojlim \calF/\Fil^i\calF. 
		\]
		Here $\OBdr$ comes with a natural filtration from $\calF$, whose $i$-th graded factors are isomorphic to $\gr^i\calF$.
		
	\end{enumerate}
\end{definition}
Here we note that slightly different from \cite{Sch16}, we need to apply a filtered completion at $\calF$ to get the structure de Rham sheaf $\OBdr$.
This is because the sheaf $\calF$ is not complete under the filtration (See Remark 3.11 in \cite{DLLZ}).
We also notice that as explained in \cite{Sch13}, locally on $X_\pe$ the element $t$ exists, is a non zero divisor, and is unique up to units. 
So the above definitions are well-defined.

\paragraph{Comparisons}
Let us assume $X_\bullet$ is a simplicial smooth rigid space.
The sheaf $\OBdrp$ over $X_{\bullet\pe}$ admits a canonical $\BB_\dR^+$-linear connection $\nabla$ induced from the differential map of $\calO_{X_\bullet}$ over $X_{\bullet \et}$, with the following diagram commutes:
\[
\xymatrix{
	\OBdrp\ar[r]^{\nabla~~~~~~~~~~} & \OBdrp\otimes_{\nu_\bullet^{-1}\calO_{X_\bullet}} \nu_\bullet^{-1}\Omega_{X_\bullet/K_0}^1\\
	\nu_\bullet^{-1}\calO_{X_\bullet}\ar[r]^d \ar[u]& \nu_\bullet^{-1}\Omega_{X_\bullet/K_0}^1\ar[u],}
\]
which is functorial among simplicial smooth rigid spaces $X_\bullet$ over $K_0$. 
Besides, the above allows us to give a natural tensor product filtration on the sequence $\OBdrp\otimes_{\nu_\bullet^{-1}\calO_{X_\bullet}}\nu_\bullet^{-1}\Omega_{X_\bullet/K_0}^\bullet$ by assigning the following subsequences:
\[
\Fil^i(\OBdrp\otimes_{\nu_\bullet^{-1}\calO_{X_\bullet}}\nu_\bullet^{-1}\Omega_{X_\bullet/K_0}^\bullet)=\sum_{j\in \ZZ} \Fil^j(\OBdrp)\otimes_{\nu_\bullet^{-1}\calO_{X_\bullet}}\nu_\bullet^{-1}\Omega_{X_\bullet/K_0}^{\geq i-j},
\]
where $\Omega_{X_\bullet/K_0}^{\geq i-j}$ is the $(i-j)$-th Hodge filtration defined by the naive truncation of the de Rham complex.
This filtration is compatible with the Hodge filtration on the de Rham complex, in the sense that the connection above induces the natural map of subcomplexes:
\[
\nu_\bullet^{-1}\Fil^i\Omega_{X_\bullet/K_0}^\bullet\rra \Fil^i(\OBdrp\otimes_{\nu_\bullet^{-1}\calO_{X_\bullet}}\nu_\bullet^{-1}\Omega_{X_\bullet/K_0}^\bullet).
\]

Furthermore, replace the complex $\OBdrp\otimes_{\nu_\bullet^{-1}\calO_{X_\bullet}}\nu_\bullet^{-1}\Omega_{X_\bullet/K_0}^\bullet$ by $\OBdr\otimes_{\nu_\bullet^{-1}\calO_{X_\bullet}}\nu_\bullet^{-1}\Omega_{X_\bullet/K_0}^\bullet$ in the above, we get the natural tensor product filtration compatible with the former:
\[
\Fil^i(\OBdr\otimes_{\nu_\bullet^{-1}\calO_{X_\bullet}}\nu_\bullet^{-1}\Omega_{X_\bullet/K_0}^\bullet)=\sum_{j\in \ZZ} \Fil^j(\OBdr)\otimes_{\nu_\bullet^{-1}\calO_{X_\bullet}}\nu_\bullet^{-1}\Omega_{X_\bullet/K_0}^{\geq i-j}.
\]

Now let us recall the following Poincar\'e Lemma on the small pro-\'etale site $X_{\bullet\pe}$:
\begin{lemma}[Poincar\'e Lemma.~\cite{Sch13} 6.13; \cite{DLLZ} 3.38]\label{Poin Bdr}
	Let $X_\bullet$ be a simplicial smooth rigid space over $K_0$. 
	Then the following natural sequence is an acyclic filtered complex of $\BB_\dR^+$-linear sheaves over $X_{\bullet\pe}$:
	\[
	\xymatrix{ 0\ar[r] &\BB_\dR^+\ar[r]& \OBdrp\ar[r]^{\nabla~~~~~~~~~~~} \ar[r]& \OBdrp\otimes_{\nu_\bullet^{-1}\calO_{X_\bullet}} \nu_\bullet^{-1}\Omega_{X_\bullet}^1 \ar[r]^{~~~~~~~~~~~~~\nabla} &\cdots\ar[r]^{\nabla~~~~~~~~~~~~~} &\OBdrp\otimes_{\nu_\bullet^{-1}\calO_{X_\bullet}} \nu_\bullet^{-1}\Omega_{X_\bullet}^n \ar[r]& \cdots .}
	\]
	Moreover, the same holds when we replace $\BB_\dR^+$ and $\OBdrp$ by $\BB_\dR$ and $\OBdr$ separately (together with their filtrations). 
\end{lemma}
\begin{proof}
	The non-simplicial version of the result is proved in \cite[6.13]{Sch13} and \cite[3.38]{DLLZ}.
	In general, the acyclicity of the natural sequence is checked by applying the restriction functor, as in Lemma \ref{sim-van}.
\end{proof}
We then recall the (pro-\'etale)-de Rham comparison for proper smooth rigid spaces in \cite{Sch13}.

Let $X_\bullet$ be a simplicial smooth rigid space over $K_0$.
As discussed above, we can form the following natural map in the filtered derived category of abelian sheaves over $X_{\bullet \et}$:
\[
\Omega_{X_\bullet/K_0}^\bullet\rra R\nu_{\bullet*}(\OBdr\otimes_{\nu_\bullet^{-1}\calO_{X_\bullet}} \nu_\bullet^{-1}\Omega_{X_\bullet/K_0}^\bullet).
\]
By taking the derived global section, we get a filtered morphism
\[
R\Gamma(X_{\bullet},\Omega_{X_\bullet/K_0}^\bullet)\rra R\Gamma(X_{\bullet K, \pe},\OBdr\otimes_{\nu_\bullet^{-1}\calO_{X_\bullet}} \nu_\bullet^{-1}\Omega_{X_\bullet/K_0}^\bullet).
\]
As the right side is $\mathrm{B_{dR}}$-linear and the map above is compatible with the filtration on $\mathrm{B_{dR}}$, the above induces the following $\mathrm{B_{dR}}$-linear map
\[
R\Gamma(X_{\bullet},\Omega_{X_\bullet/K_0}^\bullet)\otimes_{K_0}\mathrm{B_{dR}}\rra R\Gamma(X_{\bullet K,\pe},\OBdr\otimes_{\nu_\bullet^{-1}\calO_{X_\bullet}} \nu_\bullet^{-1}\Omega_{X_\bullet/K_0}^\bullet).
\]
Moreover, by endowing the complex $R\Gamma(X_{\bullet},\Omega_{X_\bullet/K_0}^\bullet)\otimes_{K_0}\mathrm{B_{dR}}$ with the (derived) tensor product filtration, the above is in fact a filtered map (see \cite[Chapter 5]{Ill71}).
The Poincar\'e Lemma \ref{Poin Bdr} implies that the right side is filtered quasi-isomorphic to $R\nu_{\bullet*}\BB_\dR$. 
So we get the following filtered morphism
\[
R\Gamma(X_{\bullet},\Omega_{X_\bullet/K_0}^\bullet)\otimes_{K_0}\mathrm{B_{dR}}\rra R\Gamma(X_{\bullet K, \pe},\BB_\dR).
\]
Here we also note that the map is $\Gal(K/K_0)$-equivariant.
This is in fact a quasi-isomorphism, assuming the proper condition:
\begin{theorem}[\cite{Sch13}, Theorem 7.11; \cite{DLLZ} 3.49]\label{et dR-proet}
	Let $X_\bullet$ be a simplicial proper smooth rigid space over $K_0$.
	Then the following two natural maps are $\Gal(K/K_0)$-equivariant filtered quasi-isomorphisms:

	\[
	\xymatrix{ R\Gamma(X_{\bullet K,\pe},\BB_\dR) \ar[r] &R\Gamma(X_{\bullet K,\pe},\OBdr\otimes_{\nu_\bullet^{-1}\calO_{X_\bullet}} \nu^{-1}\Omega_{X_\bullet/K_0}^\bullet) & R\Gamma(X_{\bullet},\Omega_{X_\bullet/K_0}^\bullet)\otimes_{K_0}\mathrm{B_\dR} \ar[l].}
	\]

\end{theorem}
\begin{proof}
	The non-simplicial version was checked in \cite[Theorem 7.11]{Sch13} and \cite[3.49]{DLLZ}.
	Namely, for each proper smooth rigid space $X_n$ over $K_0$, the natural maps below are filtered quasi-isomorphisms:
	\[
	\xymatrix{ R\Gamma(X_{n K,\pe},\BB_\dR) \ar[r] &R\Gamma(X_{n K,\pe},\OBdr\otimes_{\nu_n^{-1}\calO_{X_n}} \nu^{-1}\Omega_{X_n/K_0}^\bullet) & R\Gamma(X_{n},\Omega_{X_n/K_0}^\bullet)\otimes_{K_0}\mathrm{B_\dR} \ar[l].}
	\]
	To get the simplicial version as in the statement, we take the homotopy limit over the simplicial diagram $\Delta$.
	The homotopy limits of the left and the middle terms above are exactly the left and middle terms as in the statement, so it suffices to check that the following natural map is a filtered quasi-isomorphism
	\[
	R\Gamma(X_{\bullet},\Omega_{X_\bullet/K_0}^\bullet)\otimes_{K_0}\mathrm{B_\dR} \rra R\varprojlim_{[n]\in \Delta} R\Gamma(X_{n},\Omega_{X_n/K_0}^\bullet)\otimes_{K_0}\mathrm{B_\dR}.
	\]
	This can be checked by looking at the graded pieces for $\mathrm{B_{dR}}$ together with the natural filtered quasi-isomorphism
	\[
	R\varprojlim_{[n]\in \Delta} R\Gamma(X_{n},\Omega_{X_n/K_0}^\bullet) \cong R\Gamma(X_{\bullet},\Omega_{X_\bullet/K_0}^\bullet).
	\]
	So we are done for the simplicial case.
\end{proof}
\begin{remark}
Here we remark that the above quasi-isomorphisms between cohomology of simplicial sites could also follow from Lemma \ref{ss}:
it suffices to show the cone of a map vanishes, which follows from the spectral sequence in Lemma \ref{ss} and the vanishing of the cone for each individual $X_n$.
\end{remark}

Moreover, we can in fact replace the de Rham complex to the $\eh$ de Rham complex to compute the cohomology:
\begin{proposition}\label{eh dR-proet sm}
	Let $X_\bullet$ be a proper smooth rigid space over $K_0$.
	Then the map of simplicial ringed sites $\pi_\bullet:X_{\bullet\eh}\ra X_{\bullet\et}$ induces a canonical $
	\Gal(K/K_0)$-equivariant filtered quasi-isomorphism
	\[
	   R\Gamma(X_{\bullet\eh},\Omega_{X_\bullet/K_0,\eh}^\bullet)\otimes_{K_0} \mathrm{B_{dR}}\rra R\Gamma(X_{\bullet K,\pe},\BB_\dR). 
	\]
	Here the filtration on the complex $R\Gamma(X_{\bullet\eh},\Omega_{X_\bullet/K_0,\eh}^\bullet)\otimes_{K_0} \mathrm{B_{dR}}$ is defined by the derived tensor product filtration
	\[
	\Fil^i(R\Gamma(X_{\bullet\eh},\Omega_{X_\bullet/K_0,\eh}^\bullet)\otimes_{K_0} \mathrm{B_{dR}})=\varinjlim_{j\in \ZZ} R\Gamma(X_{X_\bullet\eh},\Omega_{X_\bullet/K_0,\eh}^{\geq i-j})\otimes_{K_0}\Fil^j\mathrm{B_{dR}}.
	\]

\end{proposition}
\begin{proof}

	By the Theorem \ref{et dR-proet} and the definition of the filtrations, to prove the map is a filtered quasi-isomorphism, it suffices to prove the filtered quasi-isomorphism for the following map on the simplicial \'etale site $X_{\bullet\et}$:
	\[
	\Omega_{X_\bullet/K_0}^\bullet\rra R\pi_{X_\bullet *}\Omega_{X_\bullet/K_0,\eh}^\bullet.
	\]
	Then by the spectral sequence for the filtration given by the naive truncation of the $\eh$ de Rham complex, it suffices to show the quasi-isomorphism for $\Omega_{X_\bullet/K_0}^i\ra R\pi_{X_\bullet *}\Omega_{X_\bullet/K_0,\eh}^i$ for $i\in \NN$ and smooth $X_\bullet/K_0$, which follows from the $\eh$ descent for differential forms (Theorem \ref{descent}) and the vanishing criterion by restrictions (Lemma \ref{sim-van}).
	
\end{proof}

Now we are ready to prove the (pro-\'etale)-de Rham comparison, for non-smooth proper rigid spaces.
\begin{theorem}\label{eh dR-proet}
	Let $X$ be a proper rigid space over $K_0$.
	Then there exists a $\Gal(K/K_0)$-equivariant filtered quasi-isomorphism
	\[
	R\Gamma(X_\eh,\Omega_{\eh, /K_0}^\bullet)\otimes_{K_0} \mathrm{B_{dR}}\rra R\Gamma(X_{K \pe}, \BB_\dR),
	\]
	which generalizes the smooth case in the Proposition \ref{eh dR-proet sm}.	
\end{theorem}
\begin{proof}
	We first notice that the pro-\'etale cohomology of $\BB_\dR$ satisfies the $\eh$-hyperdescent.
	More precisely, let $\rho:X_\bullet \ra X$ be  a $\eh$-hypercovering such that each $X_n$ is built from blowing-ups and is smooth and proper over $K_0$.
	Then we claim that the natural filtered map below is a filtered quasi-isomorphism
	\[
	R\Gamma(X_{K \pe},\BB_\dR) \rra R\Gamma(X_{\bullet K, \pe}, \BB_\dR).
	\]
	To see this, let us first notice that since $\BB_\dR$ is filtered complete, it suffices to check that the quasi-isomorphism for each graded pieces as below
	\[
	R\Gamma(X_{K \pe},\wh\calO_X(i)) \rra R\Gamma(X_{\bullet K, \pe}, \wh\calO_{X_\bullet}(i)).
	\]
	This then follows from the $\eh$-hyperdescent ($v$-hyperdescent) of the cohomology of pro-\'etale structure sheaf in Proposition \ref{pe-v}. 
	On the  other hand, similar to the proof of Theorem \ref{et dR-proet}, the natural map of sites $X_\bullet\ra X$ induces the following filtered quasi-isomorphism
	\[
	R\Gamma(X_\eh,\Omega_{\eh, /K_0}^\bullet)\otimes_{K_0} \mathrm{B_{dR}} \rra R\Gamma(X_{\bullet\eh},\Omega_{X_\bullet/K_0, \eh}^\bullet)\otimes_{K_0} \mathrm{B_{dR}}.
	\]
	So the rest follows from the natural filtered quasi-isomorphism as in Proposition \ref{eh dR-proet sm}
	\[
	R\Gamma(X_{\bullet\eh},\Omega_{X_\bullet/K_0,\eh}^\bullet)\otimes_{K_0} \mathrm{B_{dR}}\rra R\Gamma(X_{\bullet K,\pe},\BB_\dR).
	\]
	Hence we are done.

\end{proof}
As the quasi-isomorphisms above are all filtered with respect to natural filtrations, by taking the $0$-th graded piece, we get back to the Hodge-Tate decomposition for proper rigid spaces over a discretely valued field.

Moreover, we can use the above to prove the degeneracy of the $\eh$ version of the Hodge-de Rham spectral sequence.
Recall that the naive truncation of $\Omega_\eh^\bullet$ gives a filtration on it, whose associated spectral sequence is 
\[
E_1^{p,q}=\cH^q(X_\eh,\Omega_{ \eh, /K_0}^p)\Longrightarrow\cH^{p+q}(X_\eh,\Omega_{ \eh, /K_0}^\bullet).
\]
This is the $\eh$ version of the Hodge-de Rham spectral sequence.
\begin{proposition}\label{HT-dR deg}
	Let $X$ be a proper rigid space over $K_0$.
	Then the $\eh$ Hodge-de Rham spectral sequence degenerates at its $E_1$-page.
\end{proposition}
\begin{proof}
	As each $\cH^q(X_\eh,\Omega_{\eh, /K_0}^p)$ is of finite dimension over $K_0$ (by the properness of $X$ and Proposition \ref{coh}), it suffices to show that 
	\[
	\sum_{p+q=n}\dim_{K_0}	\cH^q(X_\eh,\Omega_{\eh, /K_0}^p)=\dim_{K_0}\cH^n(X_\eh,\Omega_{\eh,/K_0}^\bullet).
	\]
	We first take the tensor product of the $K_0$-vector space $\cH^n(X_\eh,\Omega_{\eh,/K_0}^\bullet)$ with the field $\mathrm{B_{dR}}$.
	Then the comparison theorem \ref{eh dR-proet} implies that the right side above is equal to 
	\[
	\dim_{\mathrm{B_{dR}}}\cH^n(X_{K \pe},\BB_\dR).
	\]
	Note that $\cH^n(X_{K \pe},\BB_\dR)$ is a $\mathrm{B_{dR}}$-vector space of finite dimension that has a filtration compatible with that of $\mathrm{B_{dR}}$, induced by the image of $\rmH^n(X_{K \pe}, \BB_\dR^+)$ in $\rmH^n(X_{K \pe}, \BB_\dR)$.
	By the Primitive Comparison Theorem in \cite{Sch12B}, we have
	\[
	\rmH^n(X_{K \pe}, \BB_\dR^+) \cong \rmH^n(X_{K \et},\QQ_p)\otimes_{\QQ_p} \Bdr.
	\]
	In particular, the cohomology group $\rmH^n(X_{K \pe}, \BB_\dR^+)$ is $\xi$-torsion sheaf and its map to $\rmH^n(X_{K \pe}, \BB_\dR)$ is injective.
	Thus $\rmH^n(X_{K \pe}, \BB_\dR)$ is a finite dimensional $\mathrm{B_{dR}}$-vector space whose dimension over $\mathrm{B_{dR}}$ is equal to $\rank_\Bdr \rmH^n(X_{K \pe}, \BB_\dR^+)$.
	In particular, the $0$-th graded piece is $\gr^0\mathrm{B_{dR}}(=K)$-vector space $\cH^n(X_{K \pe},\wh\calO_X)$ whose $K$-dimension is equal to $\dim_{\mathrm{B_{dR}}}\rmH^n(X_{K \pe}, \BB_\dR)$.
	In other words, we get the equality:
	\[
	\dim_{\mathrm{B_{dR}}}\cH^n(X_{K \pe},\BB_\dR)=\dim_K\cH^n(X_{K \pe},\wh\calO_X).
	\]
	In this way, by the degeneracy theorem for the derived direct image $R\nu_*\wh\calO_X$ (Corollary \ref{direct sum}), we get
	\begin{align*}
	\dim_K\cH^n(X_{K \pe},\wh\calO_X)&=\sum_{p+q=n}\dim_K\cH^p(X_{K \eh},\Omega_{ \eh, /K}^p(-p))\\
	&=\sum_{p+q=n}\dim_K\cH^p(X_{\eh},\Omega_{ \eh, /K_0}^p),
	\end{align*}
	where the last equality follows from the coherence of $R\pi_*\Omega_{ \eh, /K_0}^p$ on the small \'etale site $X_\et$ (Proposition \ref{coh}).
\end{proof}

\begin{remark}
	Pointed out by David Hansen,  the (pro-\'etale)-$\eh$ de Rham comparison could be  extended to general de Rham local systems (in the sense of \cite{Sch13}) for non-smooth proper rigid spaces, not just the trivial local system. 
	As the de Rham complex for local systems over a non-smooth rigid space does not behave very well, we need to consider locally free $\BB_\dR$ sheaves over the $v$-site, instead of over pro-\'etale site.
	Then a comparison of $v$-cohomology and ($\eh$) de Rham cohomology will follow from the analogous hyperdescent argument as in Theorem \ref{eh dR-proet}.
	
\end{remark}

\section{Hodge-Tate decomposition for non-smooth spaces}\label{sec9}
At last, we give the application of our results to the Hodge-Tate decomposition for non-smooth spaces, as mentioned in the introduction.
Throughout the section, let $X$ be a proper rigid space over a complete algebraically closed non-archimedean field $K$ over $\Q_p$.

Recall that by the Primitive Comparison (\cite[Theorem 3.17]{Sch12B}), we have
\[
\cH^n(X_\et,\Q_p)\otimes_{\Q_p} K=\cH^n(X_{\pe},\wh{\calO}_X).
\]
The equality enables us to compute the $p$-adic \'etale cohomology by studying the pro-\'etale cohomology.
In particular, by taking the associated derived version, the right side above can be obtained by
\begin{align*}
R\Gamma(X_\PE, \wh\calO_X)&=R\Gamma(X_\et,R\nu_*\wh\calO_X).
\end{align*}

Then we recall the following diagram of topoi associated to $X$ in Section \ref{sec4}:
\[
\xymatrix{
	\Sh(X_\PE) \ar[r]^\nu &  \Sh(X_\et)\\
	\Sh(\Pf_v|_{X^\diamond}) \ar[r]_\alpha \ar[u]^\lambda& \Sh(X_\eh) \ar[u]_{\pi_X}.}
\]
The (pro-\'etale)-v comparison (see the Proposition \ref{pe-v}) allows us to replace $R\nu_*\wh\calO_X$ by the derived direct image  $R\pi_{X *}R\alpha_*\wh\calO_v$ of the untilted complete $v$-structure sheaf.
So we have
\begin{align*}
R\Gamma(X_\PE, \wh\calO_X)&=R\Gamma(X_\et,R\pi_{X *}R\alpha_*\wh\calO_v)\\
&=R\Gamma(X_\eh,R\alpha_*\wh\calO_v).
\end{align*}
By the discussion in Section \ref{sec4}, we have
\[
R^j\alpha_*\wh\calO_v=\Omega^j_\eh(-j).
\]
So by replacing the above equality into the Leray spectral sequence for the composition of derived functors, we get
\[
E_2^{i,j}=\cH^i(X_\eh,\Omega_{ \eh}^j)(-j)\Longrightarrow \cH^{i+j}(X_\pe,\wh\calO_X).
\]
This together with the Primitive Comparison leads to the \emph{Hodge-Tate spectral sequence} for proper rigid space $X$
\[
E_2^{i,j}=\cH^i(X_\eh,\Omega_{ \eh}^j)(-j)\Longrightarrow \cH^{i+j}(X_\et,\Q_p)\otimes_{\Q_p}K.
\]
The name is justified by the special case of the $\eh$ differential in the Theorem \ref{descent}: when $X$ is smooth, the higher direct image of the $\eh$-differential vanishes, and the spectral sequence degenerates into
\[
E_2^{i,j}=\cH^i(X,\Omega_{X/K}^j)(-j)\Longrightarrow \cH^{i+j}(X_\et,\Q_p)\otimes_{\Q_p} K,
\]
with each $\cH^i(X_\eh,\Omega_{ \eh}^j)$ identified with $\cH^i(X,\Omega_{X/K}^j)$.

Now by the strong liftability of $X$ (Proposition \ref{proper lift}) and the Degeneracy Theorem \ref{deg}, the derived direct image $R\nu_*\wh\calO_X$ is non-canonically quasi-isomorphic to the direct sum \[\bigoplus_{i=0}^{\dim(X)} R\pi_{X *}(\Omega_{ \eh}^j(-j)[-j]).\]
Replace the $R\nu_*\wh\calO_X$ by this direct sum, we have
\begin{align*}
R\Gamma(X_\pe, \wh\calO_X)&=\bigoplus_{i=0}^{\dim(X)} R\Gamma(X_\eh, \Omega_{ \eh}^j(-j)[-j]).
\end{align*}
So after taking the $n$-th cohomology, we see the Hodge-Tate spectral sequence degenerates at its $E_2$-page. 
\begin{theorem}[Hodge-Tate decomposition]
	Let $X$ be a proper rigid space over a complete algebraically closed non-archimedean field $K$ of characteristic $0$.
	Then there exists a natural spectral sequence to its $p$-adic \'etale cohomology
	\[
	E_2^{i,j}=\cH^i(X_\eh,\Omega_{ \eh}^j)(-j)\Longrightarrow \cH^{i+j}(X_\et,\Q_p)\otimes_{\Q_p}K.
	\]
	Here the spectral sequence degenerates at its $E_2$-page, and $\cH^i(X_\eh, \Omega_{ \eh}^j)(-j)$ is a finite dimensional $K$-vector space that vanishes unless $0\leq i,j\leq n$.
	
	When $X$ is a smooth rigid space, $\cH^i(X_\eh, \Omega_{ \eh}^j)(-j)$ is isomorphic to $\cH^i(X, \Omega_{X/K}^j)(-j)$, and the spectral sequence is the same as the Hodge-Tate spectral sequence for smooth proper rigid space (in the sense of \cite{Sch12B}).
	
\end{theorem}
\begin{proof}
	The cohomological boundedness of $\cH^i(X_\eh, \Omega_{ \eh}^j)(-j)$ is given by the Theorem \ref{coh-bound0}.
	The finite dimensionality is given by the properness of $X$, the coherence of the $R^i\pi_{X *}\Omega_{ \eh}^j$ (the Proposition \ref{coh}), and the following equality 
	\[
	R\Gamma(X_\eh, \Omega_\eh^j)=R\Gamma(X_\et, R\pi_{X *}\Omega_\eh^j).
	\]
	Moreover when $X$ is smooth, the isomorphism between $\cH^i(X_\eh, \Omega_{ \eh}^j)(-j)$ and $\cH^i(X, \Omega_{X/K}^j)(-j)$ follows from the $\eh$-decent of differential by the Theorem \ref{descent}
\end{proof}
At last, when $X$ is defined over a discretely valued subfield $K_0$ of $K$ that has a perfect residue field, the above spectral sequence is Galois equivariant.
In particular, since $K(-j)^{\Gal(K/K_0)}=0$ for $j\neq 0$, the boundary map from $\cH^i(X_\eh,\Omega_{ \eh}^j)(-j)$ to $\cH^{i+2}(X_\eh,\Omega_{ \eh}^{j-1})(-j+1)$ is zero.
In this way, the Hodge-Tate spectral sequence degenerates canonically, and the $p$-adic \'etale cohomology splits into the direct sum of distinct Hodge-Tate weights 
\[
\cH^{n}(X_\et,\Q_p)\otimes_{\Q_p} K=\bigoplus_{i+j=n}\cH^i(X_\eh,\Omega_{ \eh}^j)(-j).
\]
This canonical (Galois equivariant) decomposition is functorial with respect to rigid spaces defined over $K_0$.

\begin{theorem}[Hodge-Tate decomposition]
	Let $Y$ be a proper rigid space over a discretely valued subfield $K_0$ of $K$ that has a perfect residue field.
	Then the spectral sequence above degenerates at its $E_2$-page.
	In fact, we have a Galois equivariant isomorphism 
	\[
	\cH^n(Y_{K\,\et},\Q_p)\otimes_{\Q_p} K=\bigoplus_{i+j=n}\cH^i(Y_\eh,\Omega_{ \eh,/K_0}^j)\otimes_{K_0} K(-j).
	\]
	The isomorphism is functorial with respect to rigid spaces $Y$ over $K_0$.
\end{theorem}

\section*{Acknowledgments}
\addcontentsline{toc}{section}{\protect\numberline{}Acknowledgments}
The author thanks David Hansen for proposing this problem, including the suggestion on the (pro-\'etale)-$\eh$ de Rham comparison in Section \ref{sec8}.
The article is part of the author's Ph.D. thesis.
The author is extremely grateful for uncountable discussions with his advisor Bhargav Bhatt; without his advisor's support, the author would not have been able to finish the article.
In addition, he would like to thank Shizhang Li, Emanuel Reinecke and Matt Stevenson for many helpful discussions.
At last, the author thanks the anonymous referees for providing many comments and suggestions and for catching imprecisions in previous drafts, which helps a lot to improve this article.

\addcontentsline{toc}{section}{\protect\numberline{}References}
\providecommand{\bysame}{\leavevmode\hbox to3em{\hrulefill}\thinspace}
\providecommand{\MR}{\relax\ifhmode\unskip\space\fi MR }
\providecommand{\MRhref}[2]{%
	\href{http://www.ams.org/mathscinet-getitem?mr=#1}{#2}
}
\providecommand{\href}[2]{#2}

\end{document}